\documentclass[12pt,a4paper]{article}
\usepackage[english]{babel}
\usepackage{verbatim}
\usepackage{bm}

\usepackage{amsmath}
\usepackage{amssymb}
\usepackage{amsfonts}
\usepackage{amsbsy}

\usepackage{graphicx}
\usepackage{graphics}
\usepackage{subcaption}
\usepackage{epsfig}
\usepackage{color}
\usepackage{array}
\usepackage{booktabs}
\usepackage{caption}
\usepackage{tabularx}
\usepackage[utf8]{inputenc}
\usepackage{listings} 

\usepackage{newlfont}
\usepackage{bm}
\usepackage{bookmark}
\usepackage{changepage}
\usepackage{geometry}
\usepackage{float}
\usepackage{nomencl}
\usepackage{cases}

\usepackage{siunitx}

\DeclareMathOperator{\dive}{\nabla\cdot}

\DeclareMathOperator{\atanh}{atanh}

\usepackage{todonotes}

\begin{document}

\title{An efficient IMEX-DG solver for the compressible Navier-Stokes equations for  non-ideal gases}

\author{Giuseppe Orlando$^{(1)}$\\ 
 Paolo Francesco Barbante$^{(1)}$, Luca Bonaventura$^{(1)}$}

\date{}

\maketitle

\begin{center}
{\small
$^{(1)}$  
MOX, Dipartimento di Matematica, Politecnico di Milano \\
Piazza Leonardo da Vinci 32, 20133 Milano, Italy\\
{\tt giuseppe.orlando@polimi.it, luca.bonaventura@polimi.it, paolo.barbante@polimi.it}
}
\end{center}


\vspace*{0.5cm}


\vspace*{0.5cm}

\pagebreak

\abstract{We propose an efficient, accurate and robust IMEX solver for the compressible  Navier-Stokes equation describing non-ideal gases with general equations of state. The method, which is based on an $h-$adaptive Discontinuos Galerkin spatial discretization and on an Additive Runge Kutta IMEX method for time discretization, is tailored for low Mach number applications and allows to simulate low Mach regimes at a significantly reduced computational cost, while maintaining full second order accuracy also for higher Mach number regimes. The method has been implemented in the framework of the \textit{deal.II} numerical library, whose adaptive mesh refinement capabilities are employed to enhance efficiency. Refinement indicators appropriate for real gas phenomena have been introduced. A number of numerical experiments on classical benchmarks for compressible flows and their extension to real gases demonstrate the properties of the proposed method.}


\section{Introduction}
\label{sec:intro} \indent

The efficient numerical solution of the compressible Navier-Stokes equations poses several major computational challenges. In particular, for flow regimes characterized by low Mach number and moderate Reynolds number values, severe time step restrictions may be required by standard explicit time discretization methods.
The use of implicit and semi-implicit methods has a long tradition in low Mach number flows, see for example the seminal papers \cite{casulli:1984, cullen:1990, robert:1982}, 
as well as many other contributions in the literature on numerical weather prediction,
see e.g. \cite{bonaventura:2000, giraldo:2005, giraldo:2013, giraldo:2010, kuehnlein:2019, restelli:2009, smolarkiewicz:2019, tumolo:2015} and the reviews in \cite{steppeler:2003, bonaventura:2012}.
Other contributions have been proposed in the literature on more classical computational fluid dynamics, see e.g. \cite{bassi:2007, bassi:2015, boscheri:2021, busto:2020,  dumbser:2016, munz:2003, tavelli:2017}. Many of these contributions focus exclusively on the equations of motion of an ideal gas and their extension to real gases is not necessarily straightforward.
Stability concerns are even more critical in these particular regimes for spatial discretizations based on the Discontinuous Galerkin (DG) method (see e.g.  \cite{giraldo:2020, karniadakis:2005} for a general presentation of this method), which is the spatial discretization used in many of the above referenced papers.

In this work, we seek to combine an accurate and flexible discontinuous DG space discretization with a second order implicit-explicit (IMEX) time discretization, see e.g. \cite{kennedy:2003, pareschi:2005, boscarino:2021}, to obtain an efficient method for compressible flow of real gases at low to moderate Mach numbers. Our goal is to derive a method that can then be easily extended to handle multiphase compressible flows, where a number of coupling and forcing terms arise that cannot be dealt with efficiently by straightforward application of conventional solvers. 
We extend to second order in time the approach of \cite{casulli:1984, dumbser:2016}, coupling implicitly the energy equation to the momentum equation, while treating the continuity equation in an explicit fashion. Our treatment also provides an outline of how a generic IMEX method based on a Diagonally Implicit Runge Kutta can be extended along the same lines. Notice that, with respect to the IMEX approach proposed for the Euler equations in \cite{zeifang:2019}, the technique presented here does not require to introduce reference solutions, does not introduce inconsistencies in the splitting with respect to a reference solution and only requires the solution of linear system of a size equal to that of the number of discrete degrees of freedom needed to describe a scalar variable, as in \cite{dumbser:2016}. 
A conceptually similar approach has been used in \cite{kuehnlein:2019, smolarkiewicz:2019} for the discretization employed in the IFS-FVM atmospheric model. 
In order to obtain a formulation that is efficient also in presence of non negligible viscous terms, we resort to an operator splitting approach, see e.g. \cite{leveque:2002}.
More specifically, as commonly done in numerical models for atmospheric physics, we split the hyperbolic part of the problem, which is treated by an IMEX extension of the method proposed in \cite{dumbser:2016}, from the diffusive terms, which are treated implicitly. Second order accuracy can then be obtained by the Strang splitting approach \cite{leveque:2002, strang:1968}.

For the spatial discretization, we rely on the DG approach implemented in the numerical library \textit{deal.II} \cite{bangerth:2007}, which is a very convenient environment to develop a reliable and easily accessible tool for large scale industrial applications.
This software also provides $h-$refinement capabilities that are exploited by the proposed method. 
For the specific case of real gases, physically based refinement criteria have been proposed and tested, which allow to track accurately convection phenomena also for more general equations of state.
The numerical experiments reported below show the ability of the proposed scheme and of its adaptive implementation to perform accurate simulations in  a range of different settings appropriate to describe non-ideal gas dynamics.

The model equations and their non-dimensional formulation are reviewed in Section \ref{sec:modeleq}. The time discretization approach is outlined and discussed in Section \ref{sec:tdisc}. The spatial discretization is presented in Section \ref{sec:sdisc}. Some implementation issues are described in Section \ref{sec:implement}, while the validation of the proposed method and its application to a number of significant benchmarks is reported in Section \ref{sec:tests}. Some conclusions and perspectives for future work are presented in Section \ref{sec:conclu}.

\section{The compressible Navier-Stokes equations}
\label{sec:modeleq} \indent

Let \(\Omega \subset \mathbb{R}^{d}, 2 \le d \le 3\) be a connected open bounded set with a sufficiently smooth boundary \(\partial\Omega\) and denote by \(\mathbf{x}\) the spatial coordinates and by \(t\) the temporal coordinate. We consider the classical unsteady compressible Navier-Stokes equations, written in flux form as:

\begin{eqnarray}
\label{eq:ns_comp}
\frac{\partial \rho}{\partial t}  
+\dive\left(\rho\mathbf{u} \right) &=& 0   \nonumber \\
\frac{\partial (\rho\mathbf{u})}{\partial t}  
+\dive\left(\rho \mathbf{u} \otimes\mathbf{u}\right) 
+ \nabla p &=& \dive \boldsymbol{\tau} + \rho\mathbf{f}  \\
\frac{\partial (\rho E)}{\partial t}  
+\dive\left[(\rho E+p)\mathbf{u} \right]&=&  \dive  (\boldsymbol{\tau}   \mathbf{u} -\mathbf{q}) +\rho\mathbf{f} \cdot  \mathbf{u}\nonumber
\end{eqnarray}
for $ \mathbf{x}\in\Omega, $ $t\in [0, T_f], $ along with suitable initial and boundary conditions to be discussed later. Here \(T_f\) is the final time, $\rho$ is the density, \(\mathbf{u}\) is the fluid velocity, $p$ is the pressure, $\mathbf{q}$ denotes the heat flux and $\mathbf{f}$ represents volumetric forces. $ \rho E $ is the total energy, which can be rewritten as  $ \rho E = \rho e + \rho k,$ where $e$ is the internal energy and $k = \|\mathbf{u}\|^2/2 $ is the kinetic energy. At this stage, no more specific assumptions are made on the fluid and the choices of equations of state will be specified in the following. We also introduce the specific enthalpy $ h = e + {p}/{\rho} $ and remark that one can also rewrite the energy flux as

$$
\left(\rho E + p\right)\mathbf{u} = \left(e+k+\frac{p}{\rho}\right)\rho\mathbf{u} = \left(h + k\right)\rho\mathbf{u}.
$$
We assume that $\mathbf{q}= -\kappa \nabla T,$ where $T$ denotes the absolute temperature and $\kappa $ the thermal conductivity. Furthermore, we assume that the linear stress constitutive equation holds and we neglect the bulk viscosity and we denote the shear viscosity as $\mu, $ so that

$$
\boldsymbol{\tau} =  
\mu \left (\nabla \mathbf{u}
+\nabla \mathbf{u}^T\right) -\frac{2\mu}{3} (\nabla \cdot \mathbf{u})\mathbf{I}.
$$
For the sake of simplicity, we also assume constant values for both \(\mu\) and \(\kappa\). This choice can be justified by considering that we aim to simulate regimes with limited variations of temperature and density and, moreover, we are mainly interested in time scales where diffusive effects play a less relevant role. The equations can then be rewritten as 

\begin{eqnarray}
\label{eq:ns_comp2}
\frac{\partial \rho}{\partial t} + \dive\left(\rho\mathbf{u} \right) &=& 0   \nonumber \\
\frac{\partial (\rho\mathbf{u})}{\partial t} +
\dive\left(\rho \mathbf{u} \otimes\mathbf{u}\right) + 
\nabla p &=&\mu\dive \left[\left (\nabla \mathbf{u} + \nabla \mathbf{u}^T\right)    
- \frac{2}{3} (\nabla \cdot \mathbf{u})\mathbf{I}\right] \nonumber \\
&+& \rho\mathbf{f} \\
\frac{\partial (\rho E)}{\partial t} + \dive\left[(h+k)\rho\mathbf{u} \right] &=& 
\mu\dive \left[ \left (\nabla \mathbf{u} + \nabla \mathbf{u}^T\right) \mathbf{u}  
-\frac{2}{3} (\nabla \cdot \mathbf{u}) \mathbf{u}\right ]\nonumber \\
&+&\kappa \Delta T+\rho\mathbf{f} \cdot\mathbf{u}.  \nonumber 
\end{eqnarray}
We now introduce reference scaling values $ {\cal L},{\cal T}, {\cal U} $ for the length, time and velocity, respectively, as well as reference values $ {\cal P},{\cal R}, {\cal T\cal T},{\cal E}, {\cal I} $ for pressure, density, temperature, total energy and internal energy, respectively. We assume unit Strouhal number $ St= {\cal L}/{{\cal U}{\cal T}}\approx 1 $, that the enthalpy scales like $ {\cal I} + {\cal P}/{\cal R} $ and that  

$$ {\cal I} \approx \frac{\cal P}{\cal R} \qquad {\cal E} = {\cal I} + {\cal U}^2. $$
The model equations can then be written in non-dimensional form as

\begin{eqnarray}
\label{eq:ns_comp_nondim}
\frac{\partial \rho}{\partial t} + \dive\left(\rho\mathbf{u} \right)&=& 0   \nonumber \\
\frac{\partial \rho\mathbf{u}}{\partial t} + 
\dive\left(\rho \mathbf{u} \otimes\mathbf{u}\right) + 
\frac{\cal P}{{\cal R} {\cal U}^2}\nabla p &=& 
\frac{\mu }{{\cal R} {\cal U}{\cal L}}
\dive \left[ \left (\nabla \mathbf{u} + 
\nabla \mathbf{u}^T\right)  - \frac{2}{3} (\nabla \cdot \mathbf{u})\mathbf{I}  \right] \nonumber \\
&+& \frac{\cal T}{{\cal U}} \rho\mathbf{f}  \\
\frac{\partial \rho E}{\partial t} + 
\dive\left[\left(h\frac{{\cal I} + {\cal P}/{\cal R}}{\cal E}+k\frac{\cal U^2}{\cal E}\right)\rho\mathbf{u} \right] &=& \frac{\mu {\cal U}}{{\cal R} {\cal E}{\cal L}}\dive \left[\left (\nabla \mathbf{u} + 
\nabla \mathbf{u}^T\right)  \mathbf{u} - \frac{2}{3} (\nabla \cdot \mathbf{u}) \mathbf{u}
\right]  \nonumber \\
&+& \frac{\kappa {\cal T\cal T}}{{\cal R}{\cal E}{\cal U}{\cal L} } \Delta T +\frac{\cal L}{{\cal E}} \rho\mathbf{f} \cdot  \mathbf{u}. \nonumber
\end{eqnarray}
Notice that, with a slight abuse of notation, we have kept the same notation for the non-dimensional variables. We then define the Reynolds, Prandtl and Mach numbers as

$$
Re= \frac{{\cal R} {\cal U}{\cal L}}{\mu } \ \ \ \kappa=\frac{c_p\mu}{Pr} \ \ \
M^2= \frac{{\cal R} {\cal U}^2}{\cal P },
$$
where $c_p$ denotes the specific heat at constant pressure, so that, for low and moderate Mach numbers, 

\begin{eqnarray}
\frac{{\cal U}^2}{\cal E}&=&\frac{1}{\frac{\cal I}{{\cal U}^2} + 1} = O(M^2)
\nonumber \\
\frac{{\cal I} + {\cal P}/{\cal R}}{\cal E} &=& \frac{\frac{\cal I}{{\cal U}^2}  + \frac{1}{M^2}}{\frac{\cal I}{{\cal U}^2} + 1} = O(1).
\end{eqnarray}
This justifies, in the above mentioned regimes, methods in which an implicit coupling
between the pressure gradient and the energy flux is enforced. This strategy has been proposed in the seminal paper \cite{casulli:1984} and in the more recent works \cite{dumbser:2016, munz:2003}. We finally assume that the only acting volumetric force is gravity, so that $\mathbf{f} = -g\mathbf{k}, $ where $g$ denotes the acceleration of gravity and $ \mathbf{k} $ the upward pointing unit vector in the standard Cartesian reference frame. It follows that

\begin{eqnarray}
\frac{\cal T}{{\cal U}} \rho g &=& \frac{g{\cal T} {\cal U}}{{\cal U}^2} \rho  = \frac{g{\cal L}}{{\cal U}^2} \rho =\frac{\rho}{Fr^2} \ \ \ \ Fr^2=\frac{{\cal U}^2}{g{\cal L}} \nonumber \\
\frac{\cal L}{{\cal E}}\rho g &=& \frac{g \rho{\cal L}}{{\cal I}+\frac{\cal P}{\cal R} + {\cal U}^2} =
\frac{g \rho{\cal L}}{{\cal U}^2}\frac1{\frac{\cal I}{\cal U}^2 + 1 + \frac{1}{M^2}} = 
\frac{\rho}{Fr^2}O(M^2).
\end{eqnarray}
As a result, we will consider the non-dimensional model equations  

\begin{eqnarray}
\label{eq:ns_comp_nondim_grav}
\frac{\partial \rho}{\partial t}  
+\dive\left(\rho\mathbf{u} \right)&=& 0   \nonumber \\
\frac{\partial \rho\mathbf{u}}{\partial t}  
+\dive\left(\rho \mathbf{u} \otimes\mathbf{u}\right) 
+ \frac{1}{M^2}\nabla p &=& 
\frac{1}{Re}
\dive \left[ \left (\nabla \mathbf{u}
+\nabla \mathbf{u}^T\right)  - \frac{2}{3} (\nabla \cdot \mathbf{u})\mathbf{I}  \right] \nonumber \\
&-&\frac{\rho}{Fr^2} \mathbf{k}\\
\frac{\partial \rho E}{\partial t}  
+\dive\left[\left(h+kM^2\right)\rho\mathbf{u} \right]&=& \frac{M^2}{Re}\dive \left[ \left (\nabla \mathbf{u}
+\nabla \mathbf{u}^T\right)  \mathbf{u}  
-\frac{2}{3} (\nabla \cdot \mathbf{u}) \mathbf{u}
\right ]  \nonumber \\
&+&   \frac{1}{PrRe}\Delta T 
-\rho\frac{M^2}{Fr^2} \mathbf{k} \cdot  \mathbf{u}, \nonumber
\end{eqnarray}
where we have taken $ c_p{\cal T\cal T}/{\cal E} \approx 1 $, which can be justified at moderate values of the Mach number. Notice that these non-dimensional equations are very similar to those derived in \cite{munz:2003}. 

\section{The equation of state}
\label{sec:EOS} \indent

The above equations must be complemented by an equation of state (EOS) for the compressible fluid. A classical choice is that of an ideal gas; in the non-dimensional variables introduced above the equation that links together pressure, density and internal
energy is given by

\begin{equation}
\label{eq:idealgas}
p = \left(\gamma - 1\right)\left(\rho E - \frac{1}{2}M^2\rho\mathbf{u}\cdot\mathbf{u}\right),
\end{equation}
with \(\gamma\) denoting the specific heats ratio. The above relation is valid only in case of constant \(\gamma\) \cite{vidal:2001}. An example of non-ideal gas equation of state is given by the general cubic equation of state, whose equation linking together internal energy, density and temperature, according to \cite{vidal:2001}, is given in dimensional form by:

\begin{equation}\label{eq:caloric_cubic_general}
e = e^{\#}\left(T\right) + \frac{a(T) - T\frac{da}{dT}}{b}U\left(\rho, b, r_{1}, r_{2}\right),
\end{equation}
where \(e^{\#}\) denotes the internal energy of an ideal gas at temperature \(T\). In case \(c_{v} = \frac{de^{\#}}{dT}\) is constant, we can write

\begin{equation}\label{eq:caloric_cubic_T}
e = c_{v}T + \frac{a(T) - T\frac{da}{dT}}{b}U\left(\rho, b, r_{1}, r_{2}\right).
\end{equation}
In case the previous hypothesis does not hold, we analogously define \(\overline{c_{v}}(T) = \frac{e^{\#}(T)}{T}\), so that \eqref{eq:caloric_cubic_general} reads as follows

\begin{equation}\label{eq:caloric_cubic_T_dip}
e = \overline{c_{v}}(T)T + \frac{a(T) - T\frac{da}{dT}}{b}U\left(\rho, b, r_{1}, r_{2}\right).
\end{equation}
The quantity \(\overline{c_{v}}(T)\) should not be understood as a real specific heat, but only as a convenient way of writing the above EOS. The function \(U\) and the constants \(r_{1}\) and \(r_{2}\) depend on the specific equation of state, whereas \(a\) and \(b\) are suitable parameters that characterize the gas behaviour. In this work, we consider the van der Waals EOS, for which \(r_{1} = r_{2} = 0\) and

\begin{equation}
U = -b\rho
\end{equation}
and the Peng-Robinson EOS, for which \(r_{1} = -1 - \sqrt{2}\), \(r_{2} = -1 + \sqrt{2}\) and

\begin{equation}
U = \frac{1}{r_{1} - r_{2}}\log\left(\frac{1 - \rho b r_{1}}{1 - \rho b r_{2}}\right).
\end{equation}
The link between pressure, density and temperature for the general cubic EOS in dimensional form can be expressed as follows:

\begin{equation}\label{eq:thermal_cubic}
p = \frac{\rho R_{g} T}{1 - \rho b} - \frac{a \rho^2}{\left(1 - \rho b r_{1}\right)\left(1 - \rho b r_{2}\right)}, 
\end{equation}
with \(R_{g}\) denoting the specific gas constant. Notice that for \(a = b = 0\), the expression for the pressure of an ideal gas is retrieved. For the sake of clarity, we introduce the following non-dimensional variables

\begin{equation}
\tilde R_{g} = \frac{\cal R \cal T \cal T}{\cal P}R_{g} \qquad \tilde a = a\frac{{\cal R}^2}{\cal P} \qquad \tilde b = {\cal R} b,
\end{equation}
so that \eqref{eq:thermal_cubic} can be rewritten in non-dimensional form as:

\begin{equation}\label{eq:theraml_cubic_adim}
p = \frac{\rho \tilde R_{g} T}{1 - \rho \tilde b} - \frac{\tilde a \rho^2}{\left(1 - \rho \tilde b r_{1}\right)\left(1 - \rho \tilde b r_{2}\right)}. 
\end{equation}
Finally, we define \(\tilde c_{v}\left(T\right) = \overline{c_{v}}\frac{\cal R \cal T \cal T}{\cal P}\), so that the non-dimensional version of \eqref{eq:caloric_cubic_T_dip} reads as follows:

\begin{equation}\label{eq:caloric_cubic_T_dip_adim}
e = \tilde c_{v}(T)T + \frac{\tilde a(T) - T\frac{d \tilde a}{dT}}{\tilde b}U\left(\rho, \tilde b, r_{1}, r_{2}\right).
\end{equation}
The last example of non-ideal gas considered is represented by the Stiffened Gas equation of state (SG-EOS) \cite{metayer:2016}, which is interesting for its convexity property and is given in dimensional variables by:

\begin{equation}\label{eq:sgs}
p = \left(\gamma - 1\right)\left(\rho E - \frac{1}{2}\rho\mathbf{u}\cdot\mathbf{u} - \rho q\right) - \gamma \pi,
\end{equation} 
where \(q\) and \(\pi\) are parameters that determine the gas characteristics. Notice that, for this equation of state, the parameters have to be taken constant \cite{metayer:2016}. We define 

\begin{equation}
\tilde q = \frac{\cal R}{\cal P}q \qquad \tilde \pi = \frac{\pi}{\cal P},
\end{equation}
so that \eqref{eq:sgs} reads in terms of non-dimensional variables as follows:

\begin{equation}\label{eq:sgs_adim}
p = \left(\gamma - 1\right)\left(\rho E - \frac{1}{2}M^2\rho\mathbf{u}\cdot\mathbf{u} - \rho \tilde q\right) - \gamma \tilde \pi.
\end{equation} 
Finally, the link between pressure, density and temperature for the SG-EOS can be written as:

\begin{equation}\label{eq:thermal_sgs}
T = \frac{p + \pi}{\rho\left(\gamma - 1\right)c_{v}}.
\end{equation}
We define \(\tilde c_{v} = c_{v}\frac{\cal R \cal T \cal T}{\cal P}\), so that the non-dimensional version of \eqref{eq:thermal_sgs} is given by:

\begin{equation}\label{eq:thermal_sgs_adim}
T = \frac{p + \tilde \pi}{\rho\left(\gamma - 1\right) \tilde c_{v}}.
\end{equation}
More accurate and general equations of state are available in literature \cite{span:2000},\cite{lemmon:2006}, but the above choices, in particular the cubic EOS, are suitable for the regimes of interest, involve non trivial non-linearities  and provide a sufficient level of complexity for the validation of the proposed numerical scheme. 

An important parameter to determine the regime in which real gas effects are relevant is the so-called compressibility factor. In terms of dimensional variables, it is given by 

\begin{equation}
z = \frac{p}{\rho R_{g} T}.
\end{equation}
When \(z \approx 1\), the gas can be treated as an ideal one, while the ideal gas law is no longer valid for values of $z$ very different from \(1\). 

We also recall here the definition of the potential temperature \(\theta\) for an ideal gas, which is commonly used in applications to atmospheric flows 
\begin{equation}
\theta = \frac{T}{\Pi},
\end{equation}
with \(\Pi = \left(\frac{p}{p_0}\right)^{\frac{\gamma - 1}{\gamma}}\) denoting the so-called Exner pressure. Here, \(p_{0}\) denotes a reference pressure value and, in this work, we consider \(p_{0} = 10^{5} \SI{}{\pascal}\). For isentropic processes, the initial condition is typically given as a perturbation with respect to a constant background potential temperature. Therefore, as discussed in Section \ref{sec:tests}, the gradient of the potential temperature is a good candidate to drive the mesh adaptation procedure.

Our goal is to employ adaptive mesh refinement also in the case of real gases, so as to enhance the computational efficiency, but for non-ideal gases the definition of a potential temperature or of quantities with similar properties is not trivial. We propose here a quantity with a simple definition, stemming from the analysis of isentropic processes, that is valid for the general cubic equation of state in the case \(\frac{da}{dT} = 0\) and \(\frac{d e^{\#}}{dT} = c_{v} = const\). For the sake of simplicity, in order to avoid the influence of reference quantities, we report the computations using dimensional variables. Let us recall the first law of thermodynamics

\begin{equation}
de = Tds - pdv = Tds + \frac{p}{\rho^2}d\rho,
\end{equation}
where \(s\) denotes the specific entropy. Dividing both  sides in the previous equation by \(T\) we obtain

\begin{equation}
\frac{1}{T}de = ds + \frac{p}{\rho^2 T}d\rho
\end{equation} 
which in an isentropic process reduces to

\begin{equation}
\frac{1}{T}de - \frac{p}{\rho^2 T}d\rho = 0.
\end{equation}
\begin{equation}\label{eq:isentropic_1}
\frac{c_{v}}{T}dT + \frac{a}{b}\frac{1}{T}\frac{\partial U}{\partial \rho}d\rho = 0.
\end{equation}
The EOS can be rewritten in the following form \cite{vidal:2001}

\begin{equation}\label{eq:EOS}
T = \left[p + \frac{a\rho^2}{\left(1 - \rho b r_{1}\right)\left(1 - \rho b r_{2}\right)}\right]\frac{\left(1 - \rho b\right)}{\rho R_{g}}.
\end{equation}
If we substitute \eqref{eq:EOS} into \eqref{eq:isentropic_1}, we obtain

\begin{align}\label{eq:isentropic_2}
&\frac{c_{v}}{T}dT + \nonumber \\
&\left(\frac{a}{b}\frac{\partial U}{\partial \rho}\rho - \frac{p}{\rho}\right)\frac{R}{\rho\left(1-\rho b\right)}\frac{\left(1 - \rho b r_{1}\right)\left(1 - \rho b r_{2}\right)}{p\left(1 - \rho b r_{1}\right)\left(1 - \rho b r_{2}\right) + a\rho^2}d\rho = 0.
\end{align}
In the case of van der Waals EOS, \(U = -b\rho\) and \(\frac{\partial U}{\partial \rho} = -b\), whereas in the case of Peng-Robinson EOS one has
$$ U = \frac{1}{r_{1} - r_{2}}\log\left(\frac{1 - \rho b r_{1}}{1 - \rho b r_{2}}\right)  \ \ \ \frac{\partial U}{\partial \rho} = -\frac{b}{\left(1 - \rho b r_{1}\right)\left(1 - \rho b r_{2}\right)}.$$ Since, for van der Waals EOS \(r_{1} = r_{2} = 0\), the expression

\begin{equation}
\frac{\partial U}{\partial \rho} = -\frac{b}{\left(1 - \rho b r_{1}\right)\left(1 - \rho b r_{2}\right)}
\end{equation}
can be applied for both van der Waals and Peng-Robinson EOS. Hence, \eqref{eq:isentropic_2} reduces to

\begin{equation}\label{eq:isentropic_3}
\frac{c_{v}}{T}dT - \frac{R_{g}}{\rho\left(1- \rho b\right)}d\rho = 0,
\end{equation}
which can then be integrated to yield

\begin{equation}
c_{v}\log(T) - 2 R_{g}\atanh\left(2\rho b - 1\right) = const
\end{equation}
or, equivalently,

\begin{equation}
\log(T) - 2\frac{R_{g}}{c_{v}}\atanh\left(2\rho b - 1\right) = const.
\end{equation}
From \eqref{eq:isentropic_3}, it is immediate to verify that, in the non-dimensional case, we obtain

\begin{equation}\label{eq:beta}
\log(T^{*}) - 2\frac{\tilde R_{g}}{\tilde c_{v}}\atanh\left(2\rho^{*} \tilde b - 1\right) = const,
\end{equation}
where the symbol \(*\) denotes non-dimensional variables. In the more general case \(\frac{da}{dT} \neq 0\), it can be shown that \cite{nederstigt:2017}

\begin{equation}\label{eq:beta_1}
\frac{p}{\rho^{\gamma_{p\rho}}} = const,
\end{equation}
where 

\begin{equation}\label{eq:gamma_pv}
\gamma_{p\rho} = \frac{c^2}{M^2}\frac{\rho}{p}.
\end{equation}
The evaluation of this quantity is less straightforward than that of \eqref{eq:beta}, since it involves the computation of non trivial derivatives, see the discussions in Appendix \ref{sec:eigenvalues} and \cite{nederstigt:2017}. Both the conserved quantities will be employed in Section \ref{sec:tests} for adaptive simulations. 

\section{The time discretization strategy}
\label{sec:tdisc} \indent

In the low Mach number limit, terms proportional to $ 1/M^2 $ in \eqref{eq:ns_comp_nondim_grav} yield stiff components of the resulting semidiscretized ODE system. Therefore, following as remarked above \cite{casulli:1984, dumbser:2016}, it is appropriate to couple implicitly the energy equation to the momentum equation, while the continuity equation can be discretized in a fully explicit fashion.  
While this would be sufficient to yield an efficient time discretization approach for the purely hyperbolic system associated to \eqref{eq:ns_comp_nondim_grav} in absence of gravity, in regimes for which

$$Pr\approx O(1), \ \ \ Fr <<1 $$ 
thermal diffusivity and gravity terms would also have to be treated implicitly for the time discretization methods to be efficient. Straightforward application of any monolithic solver would then yield large algebraic systems with multiple couplings between discrete DOF associated to different physical variables. To avoid this, we resort to an operator splitting approach, see e.g. \cite{leveque:2002}, as commonly done in numerical models for atmospheric physics. More specifically, after spatial discretization, all diffusive terms on the right hand side of \eqref{eq:ns_comp_nondim_grav} are split from the hyperbolic part on the left hand side.
The hyperbolic part is treated in a similar fashion to what outlined in \cite{dumbser:2016}, while the diffusive terms are treated implicitly. For simplicity, the gravity terms will be treated explictly in this first attempt and only a basic, first order splitting will be described, which can be easily improved to second order accuracy by the Strang splitting approach \cite{leveque:2002, strang:1968}.

For the time discretization, an IMplicit EXplicit (IMEX) Additive Runge Kutta method (ARK) \cite{kennedy:2003} method will be used. These methods are useful for time dependent problems that can be formulated as
$ \mathbf{y}' = \mathbf{f}_{S}(\mathbf{y},t) +\mathbf{f}_{NS}(\mathbf{y},t) $, where the $S$ and $NS$ subscripts denote the stiff and non-stiff components of the system, to which the implicit and explicit companion methods are applied, respectively. 
If $\mathbf{v}^{n}\approx y(t^n), $ the generic $ s- $ stage IMEX-ARK method can be defined as 

\begin{equation}
\label{imex-ark}
\begin{array}{rcl}
\mathbf{v}^{(n,l)} = \mathbf{v}^n &+& 
\Delta t  \sum\limits_{m=1}^{s-1} \bigg( a_{lm}\mathbf{f}_{NS}(\mathbf{v}^{(n,m)},t+c_m \Delta t) \\[4mm] &+&
\tilde{a}_{lm}\mathbf{f}_{S}(\mathbf{v}^{(n,m)},t+c_m\Delta t) \bigg) + \Delta t \, \tilde{a}_{ll} \, \mathbf{f}_{S}(\mathbf{v}^{(n,l)},t+c_l\Delta t),
\end{array}
\nonumber 
\end{equation}
where $l=1,\dots, s $. After computation of the intermediate stages, $\mathbf{v}^{n+1}$ is computed as

$$
\mathbf{v}^{n+1} = \mathbf{v}^n+\Delta t \sum_{l=1}^{s}b_{l}\left [\mathbf{f}_{NS}(\mathbf{v}^{(n,l)},t+c_l \Delta t)+\mathbf{f}_{S}(\mathbf{v}^{(n,l)},t+c_l\Delta t)\right].
$$

Coefficients $a_{lm}, \tilde{a}_{lm}, c_l$ and $b_l$  are determined so that the method is consistent of a given order. In particular, in addition to the order conditions specific to each sub-method, the coefficients should respect coupling conditions.
Here, we consider a variant of the IMEX method proposed in \cite{giraldo:2013}, whose coefficients are presented in the Butcher tableaux reported in Tables \ref{ark2_butch_e} and \ref{ark2_butch_i} for the explicit and implicit method, respectively, where $\gamma=2 - \sqrt{2}.$ 
The coefficients of the explicit method were proposed in \cite{giraldo:2013}, while the implicit method, also employed in the same paper, coincides indeed for the above choice of $\gamma $ with the TR-BDF2 method proposed in \cite{bank:1985, hosea:1996} and applied to the Euler equations in \cite{tumolo:2015}. Notice that the implicit method was also employed in \cite{orlando:2021} to provide a robust time discretization for the incompressible Navier-Stokes equations, which further guarantees the robustness of the proposed approach in the low Mach number limit.

Finally, even though we focus here on this specific second order method, the same strategy we outline is applicable to a generic DIRK method. In particular, higher order methods could be considered for coupling to high order spatial discretization, even though the effective overall accuracy would be limited by the splitting procedure if gravity and viscous terms are present.

\begin{table}[!h]
	\begin{center}
		\begin{tabular}{c|ccc}
			0 & 0 & &  \\
			$\gamma$ & $\gamma$ & 0 & \\
			1 & $1-\alpha$ & $ \alpha$ & 0 \\
			\hline
			& $\frac{1}{2}-\frac{\gamma}4$ & $ \frac{1}{2}-\frac{\gamma}4$ & $\frac{\gamma}2$
		\end{tabular}
	\end{center}
	\caption{\it Butcher tableaux of the explicit ARK2 method}
	\label{ark2_butch_e}
\end{table}

\begin{table}[h!]
	\begin{center}
		\begin{tabular}{c|ccc}
			0 & 0 & & \\
			$\gamma$ &  $\frac{\gamma}2$ &  $\frac{\gamma}2$ & \\
			1 & $  \frac{1}{2\sqrt{2}}$ & $  \frac{1}{2\sqrt{2}}$ & $1-\frac{1}{\sqrt{2}}$ \\
			\hline
			& $ \frac{1}{2}-\frac{\gamma}4$ & $  \frac{1}{2}-\frac{\gamma}4$ & $\frac{\gamma}2$
		\end{tabular}
	\end{center}
	\caption{\it Butcher tableaux of the implicit ARK2 method}
	\label{ark2_butch_i}
\end{table}
\noindent
Notice that, as discussed in \cite{giraldo:2013}, the choice of the coefficients

$$
\alpha = \frac{7-2\gamma}6 \ \ \ \  1- \alpha = \frac{2\gamma-1}{6} 
$$
in the third stage of the explicit part of the method is arbitrary. In \cite{giraldo:2013}, the above value of $\alpha $  was chosen with the aim of maximizing the stability region of the method. However, if a stability and absolute monotonicity analysis is carried out, as discussed in detail in Appendix \ref{sec:mono_analysis}, it can be seen that different choices might be more advantageous, in order to improve the monotonicity of the method without compromising its stability. In particular, the value of  $ \alpha = 1/2 $ appears to be a more appropriate choice, as also demonstrated by the numerical experiments reported in Section \ref{sec:tests}.

We now describe the application of this IMEX method and of the splitting approach outlined above to equations \eqref{eq:ns_comp_nondim_grav}.
Notice that, for simplicity, we first present the time semi-discretization only, while maintaining the continuous form of \eqref{eq:ns_comp_nondim_grav} with respect to the spatial variables. The detailed description of the algebraic problems resulting from the full space and time discretization according to the method outlined here will be presented in Section \ref{sec:sdisc}.

For each time step, we first consider the discretization of the hyperbolic and forcing terms. For the first stage of the method, one simply has

$$\rho^{(n,1)} = \rho^{n} \ \ \ \  \mathbf{u}^{(n,1)} = \mathbf{u}^n \ \ \ \  E^{(n,1)} = E^n.$$
For the second stage, we can write formally

\begin{eqnarray}
\label{eq:ns_comp_hyp_tdisc1}
\rho^{(n,2)} &=& \rho^n - a_{21} \Delta t\dive\left(\rho^n\mathbf{u}^n \right)   \nonumber \\
\rho^{(n,2)}\mathbf{u}^{(n,2)} &+& 
\tilde a_{22}\frac{\Delta t}{M^2}\nabla p^{(n,2)} = \mathbf{m}^{(n,2)}    \\
\rho^{(n,2)}E^{(n,2)} &+& 
\tilde a_{22} \Delta t \dive\left(h^{(n,2)}\rho^{(n,2)} \mathbf{u}^{(n,2)} \right) = \hat{e}^{(n,2)},
\nonumber
\end{eqnarray}
where we have set

\begin{eqnarray}
\label{eq:ns_comp_hyp_tdisc1_rhss}
\mathbf{m}^{(n,2)} &=& \rho^n\mathbf{u}^n \nonumber \\
&-& a_{21} \Delta t\dive\left(\rho^n \mathbf{u}^n \otimes\mathbf{u}^n\right) 
- \tilde a_{21}\frac{\Delta t}{M^2}\nabla p^n - a_{21} \frac{\Delta t}{Fr^2} \rho^n\mathbf{k}   \\
& \hskip 1cm & \nonumber \\
\hat{e}^{(n,2)} &=& \rho^nE^n - \tilde a_{21} \Delta t \dive \left(h^n  \rho^n\mathbf{u}^n \right)
- a_{21} \Delta t M^2 \dive\left(k^n\rho^n\mathbf{u}^n \right) \nonumber \\
&-&  a_{21} \frac{ \Delta t M^2}{Fr^2}\rho^n \mathbf{k} \cdot  \mathbf{u}^n. \nonumber 
\end{eqnarray}
Notice that, substituting formally $ \rho^{(n,2)}\mathbf{u}^{(n,2)} $ into the energy equation and taking into account the definitions $ \rho E = \rho e +\rho k $ and $ h = e + p/\rho, $ the above system can be solved by computing the solution of

\begin{eqnarray}
\label{eq:ns_comp_hyp_tdisc1_imp}
&&\hspace{0.6cm}\rho^{(n,2)}[e(p^{(n,2)},\rho^{(n,2)} ) +  k^{(n,2)} ] \nonumber \\
&&- \hspace{0.2cm}\tilde a^2_{22}\frac{\Delta t^2}{M^2}\dive \left[\left (e(p^{(n,2)},\rho^{(n,2)})  + \frac{p^{(n,2)}}{\rho^{(n,2)}} \right) \nabla p^{(n,2)}   \right] \\
&&+ \hspace{0.2cm}\tilde a_{22} \Delta t \dive\left[\left (e(p^{(n,2)},\rho^{(n,2)})  + \frac{p^{(n,2)}}{\rho^{(n,2)}} \right ) \mathbf{m}^{(n,2)}  \right]
=  \hat{e}^{(n,2)}\nonumber 
\end{eqnarray}
in terms of $ p^{(n,2)} $ according to the fixed point procedure described in \cite{dumbser:2016}. More specifically, setting $ \xi^{(0)} = p^{(n,2)}, k^{(n,2,0)} = k^{(n,1)}, $ one solves for $l = 1,\dots, L$ the equation

\begin{eqnarray}
\label{eq:nlineq1}
\rho^{(n,2)}e(\xi^{(l+1)},\rho^{(n,2)})   
&-& \tilde a^2_{22}\frac{\Delta t^2}{M^2}\dive \left[\left(e(\xi^{(l)},\rho^{(n,2)})  + \frac{\xi^{(l)}}{\rho^{(n,2)}}\right)  \nabla \xi^{(l+1)}\right] \nonumber\\
&=& \hat{e}^{(n,2)}-  \rho^{(n,2)}k^{(n,2,l)} \\
&-&\tilde a_{22} \Delta t \dive\left[\left (e(\xi^{(l)},\rho^{(n,2)})  + \frac{\xi^{(l)}}{\rho^{(n,2)}} \right ) \mathbf{m}^{(n,2)}  \right] \nonumber
\end{eqnarray}
and updates the velocity as 

$$
\mathbf{u}^{(n,2,l+1)} + \frac{\tilde a_{22} \Delta t}{\rho^{(n,2)}M^2}\nabla \xi^{(l+1)} = \mathbf{m}^{(n,2)}.
$$
In the case of SG-EOS, \(\rho^{(n,2)}e\left(\xi^{(l+1)},\rho^{(n,2)}\right)\) contains a term that only depends on the density, as evident from Equation \eqref{eq:sgs} and, therefore, it has to be properly considered in the right-hand side of \eqref{eq:nlineq1}. On the other hand, the general cubic EOS \eqref{eq:caloric_cubic_T_dip} contains products of quantities depending on temperature and on density. For the sake of simplicity, in order to avoid the solution of a nonlinear equation for each quadrature node, in these cases we keep the temperature at the value in the previous iteration of the fixed point procedure, so as to obtain: 

\begin{eqnarray}
\label{eq:nlineq1_cubic}
\frac{\tilde c_{v} \left( T\left(\xi^{(l)}, \rho^{(n,2)} \right) \right)}{\tilde R_{g}}\xi^{(l+1)}\left(1 - \rho^{(n,2)} \tilde b\right) &-& \nonumber \\
\tilde a^2_{22}\frac{\Delta t^2}{M^2} \dive \left[\left( e\left(\xi^{(l)},\rho^{(n,2)}\right)  + \frac{\xi^{(l)}}{\rho^{(n,2)}} \right) \nabla \xi^{(l+1)} \right] \nonumber &=& \\
\hat{e}^{(n,2)} - \rho^{(n,2)}k^{(n,2,l)} &-& \nonumber \\
\frac{\tilde c_{v}(T\left(\xi^{(l)}, \rho^{(n,2)} \right))}{\tilde R_{g}}\frac{\tilde a\left( T\left(\xi^{(l)}, \rho^{(n,2)} \right)\right) \left(\rho^{(n,2)}\right)^2}{\left(1 - \rho^{(n,2)} \tilde b r_{1}\right)\left(1 - \rho^{(n,2)} \tilde b r_{2}\right)}\left(1 - \rho^{(n,2)} \tilde b\right) &-& \nonumber \\
\frac{\rho^{(n,2)}}{\tilde b}\left[\tilde a \left( T\left(\xi^{(l)}, \rho^{(n,2)} \right) \right) - T \left( \xi^{(l)}, \rho^{(n,2)}\right)\frac{d\tilde a}{dT}\left(\xi^{(l)}, \rho^{(n,2)}\right)\right]U\left(\rho^{(n,2)}\right) &-& \nonumber \\
\tilde a_{22} \Delta t \dive\left[\left( e\left(\xi^{(l)},\rho^{(n,2)}\right)  + \frac{\xi^{(l)}}{\rho^{(n,2)}} \right) \mathbf{m}^{(n,2)}  \right]. 
\end{eqnarray}  
The same considerations as in \cite{dumbser:2016} apply concerning the favourable properties of the weakly nonlinear system resulting from the discrete form of \eqref{eq:nlineq1}. Once the iterations have been completed, one sets
$ \mathbf{u}^{(n,2)} = \mathbf{u}^{(n,2,L+1)}$  and $E^{(n,2)} $ accordingly.
For the third stage, one can write formally  

\begin{eqnarray}
\label{eq:ns_comp_hyp_tdisc2}
\rho^{(n,3)} &=& \rho^n - a_{31} \Delta t\dive\left(\rho^n\mathbf{u}^n \right)
- a_{32} \Delta t\dive\left(\rho^{(n,2)} \mathbf{u}^{(n,2)}  \right)   \nonumber \\
\rho^{(n,3)}\mathbf{u}^{(n,3)} &+&
\tilde a_{33}\frac{\Delta t}{M^2}\nabla p^{(n,3)} =
\mathbf{m}^{(n,3)}  \\
\rho^{(n,3)}E^{(n,3)} &+&
\tilde a_{33} \Delta t \dive\left(h^{(n,3)}\rho^{(n,3)} \mathbf{u}^{(n,3)}  \right)= \hat{e}^{(n,3)}, \nonumber
\nonumber
\end{eqnarray}
where the right hand sides are defined as  

\begin{eqnarray}
\label{eq:ns_comp_tdisc2_rhss}
\mathbf{m}^{(n,3)} &=& \rho^n\mathbf{u}^n  \nonumber \\
&-&a_{31} \Delta t\dive\left(\rho^n \mathbf{u}^n \otimes\mathbf{u}^n\right) 
- \tilde a_{31}\frac{\Delta t}{M^2}\nabla p^n - a_{31} \frac{\Delta t}{Fr^2}\rho^n \mathbf{k}  \\
&-&a_{32} \Delta t\dive\left(\rho^{(n,2)} \mathbf{u}^{(n,2)} \otimes\mathbf{u}^{(n,2)}\right) 
- \tilde a_{32}\frac{\Delta t}{M^2}\nabla p^{(n,2)} - a_{32}\frac{ \Delta t}{Fr^2} \rho^{(n,2)}\mathbf{k} \nonumber \\
& \hskip 1cm & \nonumber \\
\hat{e}^{(n,3)} &=& \rho^nE^n -\tilde a_{31} \Delta t \dive \left(h^n  \rho^n\mathbf{u}^n \right)
- a_{31} \Delta t M^2\dive\left(k^n \rho^n\mathbf{u}^n \right) \nonumber \\
&-&a_{31} \Delta t \frac{M^2}{Fr^2} \rho^n\mathbf{k} \cdot  \mathbf{u}^n \nonumber \\
&-&\tilde a_{32} \Delta t \dive \left(h^{(n,2)} \rho^{(n,2)}\mathbf{u}^{(n,2)} \right)
- a_{32} \Delta t M^2\dive\left(k^{(n,2)} \rho^{(n,2)}\mathbf{u}^{(n,2)} \right) \nonumber \\
&-& a_{32} \Delta t \frac{M^2}{Fr^2} \rho^{(n,2)}\mathbf{k} \cdot  \mathbf{u}^{(n,2)}. \nonumber 
\end{eqnarray}
Again, the solution of this stage is computed by substituting formally $ \rho^{(n,3)}\mathbf{u}^{(n,3)} $ into the energy equation and computing the solution of

\begin{eqnarray}
\label{eq:ns_comp_hyp_tdisc2_imp}
&&\hspace{0.6cm}\rho^{(n,3)}[e \left(p^{(n,3)},\rho^{(n,3)}\right) +  k^{(n,3)} ]\\
&&- \hspace{0.2cm}\tilde a^2_{33}\frac{\Delta t^2}{M^2}\dive \left[ \left (e \left(p^{(n,3)},\rho^{(n,3)}\right) + \frac{p^{(n,3)}}{\rho^{(n,3)}} \right) \nabla p^{(n,3)} \right]\nonumber\\
&&+ \hspace{0.2cm}\tilde a_{33} \Delta t \dive\left[\left( e\left(p^{(n,3)},\rho^{(n,3)}\right)  + \frac{p^{(n,3)}}{\rho^{(n,3)}} \right ) \mathbf{m}^{(n,3)} \right] = \hat{e}^{(n,3)}. \nonumber 
\end{eqnarray}
More specifically, setting $ \xi^{(0)} = p^{(n,3)}, k^{(n,3,0)} = k^{(n,2)}, $
one solves for $ l=1,\dots, L $ the equation

\begin{eqnarray}
\label{eq:nlineq2}
\rho^{(n,3)} e\left(\xi^{(l+1)},\rho^{(n,3)}\right)  
&-& \frac{\tilde a^2_{33} \Delta t^2}{M^2}\dive  \left[  \left (e(\xi^{(l)},\rho^{(n,3)})  + \frac{\xi^{(l)}}{\rho^{(n,3)}} \right )  \nabla \xi^{(l+1)}  \right]
\nonumber\\
&=& \hat{e}^{(n,3)}-  \rho^{(n,3)}k^{(n,3,l)} \\
&-&\tilde a_{33} \Delta t \dive\left[\left( e\left(\xi^{(l)},\rho^{(n,3)}\right)  + \frac{\xi^{(l)}}{\rho^{(n,3)}} \right ) \mathbf{m}^{(n,3)}  \right] \nonumber
\end{eqnarray} 
and updates the velocity as 

$$
\mathbf{u}^{(n,3,l+1)} + \tilde a_{33}\frac{\Delta t}{\rho^{(n,3)}M^2}\nabla \xi^{(l+1)}  = \mathbf{m}^{(n,3)}.
$$
Once again, in case of a non-ideal gas equation of state, the expression of \(\rho^{(n,3)} e\left(\xi^{(l+1)},\rho^{(n,3)}\right)\) is slightly different. For the sake of clarity, we report also for this stage the fixed point equation for the general cubic EOS:

\begin{eqnarray}
\label{eq:nlineq2_cubic}
\frac{\tilde c_{v} \left( T\left(\xi^{(l)}, \rho^{(n,3)} \right) \right)}{\tilde R_{g}}\xi^{(l+1)}\left(1 - \rho^{(n,3)} \tilde b\right) &-& \nonumber \\
\tilde a^2_{33}\frac{\Delta t^2}{M^2} \dive \left[\left( e\left(\xi^{(l)},\rho^{(n,3)}\right)  + \frac{\xi^{(l)}}{\rho^{(n,3)}} \right) \nabla \xi^{(l+1)} \right] &=& \nonumber \\
\hat{e}^{(n,3)} - \rho^{(n,3)}k^{(n,3,l)} &-& \nonumber\\
\frac{\tilde c_{v} \left( T\left(\xi^{(l)}, \rho^{(n,3)} \right) \right)}{\tilde R_{g}}\frac{\tilde a\left( T\left(\xi^{(l)}, \rho^{(n,3)} \right)\right) \left(\rho^{(n,3)}\right)^2}{\left(1 - \rho^{(n,3)} \tilde b r_{1}\right)\left(1 - \rho^{(n,3)} \tilde b r_{2}\right)}\left(1 - \rho^{(n,3)} \tilde b\right) &-& \nonumber \\
\frac{\rho^{(n,3)}}{\tilde b}\left[\tilde a \left( T\left(\xi^{(l)}, \rho^{(n,3)} \right) \right) - T \left( \xi^{(l)}, \rho^{(n,3)}\right)\frac{d\tilde a}{dT}\left(\xi^{(l)}, \rho^{(n,3)}\right)\right]U\left(\rho^{(n,3)}\right) &-& \nonumber \\
\tilde a_{33} \Delta t \dive\left[\left( e\left(\xi^{(l)},\rho^{(n,3)}\right)  + \frac{\xi^{(l)}}{\rho^{(n,3)}} \right) \mathbf{m}^{(n,3)}  \right].
\end{eqnarray} 
\
\\ 
Let us consider now the diffusive part of the Navier-Stokes equations that, as already mentioned in Section \ref{sec:intro}, will be treated with an operator splitting technique. For the sake of clarity, we denote with \(\sim\) the quantities computed in this part of the scheme; hence, we define 

\[\mathbf{\tilde u}^{(n,1)} = \mathbf{u}^{(n,3)} \qquad {\tilde E}^{(n,1)} = E^{(n,3)}\]
and we proceed to the discretization of the viscous terms, which is carried out by the implicit part of the IMEX method previously described:

\begin{eqnarray}
\label{eq:ns_comp_tdisc_split_1}
\rho^{n+1}\mathbf{\tilde u}^{(n,2)}  
&-&\tilde a_{22}\frac {\Delta t}{Re} \dive \left[ \left (\nabla \mathbf{\tilde u}
+\nabla \mathbf{\tilde u}^T\right)  - \frac{2}{3} (\nabla \cdot \mathbf{\tilde u})\mathbf{I}  \right]^{(n,2)}
=\mathbf{\tilde m}^{(n,2)}    \\
\rho^{n+1}\tilde E^{(n,2)}
&-& \tilde a_{22}  \frac{\Delta tM^2}{Re}\dive \left[ \left (\nabla \mathbf{\tilde u}
+\nabla \mathbf{\tilde u}^T\right)  \mathbf{\tilde u}  
-\frac{2}{3} (\nabla \cdot \mathbf{\tilde u}) \mathbf{\tilde u}
\right ]^{(n,2)}   \nonumber \\
&-&   \tilde a_{22}  \frac{\Delta t }{PrRe}\Delta \tilde T^{(n,2)}
= {\tilde e}^{(n,2)},
\nonumber
\end{eqnarray}
where we have set

\begin{eqnarray}
\label{eq:ns_comp_tdisc1_split_rhss}
\mathbf{\tilde m}^{(n,2)} &=& \rho^{n+1}\mathbf{\tilde u}^{(n,1)} + \tilde a_{21}\frac{\Delta t}{Re} \dive \left[ \left (\nabla \mathbf{\tilde u}
+\nabla \mathbf{\tilde u}^T\right)  - \frac{2}{3} (\nabla \cdot \mathbf{\tilde u})\mathbf{I}  \right]^{(n,1)} \nonumber \\
& \hskip 1cm & \nonumber \\
{\tilde e}^{(n,2)} &=&\rho^{n+1}\tilde E^{(n,1)}    \nonumber \\
&+& \tilde a_{21}  \frac{\Delta tM^2}{Re}\dive \left[ \left (\nabla \mathbf{u}
+\nabla \mathbf{u}^T\right)  \mathbf{u}  
-\frac{2}{3} (\nabla \cdot \mathbf{u}) \mathbf{u}
\right ]^{(n,1)}   \nonumber \\
&+&   \tilde a_{21}  \frac{\Delta t }{PrRe}\Delta \tilde T^{(n,1)}. \nonumber 
\end{eqnarray}
Notice that the momentum equation in \eqref{eq:ns_comp_tdisc_split_1} is decoupled from the energy equation and can be solved independently, so that in a subsequent step the equation for $\tilde E^{(n,2)}$ can be solved using temperature as an unknown. It is worth to mention that, in case \(\frac{da}{dT} \neq 0\) or \(\frac{d c_{v}}{dT} \neq 0 \), for the cubic EOS, we end up with a non-linear equation. The following fixed point procedure is considered: setting $ \xi^{(0)} = \tilde T^{(n,1)}, $ one solves for $l = 1,\dots, L$

\begin{eqnarray}
&&\tilde c_{v}\left( \xi^{(l)} \right) \xi^{(l+1)} + \frac{\tilde a \left(\xi^{(l)}\right) - \xi^{(l+1)} \frac{d \tilde a}{dT}\left( \xi^{(l)} \right)}{b} \nonumber \\
&&- \tilde a_{22} \frac{\Delta t M^2}{Re}\dive \left[\left (\nabla \mathbf{\tilde u}
+\nabla \mathbf{\tilde u}^T\right) \mathbf{\tilde u}  
-\frac{2}{3} (\nabla \cdot \mathbf{\tilde u}) \mathbf{\tilde u}
\right ]^{(n,2)}   \nonumber \\
&&- \tilde a_{22}  \frac{\Delta t}{Pr Re}\Delta \xi^{(l+1)}
= {\tilde e}^{(n,2)}.
\end{eqnarray}
For the third stage, one can write formally  

\begin{eqnarray}
\label{eq:ns_comp_tdisc2_split}
\rho^{n+1}\mathbf{\tilde u}^{(n,3)}  
&-&\tilde a_{33} \frac { \Delta t}{Re} \dive \left[ \left (\nabla \mathbf{\tilde u}
+\nabla \mathbf{\tilde u}^T\right)  - \frac{2}{3} (\nabla \cdot \mathbf{\tilde u})\mathbf{I}  \right]^{(n,3)}
=\mathbf{\tilde m}^{(n,3)} \nonumber \\
\rho^{n+1}\tilde E^{(n,3)}
&-& \tilde a_{33}  \frac{\Delta tM^2}{Re}\dive \left[ \left (\nabla \mathbf{\tilde u}
+\nabla \mathbf{\tilde u}^T\right)  \mathbf{\tilde u}  
-\frac{2}{3} (\nabla \cdot \mathbf{\tilde u}) \mathbf{\tilde u}
\right ]^{(n,3)}   \nonumber \\
&-& \tilde a_{33}  \frac{\Delta t }{PrRe}\Delta \tilde T^{(n,3)} = \mathbf{\tilde e}^{(n,3)}, \nonumber
\nonumber
\end{eqnarray}
where the right hand sides are defined as  

\begin{eqnarray}
\label{eq:ns_comp_tdisc2_rhss_split}
\mathbf{\tilde m}^{(n,3)} &=& \rho^{n+1}\mathbf{\tilde u}^{(n,1)} + \tilde a_{31}\frac { \Delta t}{Re} \dive \left[ \left (\nabla \mathbf{\tilde u}
+\nabla \mathbf{\tilde u}^T\right)  - \frac{2}{3} (\nabla \cdot \mathbf{\tilde u})\mathbf{I}  \right]^{(n,1)} \nonumber \\
&+& \tilde a_{32}\frac { \Delta t}{Re} \dive \left[ \left (\nabla \mathbf{ \tilde u}
+\nabla \mathbf{ \tilde u}^T\right)  - \frac{2}{3} (\nabla \cdot \mathbf{ \tilde u})\mathbf{I}  \right]^{(n,2)} \\
& \hskip 1cm & \nonumber \\
\mathbf{\tilde e}^{(n,3)} &=& \rho^{n+1}\tilde E^{(n,1)} \\
&+& \tilde a_{31}  \frac{\Delta t M^2}{Re}\dive \left[ \left (\nabla \mathbf{ \tilde u}
+\nabla \mathbf{\tilde u}^T\right)  \mathbf{\tilde u}  
-\frac{2}{3} (\nabla \cdot \mathbf{\tilde u}) \mathbf{\tilde u}
\right ]^{(n,1)}  \nonumber \\
&+& \tilde a_{32}  \frac{\Delta t M^2}{Re}\dive \left[ \left (\nabla \mathbf{ \tilde u}
+\nabla \mathbf{\tilde u}^T\right)  \mathbf{\tilde u}  
-\frac{2}{3} (\nabla \cdot \mathbf{\tilde u}) \mathbf{\tilde u}
\right ]^{(n,2)}  \nonumber \\
&+&  \tilde a_{31} \frac{ \Delta t }{PrRe}\Delta \tilde T^{(n,1)} + \tilde a_{32} \frac{ \Delta t }{PrRe}\Delta \tilde T^{(n,2)} . \nonumber 
\end{eqnarray}
Again, the momentum equation in \eqref{eq:ns_comp_tdisc2_split} is decoupled from the energy equation and can be solved independently, so that in a subsequent step the equation for $E^{(n,3)}$ can be solved using temperature as an unkown.
Finally, one sets 

$$ \mathbf{u}^{n+1} = \mathbf{\tilde u}^{(n,3)}  \ \ \ \  E^{n+1} =\tilde E^{(n,3)}$$
and the computation of the $n-$th time step is completed.

\section{The spatial discretization strategy}
\label{sec:sdisc} \indent

We consider a decomposition of the domain \(\Omega\) into a family of hexahedra \(\mathcal{T}_h\) (quadrilaterals in the two-dimensional case) and denote each element by \(K\). The skeleton \(\mathcal{E}\) denotes the set of all element faces and \(\mathcal{E} = \mathcal{E}^{I} \cup \mathcal{E}^{B}\), where \(\mathcal{E}^{I}\) is the subset of interior faces and \(\mathcal{E}^{B}\) is the subset of boundary faces. Suitable jump and average operators can then be defined as customary for finite element discretizations. A face \(\Gamma \in \mathcal{E}^{I}\) shares two elements that we denote by \(K^{+}\) with outward unit normal \(\mathbf{n}^{+}\) and \(K^{-}\) with outward unit normal \(\mathbf{n}^{-}\), whereas for a face \(\Gamma \in \mathcal{E}^{B}\) we denote by \(\mathbf{n}\) the outward unit normal. For a scalar function \(\varphi\) the jump is defined as

\[\left[\left[\varphi\right]\right] = \varphi^{+}\mathbf{n}^{+} + \varphi^{-}\mathbf{n}^{-} \quad \text{if }\Gamma \in \mathcal{E}^{I} \qquad \left[\left[\varphi\right]\right] = \varphi\mathbf{n} \quad \text{if }\Gamma \in \mathcal{E}^{B}.\]
The average is defined as

\[\left\{\left\{\varphi\right\}\right\} = \frac{1}{2}\left(\varphi^{+} + \varphi^{-}\right) \quad \text{if }\Gamma \in \mathcal{E}^{I} \qquad \left\{\left\{\varphi\right\}\right\} = \varphi \quad \text{if }\Gamma \in \mathcal{E}^{B}.\]
Similar definitions apply for a vector function \(\boldsymbol{\varphi}\):

\begin{align*}
&\left[\left[\boldsymbol{\varphi}\right]\right] = \boldsymbol{\varphi}^{+}\cdot\mathbf{n}^{+} 
+\boldsymbol{\varphi}^{-}\cdot\mathbf{n}^{-} \quad \text{if }\Gamma \in \mathcal{E}^{I} \qquad 
\left[\left[\boldsymbol{\varphi}\right]\right] = \boldsymbol{\varphi}\cdot\mathbf{n} \quad \text{if }\Gamma \in \mathcal{E}^{B} \\
&\left\{\left\{\boldsymbol{\varphi}\right\}\right\} = \frac{1}{2}\left(\boldsymbol{\varphi}^{+} + \boldsymbol{\varphi}^{-}\right) \quad \text{if }\Gamma \in \mathcal{E}^{I} \qquad \left\{\left\{\boldsymbol{\varphi}\right\}\right\} = \boldsymbol{\varphi} \quad \text{if }\Gamma \in \mathcal{E}^{B}.
\end{align*}
For vector functions, it is also useful to define a tensor jump as:

\[\left<\left<\boldsymbol{\varphi}\right>\right> = \boldsymbol{\varphi}^{+}\otimes\mathbf{n}^{+} 
+ \boldsymbol{\varphi}^{-}\otimes\mathbf{n}^{-} \quad \text{if }\Gamma \in \mathcal{E}^{I} 
\qquad \left<\left<\boldsymbol{\varphi}\right>\right> = \boldsymbol{\varphi}\otimes\mathbf{n} \quad \text{if }\Gamma \in \mathcal{E}^{B}.\]
We also introduce the following finite element spaces

\[Q_{r} = \left\{v \in L^2(\Omega) : v\rvert_K \in \mathbb{Q}_{r} \quad \forall K \in \mathcal{T}_h\right\} \ \ \
\mathbf{V}_{r} = \left[Q_{r}\right]^d,\]
where \(\mathbb{Q}_{r}\) is the space of polynomials of degree \(r\) in each coordinate direction. We then denote by \(\boldsymbol{\varphi}_i(\mathbf{x})\) the basis functions for the space \(\mathbf{V}_{r}\) and by \(\psi_i(\mathbf{x})\) the basis functions for the space \(Q_{r}\), the finite element spaces chosen for the discretization of the velocity and of the pressure (as well as the density), respectively.

\begin{equation*}
\mathbf{u}\approx \sum_{j = 1}^{\text{dim}(\mathbf{V}_{r})}u_j(t)\boldsymbol{\varphi}_j(\mathbf{x}) \qquad p \approx \sum_{j = 1}^{\text{dim}(Q_{r})}p_j(t)\psi_j(\mathbf{x}).
\end{equation*}
Given these definitions, the weak formulation for the momentum equation of the second stage \eqref{eq:ns_comp_hyp_tdisc1}
reads as follows:

\begin{align}
\label{eq:momentum_stage_2}
&\qquad \sum_K\int_K\rho^{(n,2)}\mathbf{u}^{(n,2)}\cdot\mathbf{v}d\Omega \nonumber\\ &-\quad \sum_K\int_K\tilde{a}_{22}\frac{\Delta t}{M^2}p^{(n,2)}\nabla\cdot\mathbf{v}d\Omega + \sum_{\Gamma\in\mathcal{E}}\int_{\Gamma}\tilde{a}_{22}\frac{\Delta t}{M^2}\left\{\left\{p^{(n,2)}\right\}\right\}\left[\left[\mathbf{v}\right]\right]d\Sigma \nonumber \\
&=\quad \sum_K\int_K\rho^{n}\mathbf{u}^{n}\cdot\mathbf{v}d\Omega - \sum_K\int_Ka_{21}\frac{\Delta t}{Fr^2}\rho^{n}\mathbf{k}\cdot\mathbf{v}d\Omega \nonumber \\
&+\quad \sum_K\int_K a_{21}\Delta t\left(\rho^n\mathbf{u}^n\otimes\mathbf{u}^n\right):\nabla\mathbf{v}d\Omega + \sum_K\int_K\tilde{a}_{21}\frac{\Delta t}{M^2}p^n\nabla\cdot\mathbf{v}d\Omega \nonumber \\
&-\quad \sum_{\Gamma\in\mathcal{E}}\int_{\Gamma}a_{21}\Delta t\left\{\left\{\rho^n\mathbf{u}^n\otimes\mathbf{u}^n\right\}\right\}:\left<\left<\mathbf{v}\right>\right>d\Sigma \nonumber \\
&-\quad \sum_{\Gamma\in\mathcal{E}}\int_{\Gamma}\tilde{a}_{21}\frac{\Delta t}{M^2}\left\{\left\{p^n\right\}\right\}\left[\left[\mathbf{v}\right]\right]d\Sigma \nonumber \\
&-\quad \sum_{\Gamma\in\mathcal{E}}\int_{\Gamma}a_{21}\Delta t\frac{\lambda^{(n,1)}}{2}\left<\left<\rho^{n}\mathbf{u}^{n}\right>\right>:\left<\left<\mathbf{v}\right>\right>d\Sigma, 
\end{align}
where

\begin{equation*}
\lambda^{(n,1)} = \max\left(\left|\mathbf{u}^{n^{+}}\cdot\mathbf{n}^{+}\right|, \left|\mathbf{u}^{n^{-}}\cdot\mathbf{n}^{-}\right|\right).
\end{equation*}
One can notice that centered flux has been employed as numerical flux for the quantities defined implicitly, whereas an upwind flux has been used for the quantities computed explicitly. In view of the implicit coupling between the momentum and the energy equations, we need to derive the algebraic formulation of \eqref{eq:momentum_stage_2} in order to formally substitute the degrees of freedom of the velocity into the algebraic formulation of the energy equation. We take \(\mathbf{v} = \boldsymbol{\varphi}_i\), \(i = 1...\text{dim}(\mathbf{V}_{r})\) and we exploit the representation introduced above to obtain

\begin{align}
&\qquad \sum_K\int_K\rho^{(n,2)}\sum_{j=1}^{\text{dim}(\mathbf{V}_{r})}u_j^{(n,2)}\boldsymbol{\varphi}_j\cdot\boldsymbol{\varphi}_id\Omega - \sum_K\int_K\tilde{a}_{22}\frac{\Delta t}{M^2}\sum_{j=1}^{\text{dim}(Q_{r})}p_j^{(n,2)}\psi_j\nabla\cdot\boldsymbol{\varphi}_id\Omega\nonumber \\ 
&+\quad \sum_{\Gamma\in\mathcal{E}}\int_{\Gamma}\tilde{a}_{22}\frac{\Delta t}{M^2}\sum_{j=1}^{\text{dim}(Q_{r})}p_j^{(n,2)}\left\{\left\{\psi_j\right\}\right\}\left[\left[\boldsymbol{\varphi}_i\right]\right]d\Sigma \nonumber \\
&=\quad \sum_K\int_K\rho^{n}\mathbf{u}^{n}\cdot\boldsymbol{\varphi}_id\Omega - \sum_K\int_Ka_{21}\frac{\Delta t}{Fr^2}\rho^{n}\mathbf{k}\cdot\boldsymbol{\varphi}_id\Omega \nonumber \\
&+\quad \sum_K\int_K a_{21}\Delta t\left(\rho^n\mathbf{u}^n\otimes\mathbf{u}^n\right):\nabla\boldsymbol{\varphi}_id\Omega + \sum_K\int_K\tilde{a}_{21}\frac{\Delta t}{M^2}p^n\nabla\cdot\boldsymbol{\varphi}_id\Omega \nonumber \\
&-\quad \sum_{\Gamma\in\mathcal{E}}\int_{\Gamma}a_{21}\Delta t\left\{\left\{\rho^n\mathbf{u}^n\otimes\mathbf{u}^n\right\}\right\}:\left<\left<\boldsymbol{\varphi}_i\right>\right>d\Sigma \nonumber \\
&-\quad \sum_{\Gamma\in\mathcal{E}}\int_{\Gamma}\tilde{a}_{21}\frac{\Delta t}{M^2}\left\{\left\{p^n\right\}\right\}\left[\left[\boldsymbol{\varphi}_i\right]\right]d\Sigma \nonumber \\
&-\quad \sum_{\Gamma\in\mathcal{E}}\int_{\Gamma}a_{21}\Delta t\frac{\lambda^{(n,1)}}{2}\left<\left<\rho^{n}\mathbf{u}^{n}\right>\right>:\left<\left<\boldsymbol{\varphi}_i\right>\right>d\Sigma, \label{eq:algebraic_momentum_stage_2}
\end{align}
which can be written in compact form as

\begin{equation}
\mathbf{A}^{(n,2)}\mathbf{U}^{(n,2)} + \mathbf{B}^{(n,2)}\mathbf{P}^{(n,2)} = \mathbf{F}^{(n,2)}
\end{equation}
where we have set

\begin{align}
A_{ij}^{(n,2)} &= \sum_K\int_K\rho^{(n,2)}\boldsymbol{\varphi}_j\cdot\boldsymbol{\varphi}_id\Omega \\
B_{ij}^{(n,2)} &= \sum_K\int_K-\tilde{a}_{22}\frac{\Delta t}{M^2}\nabla\cdot\boldsymbol{\varphi}_i\psi_jd\Omega + \sum_{\Gamma\in\mathcal{E}}\int_{\Gamma}\tilde{a}_{22}\frac{\Delta t}{M^2}\left\{\left\{\psi_j\right\}\right\}\left[\left[\boldsymbol{\varphi}_i\right]\right]d\Sigma 
\end{align}
with \(\mathbf{U}^{(n,2)}\) denoting the vector of the degrees of freedom associated to the velocity field and \(\mathbf{P}^{(n,2)}\) denoting the vector of the degrees of freedom associated to the pressure. 
Consider now the weak formulation for the energy equation of the second stage \eqref{eq:ns_comp_hyp_tdisc1}

\begin{align}
&\qquad \sum_K\int_K\rho^{(n,2)}E^{(n,2)}wd\Omega - \sum_K\int_K\tilde{a}_{22}\Delta th^{(n,2)}\rho^{(n,2)}\mathbf{u}^{(n,2)}\cdot\nabla wd\Omega \nonumber \\ &+\quad \sum_{\Gamma\in\mathcal{E}}\int_{\Gamma}\tilde{a}_{22}\Delta t\left\{\left\{h^{(n,2)}\rho^{(n,2)}\mathbf{u}^{(n,2)}\right\}\right\}\cdot\left[\left[w\right]\right]d\Sigma \nonumber \\
&=\quad \sum_K\int_K\rho^{n}E^{n}wd\Omega - \sum_K\int_Ka_{21}\frac{\Delta tM^2}{Fr^2}\rho^{n}\mathbf{k}\cdot\mathbf{u}^{n}wd\Omega \nonumber \\
&+\quad \sum_K\int_Ka_{21}\Delta tM^2\left(k^{n}\rho^n\mathbf{u}^n\right)\cdot\nabla wd\Omega + \sum_K\int_K\tilde{a}_{21}\Delta t\left(h^{n}\rho^n\mathbf{u}^n\right)\cdot\nabla wd\Omega \nonumber \\
&-\quad \sum_{\Gamma\in\mathcal{E}}\int_{\Gamma}a_{21}\Delta tM^2\left\{\left\{k^n\rho^n\mathbf{u}^n\right\}\right\}\cdot\left[\left[w\right]\right]d\Sigma \nonumber \\
&-\quad \sum_{\Gamma\in\mathcal{E}^{I}}\int_{\Gamma}\tilde{a}_{21}\Delta t\left\{\left\{h^n\rho^n\mathbf{u}^n\right\}\right\}\cdot\left[\left[w\right]\right]d\Sigma \nonumber \\
&-\quad \sum_{\Gamma\in\mathcal{E}}\int_{\Gamma}a_{21}\Delta t\frac{\lambda^{(n,1)}}{2}\left[\left[\rho^{n}E^{n}\right]\right] \cdot \left[\left[w\right]\right]d\Sigma. \label{eq:energy_stage_2}
\end{align}
Notice that, while the fully discrete formulation is presented here for the case of an ideal gas, in the more general case it has to be modified properly as already shown in \eqref{eq:nlineq1_cubic} for the semi discrete formulation. Take \(w = \psi_i\) and consider the expansion for \(\mathbf{u}^{(n,2)}\) in \eqref{eq:energy_stage_2} to get

\begin{align}
&\qquad\sum_K\int_K\rho^{(n,2)}E^{(n,2)}\psi_id\Omega - \sum_K\int_K\tilde{a}_{22}\Delta th^{(n,2)}\rho^{(n,2)}\sum_{j=1}^{\text{dim}(\mathbf{V}_{r})}u_j^{(n,2)}\boldsymbol{\varphi}_j\cdot\nabla\psi_id\Omega \nonumber \\ 
&+\quad\sum_{\Gamma\in\mathcal{E}}\int_{\Gamma}\tilde{a}_{22}\Delta t\sum_{j=1}^{\text{dim}(\mathbf{V}_{r})} u_j^{(n,2)}\left\{\left\{h^{(n,2)}\rho^{(n,2)}\boldsymbol{\varphi}_j\right\}\right\}\cdot\left[\left[\psi_i\right]\right]d\Sigma \nonumber \\
&=\quad\sum_K\int_K\rho^{n}E^{n}\psi_id\Omega - \sum_K\int_Ka_{21}\frac{\Delta tM^2}{Fr^2}\rho^{n}\mathbf{k}\cdot\mathbf{u}^{n}\psi_id\Omega \nonumber \\
&+\quad\sum_K\int_Ka_{21}\Delta tM^2\left(k^{n}\rho^n\mathbf{u}^n\right)\cdot\nabla \psi_id\Omega + \sum_K\int_K\tilde{a}_{21}\Delta t\left(h^{n}\rho^n\mathbf{u}^n\right)\cdot\nabla \psi_id\Omega \nonumber \\
&-\quad\sum_{\Gamma\in\mathcal{E}}\int_{\Gamma}a_{21}\Delta tM^2\left\{\left\{k^n\rho^n\mathbf{u}^n\right\}\right\}\cdot\left[\left[\psi_i\right]\right]d\Sigma \nonumber \\
&-\quad\sum_{\Gamma\in\mathcal{E}}\int_{\Gamma}\tilde{a}_{21}\Delta t\left\{\left\{h^n\rho^n\mathbf{u}^n\right\}\right\}\cdot\left[\left[\psi_i\right]\right]d\Sigma \nonumber \\
&-\quad\sum_{\Gamma\in\mathcal{E}}\int_{\Gamma}a_{21}\Delta t\frac{\lambda^{(n,1)}}{2}\left[\left[\rho^{n}E^{n}\right]\right] \cdot \left[\left[\psi_{i}\right]\right]d\Sigma, \label{eq:albegraic_energy_stage_2}
\end{align}
which can be expressed in compact form as

\begin{equation}
\mathbf{C}^{(n,2)}\mathbf{U}^{(n,2)} = \mathbf{G}^{(n,2)}
\end{equation}
where we have set

\begin{align}
C_{ij}^{(n,2)} &= \sum_K\int_K-\tilde{a}_{22}\Delta th^{(n,2)}\rho^{(n,2)}\boldsymbol{\varphi}_j\cdot\nabla\psi_id\Omega \nonumber \\ &+\sum_{\Gamma\in\mathcal{E}}\int_{\Gamma}\tilde{a}_{22}\Delta t \left\{\left\{h^{(n,2)}\rho^{(n,2)}\boldsymbol{\varphi}_j\right\}\right\}\cdot\left[\left[\psi_i\right]\right]d\Sigma. 
\end{align}
Formally we can then derive \(\mathbf{U}^{(n,2)} = (\mathbf{A}^{(n,2)})^{-1}\left(\mathbf{F}^{(n,2)} - \mathbf{B}^{(n,2)}\mathbf{P}^{(n,2)}\right)\) and obtain the following relation 

\begin{equation}
\mathbf{C}^{(n,2)}(\mathbf{A}^{(n,2)})^{-1} \left(\mathbf{F}^{(n,2)} - \mathbf{B}^{(n,2)}\mathbf{P}^{(n,2)}\right) = \mathbf{G}^{(n,2)}.
\end{equation}
Taking into account that

$$
\rho^{(n,2)}E^{(n,2)} = \rho^{(n,2)}e^{(n,2)}(p^{(n,2)}) + M^2\rho^{(n,2)}k^{(n,2)},
$$ 
we finally obtain 

\begin{equation}
\mathbf{C}^{(n,2)}(\mathbf{A}^{(n,2)})^{-1}\left(\mathbf{F}^{(n,2)} - \mathbf{B}^{(n,2)}\mathbf{P}^{(n,2)}\right) = -\mathbf{D}^{(n,2)}\mathbf{P}^{(n,2)} +  \tilde{\mathbf{G}}^{(n,2)}
\end{equation}
where we have set

\begin{align}
D_{ij}^{(n,2)} = &\sum_K\int_K\rho^{(n,2)}e^{(n,2)}(\psi_j)\psi_id\Omega
\end{align}
and \(\tilde{\mathbf{G}}^{(n,2)}\) takes into account all the other terms (the one at previous stage and the kinetic energy). The above system can be solved in terms of \(\mathbf{P}^{(n,2)}\) according to the fixed point procedure described in \cite{dumbser:2016}. More specifically, setting \(\mathbf{P}^{(n,2,0)} = \mathbf{P}^{(n,1)}, k^{(n,2,0)} = k^{(n,1)}\), for \(l = 1,\dots, M\) one solves the equation

\begin{equation}
\left(\mathbf{D}^{(n,2,l)} - \mathbf{C}^{(n,2,l)}(\mathbf{A}^{(n,2)})^{-1}
\mathbf{B}^{(n,2)}\right)\mathbf{P}^{(n,2,l+1)} =  
\tilde{\mathbf{G}}^{(n,2,l)} - \mathbf{C}^{(n,2,l)}(\mathbf{A}^{(n,2)})^{-1}\mathbf{F}^{(n,2,l)}
\nonumber
\end{equation}
and updates the velocity solving

\begin{equation}
\mathbf{A}^{(n,2)}\mathbf{U}^{(n,2,l+1)} = \mathbf{F}^{(n,2,l)} - \mathbf{B}^{(n,2)}\mathbf{P}^{(n,2,l+1)}. \nonumber
\end{equation}
For the third stage, we proceed in a similar manner. We start with the weak formulation of the momentum equation in \eqref{eq:ns_comp_hyp_tdisc2}:

\begin{align}
&\qquad\sum_K\int_K\rho^{(n,3)}\mathbf{u}^{(n,3)}\cdot\mathbf{v}d\Omega \nonumber \\
&-\quad\sum_K\int_K\tilde{a}_{33}\frac{\Delta t}{M^2}p^{(n,3)}\nabla\cdot\mathbf{v}d\Omega + \sum_{\Gamma\in\mathcal{E}}\int_{\Gamma}\tilde{a}_{33}\frac{\Delta t}{M^2}\left\{\left\{p^{(n,3)}\right\}\right\}\left[\left[\mathbf{v}\right]\right]d\Sigma \nonumber \\
&=\quad\sum_K\int_K\rho^{n}\mathbf{u}^{n}\cdot\mathbf{v}d\Omega - \sum_K\int_Ka_{31}\frac{\Delta t}{Fr^2}\rho^{n}\mathbf{k}\cdot\mathbf{v}d\Omega - \sum_K\int_Ka_{32}\frac{\Delta t}{Fr^2}\rho^{(n,2)}\mathbf{k}\cdot\mathbf{v}d\Omega \nonumber \\
&+\quad\sum_K\int_Ka_{31}\Delta t\left(\rho^n\mathbf{u}^n\otimes\mathbf{u}^n\right):\nabla\mathbf{v}d\Omega + \sum_K\int_K\tilde{a}_{31}\frac{\Delta t}{M^2}p^n\nabla\cdot\mathbf{v}d\Omega \nonumber \\
&+\quad\sum_K\int_Ka_{32}\Delta t\left(\rho^{(n,2)}\mathbf{u}^{(n,2)}\otimes\mathbf{u}^{(n,2)}\right):\nabla\mathbf{v}d\Omega + \sum_K\int_K\tilde{a}_{32}\frac{\Delta t}{M^2}p^{(n,2)}\nabla\cdot\mathbf{v}d\Omega \nonumber \\
&-\quad\sum_{\Gamma\in\mathcal{E}}\int_{\Gamma}a_{31}\Delta t\left\{\left\{\rho^n\mathbf{u}^n\otimes\mathbf{u}^n\right\}\right\}:\left<\left<\mathbf{v}\right>\right>d\Sigma \nonumber \\ 
&-\quad\sum_{\Gamma\in\mathcal{E}}\int_{\Gamma}\tilde{a}_{31}\frac{\Delta t}{M^2}\left\{\left\{p^n\right\}\right\}\left[\left[\mathbf{v}\right]\right]d\Sigma \nonumber \\
&-\quad\sum_{\Gamma\in\mathcal{E}}\int_{\Gamma}a_{32}\Delta t\left\{\left\{\rho^{(n,2)}\mathbf{u}^{(n,2)}\otimes\mathbf{u}^{(n,2)}\right\}\right\}:\left<\left<\mathbf{v}\right>\right>d\Sigma \nonumber \\
&-\quad\sum_{\Gamma\in\mathcal{E}}\int_{\Gamma}\tilde{a}_{32}\frac{\Delta t}{M^2}\left\{\left\{p^{(n,2)}\right\}\right\}\left[\left[\mathbf{v}\right]\right]d\Sigma \nonumber \\
&-\quad\sum_{\Gamma\in\mathcal{E}}\int_{\Gamma}a_{31}\Delta t\frac{\lambda^{(n,1)}}{2}\left<\left<\rho^{n}\mathbf{u}^{n}\right>\right>:\left<\left<\mathbf{v}\right>\right>d\Sigma \nonumber \\
&-\quad\sum_{\Gamma\in\mathcal{E}}\int_{\Gamma}a_{32}\Delta t\frac{\lambda^{(n,2)}}{2}\left<\left<\rho^{(n,2)}\mathbf{u}^{(n,2)}\right>\right>:\left<\left<\mathbf{v}\right>\right>d\Sigma, \label{eq:momentum_stage_3}
\end{align}
where

\begin{align*}
\lambda^{(n,2)} &= \max\left(\left|\mathbf{u}^{(n,2)^{+}}\cdot\mathbf{n}^{+}\right|,\left|\mathbf{u}^{(n,2)^{-}}\cdot\mathbf{n}^{-}\right|\right)
\end{align*}
is employed for the upwind flux. Now, taking \(\mathbf{v} = \boldsymbol{\varphi}_i\) and exploiting the representation of \(\mathbf{u}^{(n,3)}\) and \(p^{(n,3)}\), we end up with the following relation

\begin{align}
&\qquad\sum_K\int_K\rho^{(n,3)}\sum_{j=1}^{\text{dim}(\mathbf{V}_{r})}u_j^{(n,3)}\boldsymbol{\varphi}_j\cdot\boldsymbol{\varphi}_id\Omega - \sum_K\int_K\tilde{a}_{33}\frac{\Delta t}{M^2}\sum_{j=1}^{\text{dim}(Q_{r})}p_j^{(n,3)}\psi_j\nabla\cdot\boldsymbol{\varphi}_id\Omega\nonumber \\ 
&+\quad\sum_{\Gamma\in\mathcal{E}}\int_{\Gamma}\tilde{a}_{33}\frac{\Delta t}{M^2}\sum_{j=1}^{\text{dim}(Q_{r})}p_j^{(n,3)}\left\{\left\{\psi_j\right\}\right\}\left[\left[\boldsymbol{\varphi}_i\right]\right]d\Sigma \nonumber \\
&=\quad\sum_K\int_K\rho^{n}\mathbf{u}^{n}\cdot\boldsymbol{\varphi}_id\Omega - \sum_K\int_Ka_{31}\frac{\Delta t}{Fr^2}\rho^{n}\mathbf{k}\cdot\boldsymbol{\varphi}_id\Omega - \sum_K\int_Ka_{32}\frac{\Delta t}{Fr^2}\rho^{(n,2)}\mathbf{k}\cdot\boldsymbol{\varphi}_id\Omega \nonumber \\
&+\quad\sum_K\int_Ka_{31}\Delta t\left(\rho^n\mathbf{u}^n\otimes\mathbf{u}^n\right):\nabla\boldsymbol{\varphi}_id\Omega + \sum_K\int_K\tilde{a}_{31}\frac{\Delta t}{M^2}p^n\nabla\cdot\boldsymbol{\varphi}_id\Omega \nonumber \\
&+\quad\sum_K\int_Ka_{32}\Delta t\left(\rho^{(n,2)}\mathbf{u}^{(n,2)}\otimes\mathbf{u}^{(n,2)}\right):\nabla\boldsymbol{\varphi}_id\Omega + \sum_K\int_K\tilde{a}_{32}\frac{\Delta t}{M^2}p^{(n,2)}\nabla\cdot\boldsymbol{\varphi}_id\Omega \nonumber \\
&-\quad\sum_{\Gamma\in\mathcal{E}}\int_{\Gamma}a_{31}\Delta t\left\{\left\{\rho^n\mathbf{u}^n\otimes\mathbf{u}^n\right\}\right\}:\left<\left<\boldsymbol{\varphi}_i\right>\right>d\Sigma \nonumber \\ 
&-\quad\sum_{\Gamma\in\mathcal{E}}\int_{\Gamma}\tilde{a}_{31}\frac{\Delta t}{M^2}\left\{\left\{p^n\right\}\right\}\left[\left[\boldsymbol{\varphi}_i\right]\right]d\Sigma \nonumber \\
&-\quad\sum_{\Gamma\in\mathcal{E}}\int_{\Gamma}a_{32}\Delta t\left\{\left\{\rho^{(n,2)}\mathbf{u}^{(n,2)}\otimes\mathbf{u}^{(n,2)}\right\}\right\}:\left<\left<\boldsymbol{\varphi}_i\right>\right>d\Sigma \nonumber \\
&-\quad\sum_{\Gamma\in\mathcal{E}}\int_{\Gamma}\tilde{a}_{32}\frac{\Delta t}{M^2}\left\{\left\{p^{(n,2)}\right\}\right\}\left[\left[\boldsymbol{\varphi}_i\right]\right]d\Sigma \nonumber \\
&-\quad\sum_{\Gamma\in\mathcal{E}}\int_{\Gamma}a_{31}\Delta t\frac{\lambda^{(n,1)}}{2}\left<\left<\rho^{n}\mathbf{u}^{n}\right>\right>:\left<\left<\boldsymbol{\varphi}_i\right>\right>d\Sigma \nonumber \\
&-\quad\sum_{\Gamma\in\mathcal{E}}\int_{\Gamma}a_{32}\Delta t\frac{\lambda^{(n,2)}}{2}\left<\left<\rho^{(n,2)}\mathbf{u}^{(n,2)}\right>\right>:\left<\left<\boldsymbol{\varphi}_i\right>\right>d\Sigma, 
\end{align}
which can be written in compact form as

\begin{equation}
\mathbf{A}^{(n,3)}\mathbf{U}^{(n,3)} + \mathbf{B}^{(n,3)}\mathbf{P}^{(n,3)} = \mathbf{F}^{(n,3)},
\end{equation}
where we have set

\begin{align}
A_{ij}^{(n,3)} &= \sum_K\int_K\rho^{(n,3)}\boldsymbol{\varphi}_j\cdot\boldsymbol{\varphi}_id\Omega \\
B_{ij}^{(n,3)} &= \sum_K\int_K-\tilde{a}_{33}\frac{\Delta t}{M^2}\nabla\cdot\boldsymbol{\varphi}_i\psi_jd\Omega \nonumber \\ &+\sum_{\Gamma\in\mathcal{E}}\int_{\Gamma}\tilde{a}_{33}\frac{\Delta t}{M^2}\left\{\left\{\psi_j\right\}\right\}\left[\left[\boldsymbol{\varphi}_i\right]\right]d\Sigma 
\end{align}
and \(\mathbf{U}^{(n,3)}\) denotes the vector of the degrees of freedom associated to the velocity field, whereas \(\mathbf{P}^{(n,3)}\) denotes the vector of the degrees of freedom associated to the pressure. 
Consider now the weak formulation for the energy equation in \eqref{eq:ns_comp_hyp_tdisc2}

\begin{align}
&\qquad\sum_K\int_K\rho^{(n,3)}E^{(n,3)}wd\Omega - \sum_K\int_K\tilde{a}_{33}\Delta th^{(n,3)}\rho^{(n,3)}\mathbf{u}^{(n,3)}\cdot\nabla wd\Omega \nonumber \\ &+\quad\sum_{\Gamma\in\mathcal{E}}\int_{\Gamma}\tilde{a}_{33}\Delta t\left\{\left\{h^{(n,3)}\rho^{(n,3)}\mathbf{u}^{(n,3)}\right\}\right\}\cdot\left[\left[w\right]\right]d\Sigma \nonumber \\
&=\quad\sum_K\int_K\rho^{n}E^{n}wd\Omega - \sum_K\int_Ka_{31}\frac{\Delta tM^2}{Fr^2}\rho^{n}\mathbf{k}\cdot\mathbf{u}^{n}wd\Omega - \sum_K\int_Ka_{32}\frac{\Delta tM^2}{Fr^2}\rho^{(n,2)}\mathbf{k}\cdot\mathbf{u}^{(n,2)}wd\Omega \nonumber \\
&+\quad\sum_K\int_Ka_{31}\Delta tM^2\left(k^{n}\rho^n\mathbf{u}^n\right)\cdot\nabla wd\Omega + \sum_K\int_K\tilde{a}_{31}\Delta t\left(h^{n}\rho^n\mathbf{u}^n\right)\cdot\nabla wd\Omega \nonumber \\
&+\quad\sum_K\int_Ka_{32}\Delta tM^2\left(k^{(n,2)}\rho^{(n,2)}\mathbf{u}^{(n,2)}\right)\cdot\nabla wd\Omega + \sum_K\int_K\tilde{a}_{32}\Delta t\left(h^{(n,2)}\rho^{(n,2)}\mathbf{u}^{(n,2)}\right)\cdot\nabla wd\Omega \nonumber \\
&-\quad\sum_{\Gamma\in\mathcal{E}}\int_{\Gamma}a_{31}\Delta tM^2\left\{\left\{k^n\rho^n\mathbf{u}^n\right\}\right\}\cdot\left[\left[w\right]\right]d\Sigma \nonumber \\
&-\quad\sum_{\Gamma\in\mathcal{E}}\int_{\Gamma}\tilde{a}_{31}\Delta t\left\{\left\{h^n\rho^n\mathbf{u}^n\right\}\right\}\cdot\left[\left[w\right]\right]d\Sigma \nonumber \\
&-\quad\sum_{\Gamma\in\mathcal{E}}\int_{\Gamma}a_{32}\Delta tM^2\left\{\left\{k^{(n,2)}\rho^{(n,2)}\mathbf{u}^{(n,2)}\right\}\right\}\cdot\left[\left[w\right]\right]d\Sigma \nonumber \\
&-\quad\sum_{\Gamma\in\mathcal{E}}\int_{\Gamma}\tilde{a}_{32}\Delta t\left\{\left\{h^{(n,2)}\rho^{(n,2)}\mathbf{u}^{(n,2)}\right\}\right\}\cdot\left[\left[w\right]\right]d\Sigma \nonumber \\
-&\quad\sum_{\Gamma\in\mathcal{E}}\int_{\Gamma}a_{31}\Delta t\frac{\lambda^{(n,1)}}{2}\left[\left[\rho^{n}E^{n}\right]\right] \cdot \left[\left[w\right]\right]d\Sigma \nonumber \\
-&\quad\sum_{\Gamma\in\mathcal{E}}\int_{\Gamma}a_{31}\Delta t\frac{\lambda^{(n,1)}}{2}\left[\left[\rho^{n}E^{n}\right]\right] \cdot \left[\left[w\right]\right]d\Sigma \nonumber \\
-&\quad\sum_{\Gamma\in\mathcal{E}}\int_{\Gamma}a_{32}\Delta t\frac{\lambda^{(n,2)}}{2}\left[\left[\rho^{(n,2)}E^{(n,2)}\right]\right] \cdot \left[\left[w\right]\right]d\Sigma. \label{eq:energy_stage_3} 
\end{align}
Take now \(w = \psi_i\) and consider the expansion for \(\mathbf{u}^{(n,3)}\)

\begin{align}
&\qquad\sum_K\int_K\rho^{(n,3)}E^{(n,3)}\psi_id\Omega - \sum_K\int_K\tilde{a}_{33}\Delta th^{(n,3)}\rho^{(n,3)}\sum_{j=1}^{\text{dim}(\mathbf{V}_{r})}u_j^{(n,3)}\boldsymbol{\varphi}_j\cdot\nabla \psi_id\Omega \nonumber \\ &+\quad\sum_{\Gamma\in\mathcal{E}^{I}}\int_{\Gamma}\tilde{a}_{33}\Delta t \sum_{j=1}^{\text{dim}(\mathbf{V}_{r})}u_j^{(n,3)}\left\{\left\{h^{(n,3)}\rho^{(n,3)}\boldsymbol{\varphi}_j\right\}\right\}\cdot\left[\left[\psi_i\right]\right]d\Sigma \nonumber \\
&=\quad\sum_K\int_K\rho^{n}E^{n}\psi_id\Omega - \sum_K\int_Ka_{31}\frac{\Delta tM^2}{Fr^2}\rho^{n}\mathbf{k}\cdot\mathbf{u}^{n}\psi_id\Omega - \sum_K\int_Ka_{32}\frac{\Delta tM^2}{Fr^2}\rho^{(n,2)}\mathbf{k}\cdot\mathbf{u}^{(n,2)}\psi_id\Omega \nonumber \\
&+\quad\sum_K\int_Ka_{31}\Delta tM^2\left(k^{n}\rho^n\mathbf{u}^n\right)\cdot\nabla \psi_id\Omega + \sum_K\int_K\tilde{a}_{31}\Delta t\left(h^{n}\rho^n\mathbf{u}^n\right)\cdot\nabla \psi_id\Omega \nonumber \\
&+\quad\sum_K\int_Ka_{32}\Delta tM^2\left(k^{(n,2)}\rho^{(n,2)}\mathbf{u}^{(n,2)}\right)\cdot\nabla \psi_id\Omega + \sum_K\int_K\tilde{a}_{32}\Delta t\left(h^{(n,2)}\rho^{(n,2)}\mathbf{u}^{(n,2)}\right)\cdot\nabla \psi_id\Omega \nonumber \\
&-\quad\sum_{\Gamma\in\mathcal{E}}\int_{\Gamma}a_{31}\Delta tM^2\left\{\left\{k^n\rho^n\mathbf{u}^n\right\}\right\}\cdot\left[\left[\psi_i\right]\right]d\Sigma \nonumber \\
&-\quad\sum_{\Gamma\in\mathcal{E}}\int_{\Gamma}\tilde{a}_{31}\Delta t\left\{\left\{h^n\rho^n\mathbf{u}^n\right\}\right\}\cdot\left[\left[\psi_i\right]\right]d\Sigma \nonumber \\
&-\quad\sum_{\Gamma\in\mathcal{E}}\int_{\Gamma}a_{32}\Delta tM^2\left\{\left\{k^{(n,2)}\rho^{(n,2)}\mathbf{u}^{(n,2)}\right\}\right\}\cdot\left[\left[\psi_i\right]\right]d\Sigma \nonumber \\
&-\quad\sum_{\Gamma\in\mathcal{E}}\int_{\Gamma}\tilde{a}_{32}\Delta t\left\{\left\{h^{(n,2)}\rho^{(n,2)}\mathbf{u}^{(n,2)}\right\}\right\}\cdot\left[\left[\psi_i\right]\right]d\Sigma \nonumber \\
-&\quad\sum_{\Gamma\in\mathcal{E}}\int_{\Gamma}a_{31}\Delta t\frac{\lambda^{(n,1)}}{2}\left[\left[\rho^{n}E^{n}\right]\right] \cdot \left[\left[\psi_{i}\right]\right]d\Sigma \nonumber \\
-&\quad\sum_{\Gamma\in\mathcal{E}}\int_{\Gamma}a_{32}\Delta t\frac{\lambda^{(n,2)}}{2}\left[\left[\rho^{(n,2)}E^{(n,2)}\right]\right] \cdot \left[\left[\psi_{i}\right]\right]d\Sigma, \label{eq:algebraic_energy_stage_3}
\end{align}
which can be expressed in compact form as

\begin{equation}
\mathbf{C}^{(n,3)}\mathbf{U}^{(n,3)} = \mathbf{G}^{(n,3)},
\end{equation}
where

\begin{align}
C_{ij}^{(n,3)} &= \sum_K\int_K-\tilde{a}_{33}\Delta th^{(n,3)}\rho^{(n,3)}\boldsymbol{\varphi}_j\cdot\nabla\psi_id\Omega \nonumber \\
&+\sum_{\Gamma\in\mathcal{E}}\int_{\Gamma}\tilde{a}_{33}\Delta t \left\{\left\{h^{(n,3)}\rho^{(n,3)}\boldsymbol{\varphi}_j\right\}\right\}\cdot\left[\left[\psi_i\right]\right]d\Sigma. 
\end{align}
Formally, one can derive \(\mathbf{U}^{(n,3)} = (\mathbf{A}^{{(n,3)}})^{-1}\left(\mathbf{F}^{(n,3)} - \mathbf{B}^{(n,3)}\mathbf{P}^{(n,3)}\right)\) and obtain the following relation 

\begin{equation}
\mathbf{C}^{(n,3)}(\mathbf{A}^{{(n,3)}})^{-1}\left(\mathbf{F}^{(n,3)} - \mathbf{B}^{(n,3)}\mathbf{P}^{(n,3)}\right) = \mathbf{G}^{(n,3)}.
\end{equation}
Taking into account that 

$$\rho^{(n,3)}E^{(n,3)} = \rho^{(n,3)}e^{(n,3)} + M^2\rho^{(n,3)}k^{(n,3)},$$
we obtain

\begin{equation}
\mathbf{C}^{(n,3)}(\mathbf{A}^{{(n,3)}})^{-1}\left(\mathbf{F}^{(n,3)} - \mathbf{B}^{(n,3)}\mathbf{P}^{(n,3)}\right) = -\mathbf{D}^{(n,3)}\mathbf{P}^{(n,3)} +  \tilde{\mathbf{G}}^{(n,3)},
\end{equation}
where 

\begin{align}
D_{ij}^{(n,3)} = &\sum_K\int_K\rho^{(n,3)}e^{(n,3)}(\psi_j)\psi_id\Omega
\end{align}
and \(\tilde{\mathbf{G}}^{(n,3)}\) takes into account all the other terms (the one at previous stage and the kinetic energy). Again, the above system is solved by a fixed point procedure. More specifically, setting \(\mathbf{P}^{(n,3,0)} = \mathbf{P}^{(n,3)}, k^{(n,3,0)} = k^{(n,2)}\) for \(l = 1,\dots, M\) one solves the equation

\begin{equation}
\left(\mathbf{D}^{(n,3,l)} - \mathbf{C}^{(n,3,l)}(\mathbf{A}^{{(n,3)}})^{-1}\mathbf{B}^{(n,3)}\right)\mathbf{P}^{(n,3,l+1)} = 
\tilde{\mathbf{G}}^{(n,3,l)} - \mathbf{C}^{(n,3,l)}(\mathbf{A}^{{(n,3)}})^{-1}\mathbf{F}^{(n,2,l)}
\nonumber 
\end{equation}
and then updates the velocity solving

\begin{equation}
\mathbf{A}^{(n,3)}\mathbf{U}^{(n,3,l+1)} = \mathbf{F}^{(n,3,l)} - \mathbf{B}^{(n,3)}\mathbf{P}^{(n,3,l+1)}. \nonumber
\end{equation}
Once the iterations have been completed, one sets
$ \mathbf{u}^{(n,3)}=\mathbf{u}^{(n,3,M+1)} $  and $ E^{(n,3)} $ accordingly.
One sets then

$$ \rho^{n+1} = \rho^{(n,3)}  \ \ \ \  \mathbf{\tilde{u}}^{(n,1)} = \mathbf{u}^{(n,3)}  \ \ \ \  \tilde{E}^{(n,1)} = E^{(n,3)} $$
and proceeds to the implicit discretization of the viscous terms, which is carried out by the implicit part of the IMEX method described above.
We would like to stress that the method outlined above does not require to introduce reference solutions, does not introduce inconsistencies in the splitting and only requires the solution of linear systems of a size equal to that of the number of discrete degrees of freedom needed to describe a scalar variable, as in \cite{dumbser:2016}. This contrasts with other low Mach approaches based on IMEX methods, such as e.g. the technique proposed for the Euler equations in \cite{zeifang:2019}.

\section{Implementation issues}
\label{sec:implement} \indent

As stated in Section \ref{sec:intro}, the proposed method has been implemented using the numerical library \textit{deal.II}, which is based on a matrix-free approach. As a consequence, no global sparse matrix is built and only the action of the linear operators defined above on a vector is actually computed. The implementation follows the operator splitting strategy previously described. Therefore, we have built two \textit{ad-hoc} structures, one for the hyperbolic part and one for the diffusive part. Another feature of the library employed during the numerical simulations is the mesh adaptation capability, as we will see in the presentation of the results. The preconditioned conjugate gradient method implemented in the function \textit{SolverCG} of the library was employed to solve the linear systems for the density and for the update of the velocity in the fixed point iterations, as well as to solve the linear systems associated to the matrices \(\mathbf{A}^{(n,2)}\) and \(\mathbf{A}^{(n,3)}\), while the GMRES solver implemented in the function \textit{SolverGMRES} of the same library was used for the remaining linear systems. 

On the other hand, in the diffusive part, a Symmetric Interior Penalty (SIP) approach has been adopted for the space discretization \cite{arnold:1982, orlando:2021}. Moreover, following \cite{fehn:2019}, we set for each face \(\Gamma\) of a cell \(K\)

\[\sigma_{\Gamma,K} = \left(r+1\right)^2\frac{\text{diam}\left(\Gamma\right)}{\text{diam}\left(K\right)}\]
and we define the penalization constant of the SIP method as

\[\overline{C} = \frac{1}{2}\left(\sigma_{\Gamma,K^{+}} + \sigma_{\Gamma,K^{-}}\right)\]
if \(\Gamma \in \mathcal{E}^{I}\) and \(\overline{C} = \sigma_{\Gamma,K}\) if \(\Gamma \in \mathcal{E}^{B}\). All the linear systems for this part are solved using the preconditioned conjugated gradient mentioned above. A geometric multigrid preconditioner has been employed for the solution of the symmetric linear systems using the procedure described in \cite{janssen:2011}, whereas a Jacobi preconditioner has been used for the non symmetric ones. 

\section{Numerical tests}
\label{sec:tests}

The numerical scheme outlined in the previous Sections has been validated in a number of  benchmarks. We set $ \cal H =\min\{\mathrm{diam}(K) |  K\in\mathcal{T}_h \} $ and we define two Courant numbers, one based on the speed of sound denoted by \(C\), the so-called acoustic Courant number, and one based on the local velocity of the flow, the so-called advective Courant number, denoted by \(C_u\):

\begin{equation}
\label{eq:Courant}
C = \frac{1}{M}rc \Delta t/{\cal H}, \ \ \ \ \    C_u = ru \Delta t/{\cal H},
\end{equation}
where \(c\) is the magnitude of the speed of sound and \(u\) is the magnitude of the flow velocity. The factor \(\frac{1}{M}\) is due to the scaling of the speed of sound, as reported in \cite{munz:2003}, for an ideal gas, and proven in \ref{sec:eigenvalues} in the monodimensional case for a general equation of state. For the tests using the ideal gas law, the value \(\gamma = 1.4\) for the specific heat ratio is employed, unless differently stated. 

\subsection{Isentropic vortex}
\label{ssec:vortex}

As a first benchmark, we consider for an ideal gas the two dimensional inviscid isentropic vortex also studied in \cite{tavelli:2017, zeifang:2019}. For this test, an analytic solution is available, that can be used to assess the convergence properties of the scheme. The initial conditions are given as a perturbation of a reference state

\[\rho(\mathbf{x},0) = \rho_\infty + \delta\rho \qquad \mathbf{u}(\mathbf{x},0) = \mathbf{u}_{\infty} + \delta\mathbf{u} \qquad p(\mathbf{x},0) = p_{\infty} + \delta p\]
and the exact solution is a propagation of the initial condition at the background velocity

\[\rho(\mathbf{x},t) = \rho(\mathbf{x} - \mathbf{u}_{\infty}t,0) \qquad \mathbf{u}(\mathbf{x},t) = \mathbf{u}(\mathbf{x} - \mathbf{u}_{\infty}t,0) \qquad p(\mathbf{x},t) = p(\mathbf{x} - \mathbf{u}_{\infty}t,0).\]
The typical perturbation is defined as

\begin{equation}
\widetilde{\delta T} = \frac{1-\gamma}{8\gamma\pi^2}\beta^2e^{1 - \tilde r^2},
\end{equation} 
with \(\tilde r^2 = \left(x-x_0\right)^2 + \left(y-y_0\right)^2\) denoting the radial coordinate and \(\beta\) being the vortex strength. As explained in \cite{zeifang:2019}, however, in order to emphasize the role of the Mach number \(M\), we define 

\begin{equation}
\delta T = \frac{1-\gamma}{8\gamma\pi^2}M^2\beta^2e^{1 - \tilde r^2}
\end{equation} 
and we set

\begin{equation}
\rho(\mathbf{x},0) = \left(1 + \delta T\right)^{\frac{1}{\gamma-1}} 
\ \ \ \
p(\mathbf{x},0) = M^2\left(1 + \delta T\right)^{\frac{\gamma}{\gamma-1}}.
\end{equation}
For what concerns the velocity the typical perturbation is defined as

\begin{equation}
\widetilde{\delta\mathbf{u}} = \beta\begin{pmatrix}
-y \\
x
\end{pmatrix}\frac{e^{\frac{1}{2}\left(1 - \tilde r^2\right)}}{2\pi}
\end{equation}
and also in this case we rescale it using \(M\) 
\begin{equation}
\delta\mathbf{u} = \beta M\begin{pmatrix}
-y \\
x
\end{pmatrix}\frac{e^{\frac{1}{2}\left(1 - \tilde r^2\right)}}{2\pi}.
\end{equation}
We apply the same reasoning also to the background velocity and therefore we define \(\mathbf{u}_\infty = M\tilde{\mathbf{u}}_\infty\) with
\(\tilde{\mathbf{u}}_\infty = [10, 10]^T\). To avoid problems related to the definition of boundary conditions, we choose a sufficiently large domain \(\Omega = \left(-10,10\right)^2\) and periodic boundary conditions and we set \(\rho_\infty = 1\), \(p_\infty = 1\), \(x_0 = y_0 = 0\), \(\beta = 10\), the final time \(T_f=1\) and \(M = 0.1\). Notice that we refrain from investigating the properties of the method in the very low Mach number regime
for this test, since this entails an almost constant solution. The numerical experiments have been carried out on Cartesian meshes of square elements with $N_{el} $ elements
in each coordinate direction, choosing for each spatial resolution time steps so that the Courant numbers remained constant (hyperbolic scaling). 

We first consider the original  IMEX-ARK scheme with $\alpha =\frac{7 - 2\gamma}{6} $ for the explicit part.
In Tables \ref{tab:vortex_k1_C_01}-\ref{tab:vortex_k1_C_15}, the  \(L^2\) errors for  density, velocity and pressure are reported at various resolutions for the $r=1 $ case, while the corresponding errors in the $r=2 $ case are given in Tables \ref{tab:vortex_k2_C_01}-\ref{tab:vortex_k2_C_15}. We observe that, in general, convergence rates of at least \(r + \frac{1}{2}\) are observed for \(r = 1\) while for \(r = 2\) the convergence rate seem to degrade at the finest resolution as soon as the Courant number grows, due to increasing influence of the dominant second order time discretization error. 

\begin{table}[H]
	\centering
	\begin{tabular}{|c|c|c|c|c|c|c|}
		\hline
		$N_{el} $ & \(L^2\) rel. error \(\rho\) & \(L^2\) rate \(\rho\) & \(L^2\) rel. error \(\mathbf{u}\) & \(L^2\) rate \(\mathbf{u}\) & \(L^2\) rel. error \(p\) & \(L^2\) rate \(p\) \\
		\hline
		\(10\) & \(1.99 \cdot 10^{-3}\) &  & \(1.19 \cdot 10^{-2}\) & & \(2.79 \cdot 10^{-3}\) & \\
		\hline
		\(20\) & \(7.87 \cdot 10^{-4}\) & 1.34 & \(3.87 \cdot 10^{-3}\) & 1.62 & \(1.11 \cdot 10^{-3}\) & 1.33 \\
		\hline
		\(40\) & \(2.56 \cdot 10^{-4}\) & 1.62 & \(1.08 \cdot 10^{-3}\) & 1.84 & \(3.62 \cdot 10^{-4}\) & 1.62 \\
		\hline
		\(80\) & \(7.22 \cdot 10^{-5}\) & 1.83 & \(2.73 \cdot 10^{-4}\) & 1.98 & \(1.01 \cdot 10^{-4}\) & 1.84 \\
		\hline 
	\end{tabular}
	\caption{Convergence test for the inviscid isentropic vortex at \(C \approx 0.01\), \(C_u \approx 0.01\) with  \(r = 1\) and $\alpha =\frac{7 - 2\gamma}{6} $ for the explicit part. Relative errors for the density, the velocity and the pressure in $L^2$ norm. \(N_{el}\) denotes the number of elements along each direction.}
	\label{tab:vortex_k1_C_01}
\end{table}

\begin{table}[H]
	\centering
	\begin{tabular}{|c|c|c|c|c|c|c|}
		\hline
		$N_{el} $ & \(L^2\) rel. error \(\rho\) & \(L^2\) rate \(\rho\) & \(L^2\) rel. error \(\mathbf{u}\) & \(L^2\) rate \(\mathbf{u}\) & \(L^2\) rel. error \(p\) & \(L^2\) rate \(p\) \\
		\hline
		\(10\) & \(1.99 \cdot 10^{-3}\) &  & \(1.20 \cdot 10^{-2}\) & & \(2.77 \cdot 10^{-3}\) & \\
		\hline
		\(20\) & \(7.93 \cdot 10^{-4}\) & 1.33 & \(3.95 \cdot 10^{-3}\) & 1.61 & \(1.11 \cdot 10^{-3}\) & 1.32 \\
		\hline
		\(40\) & \(2.61 \cdot 10^{-4}\) & 1.60 & \(1.13 \cdot 10^{-3}\) & 1.81 & \(3.65 \cdot 10^{-4}\) & 1.60 \\
		\hline
		\(80\) & \(7.50 \cdot 10^{-5}\) & 1.80 & \(2.98 \cdot 10^{-4}\) & 1.92 & \(1.03 \cdot 10^{-4}\) & 1.83 \\
		\hline 
	\end{tabular}
	\caption{Convergence test for the inviscid isentropic vortex at \(C \approx 0.05\), \(C_u \approx 0.05\) with  \(r = 1\) and $\alpha =\frac{7 - 2\gamma}{6} $ for the explicit part. Relative errors for the density, the velocity and the pressure in $L^2$ norm. \(N_{el}\) denotes the number of elements along each direction.}
	\label{tab:vortex_k1_C_05}
\end{table}

\begin{table}[H]
	\centering
	\begin{tabular}{|c|c|c|c|c|c|c|}
		\hline
		$N_{el} $ & \(L^2\) rel. error \(\rho\) & \(L^2\) rate \(\rho\) & \(L^2\) rel. error \(\mathbf{u}\) & \(L^2\) rate \(\mathbf{u}\) & \(L^2\) rel. error \(p\) & \(L^2\) rate \(p\) \\
		\hline
		\(10\) & \(2.11 \cdot 10^{-3}\) &  & \(1.20 \cdot 10^{-2}\) & & \(2.72 \cdot 10^{-3}\) & \\
		\hline
		\(20\) & \(8.01 \cdot 10^{-4}\) & 1.40 & \(4.09 \cdot 10^{-3}\) & 1.56 & \(1.09 \cdot 10^{-3}\) & 1.32 \\
		\hline
		\(40\) & \(2.67 \cdot 10^{-4}\) & 1.58 & \(1.22 \cdot 10^{-3}\) & 1.75 & \(3.63 \cdot 10^{-4}\) & 1.59 \\
		\hline
		\(80\) & \(7.99 \cdot 10^{-5}\) & 1.74 & \(3.44 \cdot 10^{-4}\) & 1.83 & \(1.06 \cdot 10^{-4}\) & 1.78 \\
		\hline 
	\end{tabular}
	\caption{Convergence test for the inviscid isentropic vortex at \(C \approx 0.15\), \(C_u \approx 0.14\) with  \(r = 1\) and $\alpha =\frac{7 - 2\gamma}{6} $ for the explicit part. Relative errors for the density, the velocity and the pressure in $L^2$ norm. \(N_{el}\) denotes the number of elements along each direction.}
	\label{tab:vortex_k1_C_15}
\end{table}

\begin{table}[H]
	\centering
	\begin{tabular}{|c|c|c|c|c|c|c|}
		\hline
		$N_{el} $ & \(L^2\) rel. error \(\rho\) & \(L^2\) rate \(\rho\) & \(L^2\) rel. error \(\mathbf{u}\) & \(L^2\) rate \(\mathbf{u}\) & \(L^2\) rel. error \(p\) & \(L^2\) rate \(p\) \\
		\hline
		\(10\) & \(6.37 \cdot 10^{-4}\) &  & \(2.61 \cdot 10^{-3}\) & & \(9.09 \cdot 10^{-4}\) & \\
		\hline
		\(20\) & \(1.18 \cdot 10^{-4}\) & 2.43 & \(3.59 \cdot 10^{-4}\) & 2.86 & \(1.64 \cdot 10^{-4}\) & 2.47 \\
		\hline
		\(40\) & \(1.83 \cdot 10^{-5}\) & 2.69 & \(4.39 \cdot 10^{-5}\) & 3.03 & \(2.55 \cdot 10^{-5}\) & 2.69 \\
		\hline
		\(80\) & \(3.08 \cdot 10^{-6}\) & 2.57 & \(6.96 \cdot 10^{-6}\) & 2.66 & \(4.21 \cdot 10^{-6}\) & 2.60 \\
		\hline 
	\end{tabular}
	\caption{Convergence test for the inviscid isentropic vortex at \(C \approx 0.01\), \(C_u \approx 0.01\) with  \(r = 2\) and $\alpha =\frac{7 - 2\gamma}{6} $ for the explicit part. Relative errors for the density, the velocity and the pressure in $L^2$ norm. \(N_{el}\) denotes the number of elements along each direction.}
	\label{tab:vortex_k2_C_01}
\end{table}

\begin{table}[H]
	\centering
	\begin{tabular}{|c|c|c|c|c|c|c|}
		\hline
		$N_{el} $ & \(L^2\) rel. error \(\rho\) & \(L^2\) rate \(\rho\) & \(L^2\) rel. error \(\mathbf{u}\) & \(L^2\) rate \(\mathbf{u}\) & \(L^2\) rel. error \(p\) & \(L^2\) rate \(p\) \\
		\hline
		\(10\) & \(6.33 \cdot 10^{-4}\) &  & \(2.65 \cdot 10^{-3}\) & & \(9.08 \cdot 10^{-4}\) & \\
		\hline
		\(20\) & \(1.19 \cdot 10^{-4}\) & 2.41 & \(3.89 \cdot 10^{-4}\) & 2.77 & \(1.65 \cdot 10^{-4}\) & 2.46 \\
		\hline
		\(40\) & \(2.00 \cdot 10^{-5}\) & 2.57 & \(6.45 \cdot 10^{-5}\) & 2.59 & \(2.67 \cdot 10^{-5}\) & 2.63 \\
		\hline
		\(80\) & \(4.87 \cdot 10^{-6}\) & 2.04 & \(2.22 \cdot 10^{-5}\) & 1.54 & \(5.48 \cdot 10^{-6}\) & 2.28 \\
		\hline 
	\end{tabular}
	\caption{Convergence test for the inviscid isentropic vortex at \(C \approx 0.05\), \(C_u \approx 0.05\) with  \(r = 2\) and $\alpha =\frac{7 - 2\gamma}{6} $ for the explicit part. Relative errors for the density, the velocity and the pressure in $L^2$ norm. \(N_{el}\) denotes the number of elements along each direction.}
	\label{tab:vortex_k2_C_05}
\end{table}

\begin{table}[H]
	\centering
	\begin{tabular}{|c|c|c|c|c|c|c|}
		\hline
		$N_{el} $ & \(L^2\) rel. error \(\rho\) & \(L^2\) rate \(\rho\) & \(L^2\) rel. error \(\mathbf{u}\) & \(L^2\) rate \(\mathbf{u}\) & \(L^2\) rel. error \(p\) & \(L^2\) rate \(p\) \\
		\hline
		\(10\) & \(6.24 \cdot 10^{-4}\) &  & \(2.75 \cdot 10^{-3}\) & & \(8.96 \cdot 10^{-4}\) & \\
		\hline
		\(20\) & \(1.24 \cdot 10^{-4}\) & 2.33 & \(4.72 \cdot 10^{-4}\) & 2.44 & \(1.65 \cdot 10^{-4}\) & 2.44 \\
		\hline
		\(40\) & \(2.60 \cdot 10^{-5}\) & 2.25 & \(1.22 \cdot 10^{-4}\) & 1.95 & \(3.07 \cdot 10^{-5}\) & 2.43 \\
		\hline
		\(80\) & \(9.71 \cdot 10^{-6}\) & 1.42 & \(5.37 \cdot 10^{-5}\) & 1.18 & \(9.42 \cdot 10^{-6}\) & 1.70 \\
		\hline 
	\end{tabular}
	\caption{Convergence test for the inviscid isentropic vortex at \(C \approx 0.15\), \(C_u \approx 0.14\) with  \(r = 2\) and $\alpha =\frac{7 - 2\gamma}{6} $ for the explicit part. Relative errors for the density, the velocity and the pressure in $L^2$ norm. \(N_{el}\) denotes the number of elements along each direction.}
	\label{tab:vortex_k2_C_15}
\end{table}

Analogous results are shown in Tables \ref{tab:vortex_k1_alpha05_C_01}-\ref{tab:vortex_k2_alpha05_C_15}
for the modified scheme with \(\alpha = 0.5\), chosen, as discussed in Appendix \ref{sec:mono_analysis}, in order to increase the region of absolute monotonicity without affecting too much  stability. It can be seen that, for $r=1, $ slightly lower errors are obtained, especially for the density, while the behaviour for \(r = 2\) is similar to that of the original scheme. 

\begin{table}[H]
	\centering
	\begin{tabular}{|c|c|c|c|c|c|c|}
		\hline
		$N_{el} $ & \(L^2\) rel. error \(\rho\) & \(L^2\) rate \(\rho\) & \(L^2\) rel. error \(\mathbf{u}\) & \(L^2\) rate \(\mathbf{u}\) & \(L^2\) rel. error \(p\) & \(L^2\) rate \(p\) \\
		\hline
		\(10\) & \(2.02 \cdot 10^{-3}\) &  & \(1.19 \cdot 10^{-2}\) & & \(2.80 \cdot 10^{-3}\) & \\
		\hline
		\(20\) & \(7.82 \cdot 10^{-4}\) & 1.37 & \(3.81 \cdot 10^{-3}\) & 1.64 & \(1.11 \cdot 10^{-3}\) & 1.33 \\
		\hline
		\(40\) & \(2.51 \cdot 10^{-4}\) & 1.64 & \(1.04 \cdot 10^{-3}\) & 1.87 & \(3.58 \cdot 10^{-4}\) & 1.63 \\
		\hline
		\(80\) & \(6.96\cdot10^{-5}\) & 1.85 & \(2.50\cdot10^{-4}\) & 2.06 & \(9.89\cdot10^{-5}\) & 1.86 \\
		\hline 
	\end{tabular}
	\caption{Convergence test for the inviscid isentropic vortex at \(C \approx 0.01\), \(C_u \approx 0.01\) with  \(r = 1\) and $\alpha = 0.5 $ for the explicit part. Relative errors for the density, the velocity and the pressure in $L^2$ norm. \(N_{el}\) denotes the number of elements along each direction.}
	\label{tab:vortex_k1_alpha05_C_01}
\end{table}

\begin{table}[H]
	\centering
	\begin{tabular}{|c|c|c|c|c|c|c|}
		\hline
		$N_{el} $ & \(L^2\) rel. error \(\rho\) & \(L^2\) rate \(\rho\) & \(L^2\) rel. error \(\mathbf{u}\) & \(L^2\) rate \(\mathbf{u}\) & \(L^2\) rel. error \(p\) & \(L^2\) rate \(p\) \\
		\hline
		\(10\) & \(2.11 \cdot 10^{-3}\) &  & \(1.19 \cdot 10^{-2}\) & & \(2.82 \cdot 10^{-3}\) & \\
		\hline
		\(20\) & \(7.66 \cdot 10^{-4}\) & 1.46 & \(3.61 \cdot 10^{-3}\) & 1.72 & \(1.10 \cdot 10^{-3}\) & 1.36 \\
		\hline
		\(40\) & \(2.38 \cdot 10^{-4}\) & 1.69 & \(9.16 \cdot 10^{-4}\) & 1.98 & \(3.46 \cdot 10^{-4}\) & 1.67 \\
		\hline
		\(80\) & \(6.43 \cdot 10^{-5}\) & 1.89 & \(1.98 \cdot 10^{-4}\) & 2.21 & \(9.34 \cdot 10^{-5}\) & 1.89 \\
		\hline 
	\end{tabular}
	\caption{Convergence test for the inviscid isentropic vortex at \(C \approx 0.05\), \(C_u \approx 0.05\) with  \(r = 1\) and $\alpha = 0.5 $ for the explicit part. Relative errors for the density, the velocity and the pressure in $L^2$ norm. \(N_{el}\) denotes the number of elements along each direction.}
	\label{tab:vortex_k1_alpha05_C_05}
\end{table}

\begin{table}[H]
	\centering
	\begin{tabular}{|c|c|c|c|c|c|c|}
		\hline
		$N_{el} $ & \(L^2\) rel. error \(\rho\) & \(L^2\) rate \(\rho\) & \(L^2\) rel. error \(\mathbf{u}\) & \(L^2\) rate \(\mathbf{u}\) & \(L^2\) rel. error \(p\) & \(L^2\) rate \(p\) \\
		\hline
		\(10\) & \(2.38 \cdot 10^{-3}\) &  & \(1.17 \cdot 10^{-2}\) & & \(2.89 \cdot 10^{-3}\) & \\
		\hline
		\(20\) & \(7.49 \cdot 10^{-4}\) & 1.67 & \(3.25 \cdot 10^{-3}\) & 1.85 & \(1.09 \cdot 10^{-3}\) & 1.41 \\
		\hline
		\(40\) & \(2.26 \cdot 10^{-4}\) & 1.73 & \(7.60 \cdot 10^{-4}\) & 2.10 & \(3.31 \cdot 10^{-4}\) & 1.72 \\
		\hline
		\(80\) & \(6.85 \cdot 10^{-5}\) & 1.72 & \(2.26 \cdot 10^{-4}\) & 1.75 & \(9.12 \cdot 10^{-5}\) & 1.86 \\
		\hline 
	\end{tabular}
	\caption{Convergence test for the inviscid isentropic vortex at \(C \approx 0.15\), \(C_u \approx 0.14\) with  \(r = 1\) and $\alpha = 0.5 $ for the explicit part. Relative errors for the density, the velocity and the pressure in $L^2$ norm. \(N_{el}\) denotes the number of elements along each direction.}
	\label{tab:vortex_k1_alpha05_C_15}
\end{table}

\begin{table}[H]
	\centering
	\begin{tabular}{|c|c|c|c|c|c|c|}
		\hline
		$N_{el} $ & \(L^2\) rel. error \(\rho\) & \(L^2\) rate \(\rho\) & \(L^2\) rel. error \(\mathbf{u}\) & \(L^2\) rate \(\mathbf{u}\) & \(L^2\) rel. error \(p\) & \(L^2\) rate \(p\) \\
		\hline
		\(10\) & \(6.41 \cdot 10^{-4}\) &  & \(2.59 \cdot 10^{-3}\) & & \(9.07 \cdot 10^{-4}\) & \\
		\hline
		\(20\) & \(1.18 \cdot 10^{-4}\) & 2.44 & \(3.41 \cdot 10^{-4}\) & 2.93 & \(1.65 \cdot 10^{-4}\) & 2.43 \\
		\hline
		\(40\) & \(1.84 \cdot 10^{-5}\) & 2.68 & \(4.46 \cdot 10^{-5}\) & 2.94 & \(2.54 \cdot 10^{-5}\) & 2.50 \\
		\hline
		\(80\) & \(3.72 \cdot 10^{-6}\) & 2.31 & \(1.33 \cdot 10^{-5}\) & 1.75 & \(4.61 \cdot 10^{-6}\) & 2.46 \\
		\hline 
	\end{tabular}
	\caption{Convergence test for the inviscid isentropic vortex at \(C \approx 0.01\), \(C_u \approx 0.01\) with  \(r = 2\) and $ \alpha = 0.5 $ for the explicit part. Relative errors for the density, the velocity and the pressure in $L^2$ norm. \(N_{el}\) denotes the number of elements along each direction.}
	\label{tab:vortex_k2_alpha05_C_01}
\end{table}

\begin{table}[H]
	\centering
	\begin{tabular}{|c|c|c|c|c|c|c|}
		\hline
		$N_{el} $ & \(L^2\) rel. error \(\rho\) & \(L^2\) rate \(\rho\) & \(L^2\) rel. error \(\mathbf{u}\) & \(L^2\) rate \(\mathbf{u}\) & \(L^2\) rel. error \(p\) & \(L^2\) rate \(p\) \\
		\hline
		\(10\) & \(6.59 \cdot 10^{-4}\) &  & \(2.55 \cdot 10^{-3}\) & & \(9.05 \cdot 10^{-4}\) & \\
		\hline
		\(20\) & \(1.27 \cdot 10^{-4}\) & 2.38 & \(3.66 \cdot 10^{-4}\) & 2.80 & \(1.70 \cdot 10^{-4}\) & 2.41 \\
		\hline
		\(40\) & \(2.88 \cdot 10^{-5}\) & 2.14 & \(1.27 \cdot 10^{-4}\) & 1.53 & \(3.22 \cdot 10^{-5}\) & 2.40 \\
		\hline
		\(80\) & \(1.20 \cdot 10^{-5}\) & 1.26 & \(6.35 \cdot 10^{-5}\) & 1.00 & \(1.15 \cdot 10^{-5}\) & 1.49 \\
		\hline 
	\end{tabular}
	\caption{Convergence test for the inviscid isentropic vortex at \(C \approx 0.05\), \(C_u \approx 0.05\) with  \(r = 2\) and $ \alpha = 0.5 $ for the explicit part. Relative errors for the density, the velocity and the pressure in $L^2$ norm. \(N_{el}\) denotes the number of elements along each direction.}
	\label{tab:vortex_k2_alpha05_C_05}
\end{table}

\begin{table}[H]
	\centering
	\begin{tabular}{|c|c|c|c|c|c|c|}
		\hline
		$N_{el} $ & \(L^2\) rel. error \(\rho\) & \(L^2\) rate \(\rho\) & \(L^2\) rel. error \(\mathbf{u}\) & \(L^2\) rate \(\mathbf{u}\) & \(L^2\) rel. error \(p\) & \(L^2\) rate \(p\) \\
		\hline
		\(10\) & \(7.18 \cdot 10^{-4}\) &  & \(2.63 \cdot 10^{-3}\) & & \(9.17 \cdot 10^{-4}\) & \\
		\hline
		\(20\) & \(1.73 \cdot 10^{-4}\) & 2.05 & \(6.40 \cdot 10^{-4}\) & 2.04 & \(2.01 \cdot 10^{-4}\) & 2.19 \\
		\hline
		\(40\) & \(6.17 \cdot 10^{-5}\) & 1.49 & \(3.19 \cdot 10^{-4}\) & 1.00 & \(5.98 \cdot 10^{-5}\) & 1.75 \\
		\hline
		\(80\) & \(2.97 \cdot 10^{-5}\) & 1.05 & \(1.60 \cdot 10^{-4}\) & 1.00 & \(2.77 \cdot 10^{-5}\) & 1.11 \\
		\hline 
	\end{tabular}
	\caption{Convergence test for the inviscid isentropic vortex at \(C \approx 0.15\), \(C_u \approx 0.14\) with  \(r = 2\) and $ \alpha = 0.5 $ for the explicit part. Relative errors for the density, the velocity and the pressure in $L^2$ norm. \(N_{el}\) denotes the number of elements along each direction.}
	\label{tab:vortex_k2_alpha05_C_15}
\end{table}

In further numerical experiments, we have observed that the lack of absolute monotonicity strongly affects the computation of density and, as a consequence, the stability of the whole numerical scheme. For Courant number around $ C\approx 0.2 $ the original method becomes unstable, while the modified scheme with \(\alpha = 0.5\) is still able to recover the expected convergence rates at least in the $r=1 $ case, as evident from Table \ref{tab:vortex_k1_alfa0,5_C0,2}, while again in the \(r=2\) reported in Table \ref{tab:vortex_k2_alfa0,5_C0,2} we observe a degradation of the convergence rates. In order to be able to run at slightly longer time steps we have then chosen to use the $ \alpha = 0.5 $ value for the IMEX scheme for the rest of the numerical simulations carried out in this paper. We notice also that, for both schemes, the results compare well with the analogous results presented in \cite{tavelli:2017} and with those obtained in \cite{zeifang:2019} with a higher order IMEX method.

\begin{table}[H]
	\centering
	\begin{tabular}{|c|c|c|c|c|c|c|}
		\hline
		$N_{el} $ & \(L^2\) rel. error \(\rho\) & \(L^2\) rate \(\rho\) & \(L^2\) rel. error \(\mathbf{u}\) & \(L^2\) rate \(\mathbf{u}\) & \(L^2\) rel. error \(p\) & \(L^2\) rate \(p\) \\
		\hline
		\(10\) & \(2.71 \cdot 10^{-3}\) &  & \(1.16 \cdot 10^{-2}\) & & \(2.96 \cdot 10^{-3}\) & \\
		\hline
		\(20\) & \(7.74 \cdot 10^{-4}\) & 1.81 & \(2.95 \cdot 10^{-3}\) & 1.98 & \(1.09 \cdot 10^{-3}\) & 1.44 \\
		\hline
		\(40\) & \(2.34 \cdot 10^{-4}\) & 1.73 & \(7.71 \cdot 10^{-4}\) & 1.94 & \(3.28 \cdot 10^{-4}\) & 1.73 \\
		\hline
		\(80\) & \(8.91 \cdot 10^{-5}\) & 1.39 & \(3.74 \cdot 10^{-4}\) & 1.04 & \(1.01 \cdot 10^{-4}\) & 1.70 \\
		\hline 
	\end{tabular}
	\caption{Convergence test for the inviscid isentropic vortex at \(C \approx 0.2\), \(C_u \approx 0.2\) with  \(r = 1\) and $ \alpha = 0.5 $ for the explicit part. Relative errors for the density, the velocity and the pressure in $L^2$ norm. \(N_{el}\) denotes the number of elements along each direction.}
	\label{tab:vortex_k1_alfa0,5_C0,2}
\end{table}

\begin{table}[H]
	\centering
	\begin{tabular}{|c|c|c|c|c|c|c|}
		\hline
		$N_{el} $ & \(L^2\) rel. error \(\rho\) & \(L^2\) rate \(\rho\) & \(L^2\) rel. error \(\mathbf{u}\) & \(L^2\) rate \(\mathbf{u}\) & \(L^2\) rel. error \(p\) & \(L^2\) rate \(p\) \\
		\hline
		\(10\) & \(8.18 \cdot 10^{-4}\) &  & \(2.92 \cdot 10^{-3}\) & & \(9.58 \cdot 10^{-4}\) & \\
		\hline
		\(20\) & \(2.42 \cdot 10^{-4}\) & 1.76 & \(1.02 \cdot 10^{-3}\) & 1.52 & \(2.57 \cdot 10^{-4}\) & 1.90 \\
		\hline
		\(40\) & \(9.88 \cdot 10^{-5}\) & 1.29 & \(5.20 \cdot 10^{-4}\) & 0.97 & \(9.36\cdot10^{-5}\) & 1.46 \\
		\hline
		\(80\) & \(4.79 \cdot 10^{-5}\) & 1.04 & \(2.58 \cdot 10^{-4}\) & 1.01 & \(4.47 \cdot 10^{-5}\) & 1.07 \\
		\hline 
	\end{tabular}
	\caption{Convergence test for the inviscid isentropic vortex at \(C \approx 0.2\), \(C_u \approx 0.2\) with \(r = 2\) and $ \alpha = 0.5 $ for the explicit part. Relative errors for the density, the velocity and the pressure in $L^2$ norm. \(N_{el}\) denotes the number of elements along each direction.}
	\label{tab:vortex_k2_alfa0,5_C0,2}
\end{table}

\begin{figure}[H]
	\includegraphics[width=0.4\textwidth]{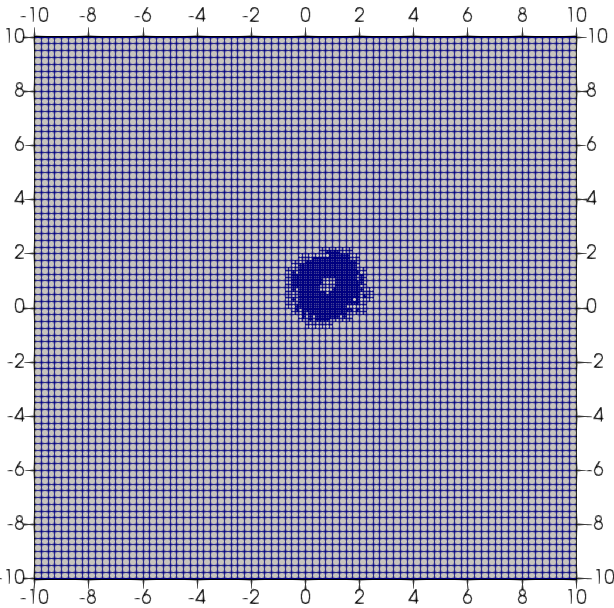} a)
	\includegraphics[width=0.45\textwidth]{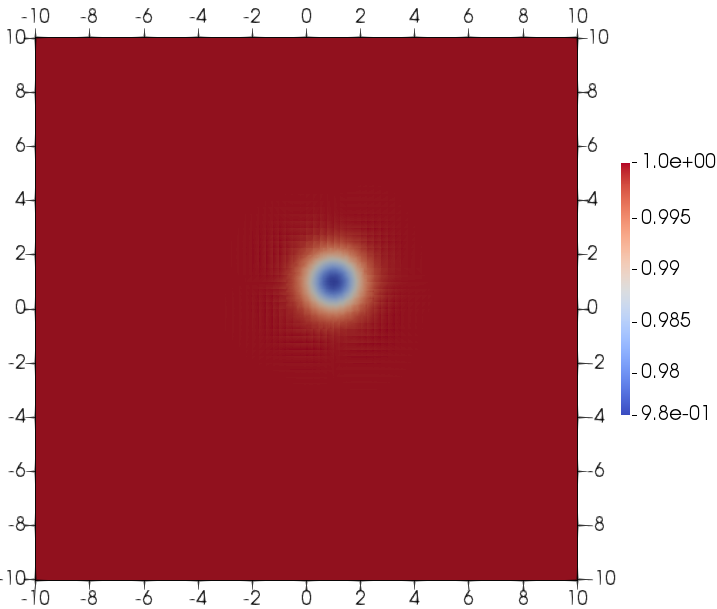} b)
	\caption{Adaptive simulation of the inviscid isentropic vortex benchmark: a) computational mesh at \(t = T_f\), b) contour plot of the density at \(t = T_f\).}
	\label{fig:vortex_adaptive}
\end{figure}

For validation purposes, we have also tested in this case the $h-$adaptive version of the method. The local refinement criterion is based on the gradient of the density. More specifically, we define for each element \(K\) the quantity

\begin{equation}\label{eq:vortex_indicator}
\eta_K = \max_{i \in \mathcal{N}_{K}} \left|\nabla \rho\right|_{i}
\end{equation}
that acts as local refinement indicator, where \(\mathcal{N}_K\) denotes the set of nodes over the element \(K\). Table \ref{tab:vortex_k1_alfa0,5_adaptive} shows the relative errors for all the quantities on a sequence of adaptive simulations keeping the maximum Courant numbers fixed. Figure \ref{fig:vortex_adaptive} shows instead the density and the adapted mesh at \(t = T_f\), from which it can be seen that the refinement criterion is able to track the vortex correctly.

\begin{table}[H]
	\centering
	\begin{tabular}{|c|c|c|c|}
		\hline
		\(N_{el}\) & \(L^2\) rel. error \(\rho\) & \(L^2\) rel. error \(\mathbf{u}\) & \(L^2\) rel. error \(p\) \\
		\hline
		271 & \(2.19 \cdot 10^{-2}\) & \(1.20 \cdot 10^{-2}\) & \(2.97 \cdot 10^{-3}\) \\
		\hline
		586 & \(6.39 \cdot 10^{-4}\) & \(3.09 \cdot 10^{-3}\) & \(9.08 \cdot 10^{-4}\) \\
		\hline
		1999 & \(1.80\cdot10^{-4}\) & \(7.93\cdot10^{-4}\) & \(2.56\cdot10^{-4}\) \\
		\hline
		7678 & \(5.16\cdot10^{-5}\) & \(1.84\cdot10^{-4}\) & \(7.42\cdot10^{-5}\) \\
		\hline 
	\end{tabular}
	\caption{Adaptive simulations of the inviscid isentropic vortex at different resolutions  with a maximum \(C \approx 0.1\), \(C_u\approx 0.1\), relative errors for the density, the velocity and the pressure in $L^2$ norm with \(k = 1\).}
	\label{tab:vortex_k1_alfa0,5_adaptive}
\end{table}

\subsection{Sod shock tube problem}
\label{ssec:Sod}

\begin{figure}[H]
	\includegraphics[width=0.29\textwidth]{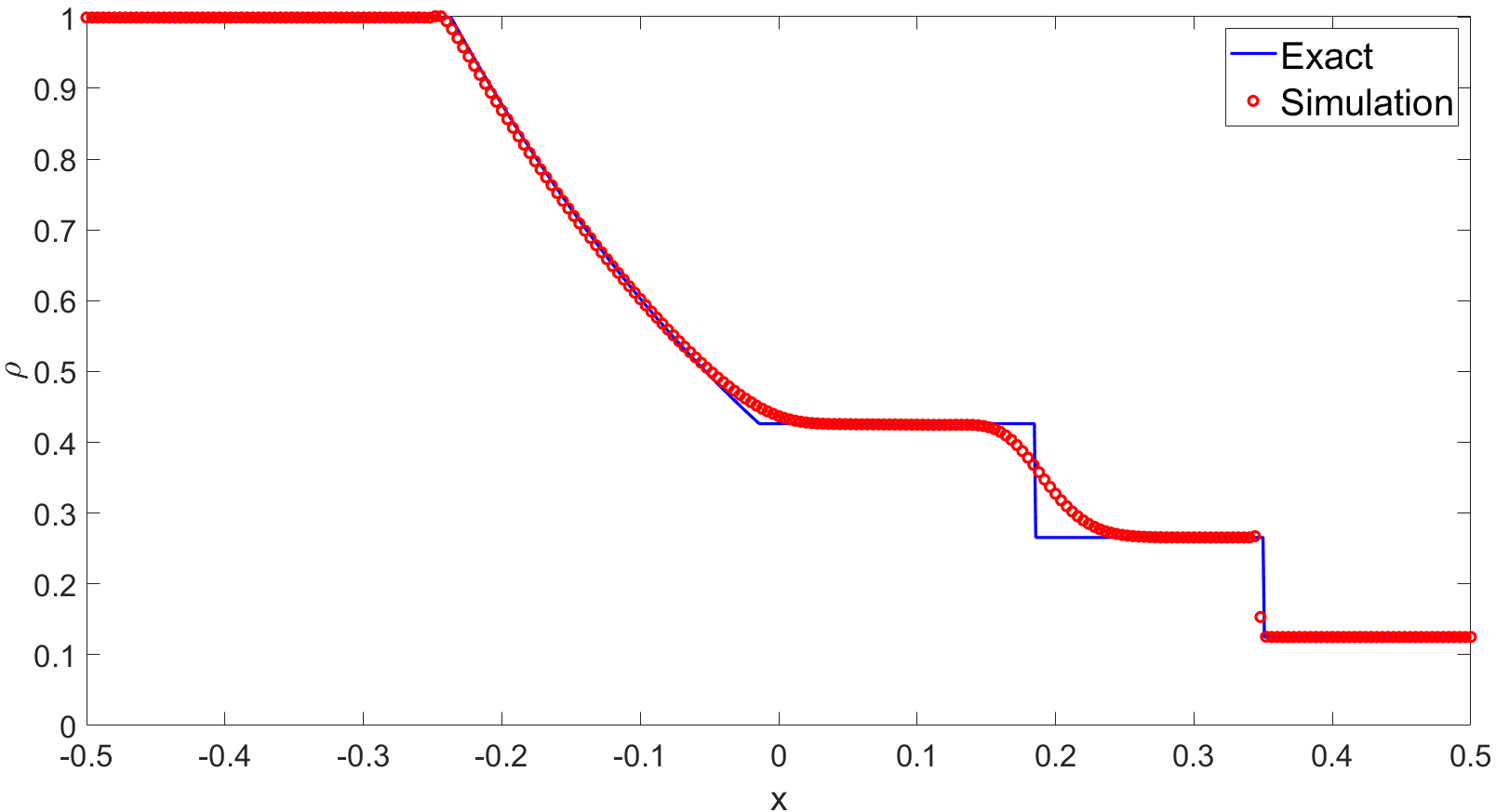} a)
	\includegraphics[width=0.29\textwidth]{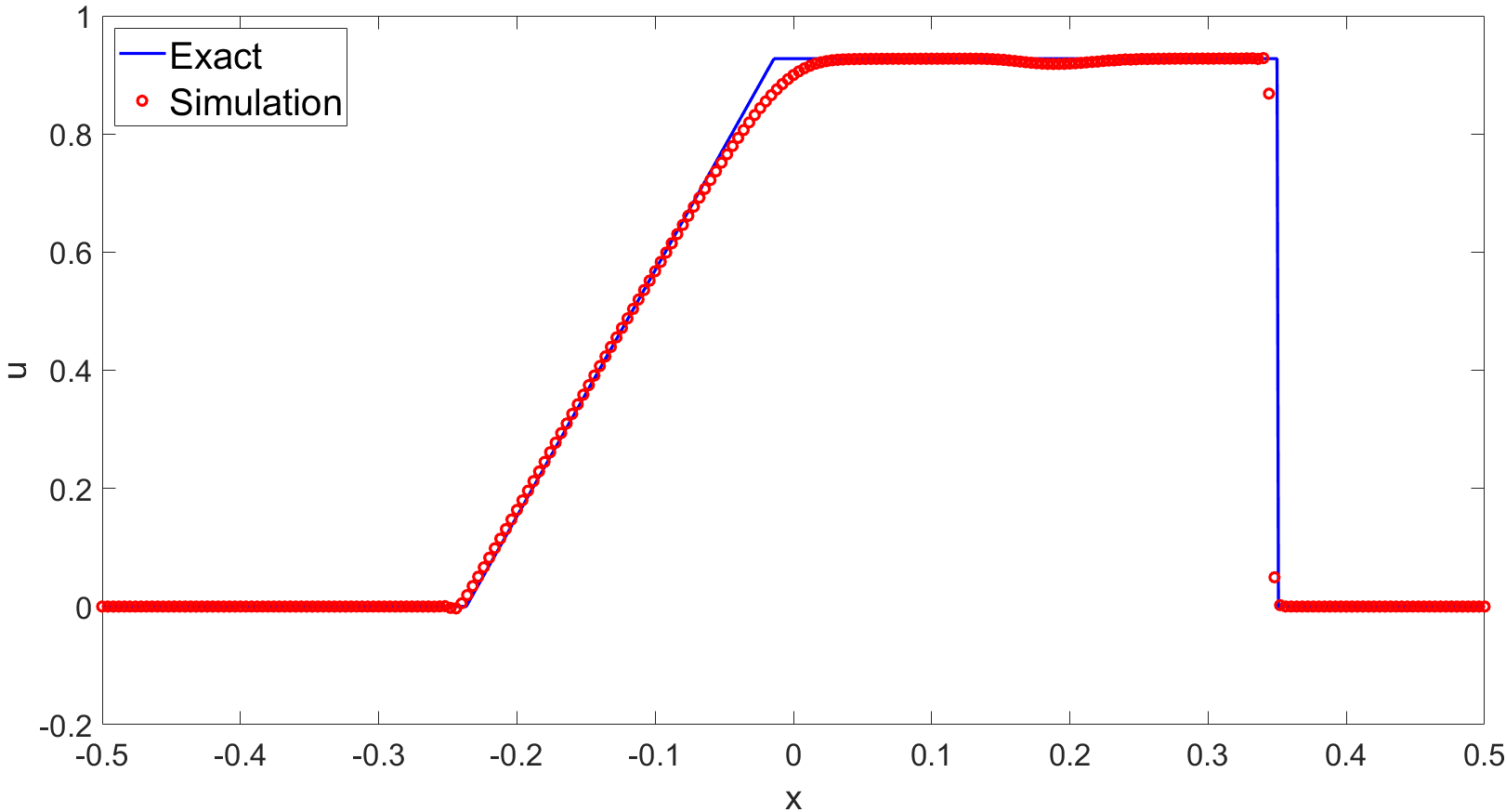} b)
	\includegraphics[width=0.29\textwidth]{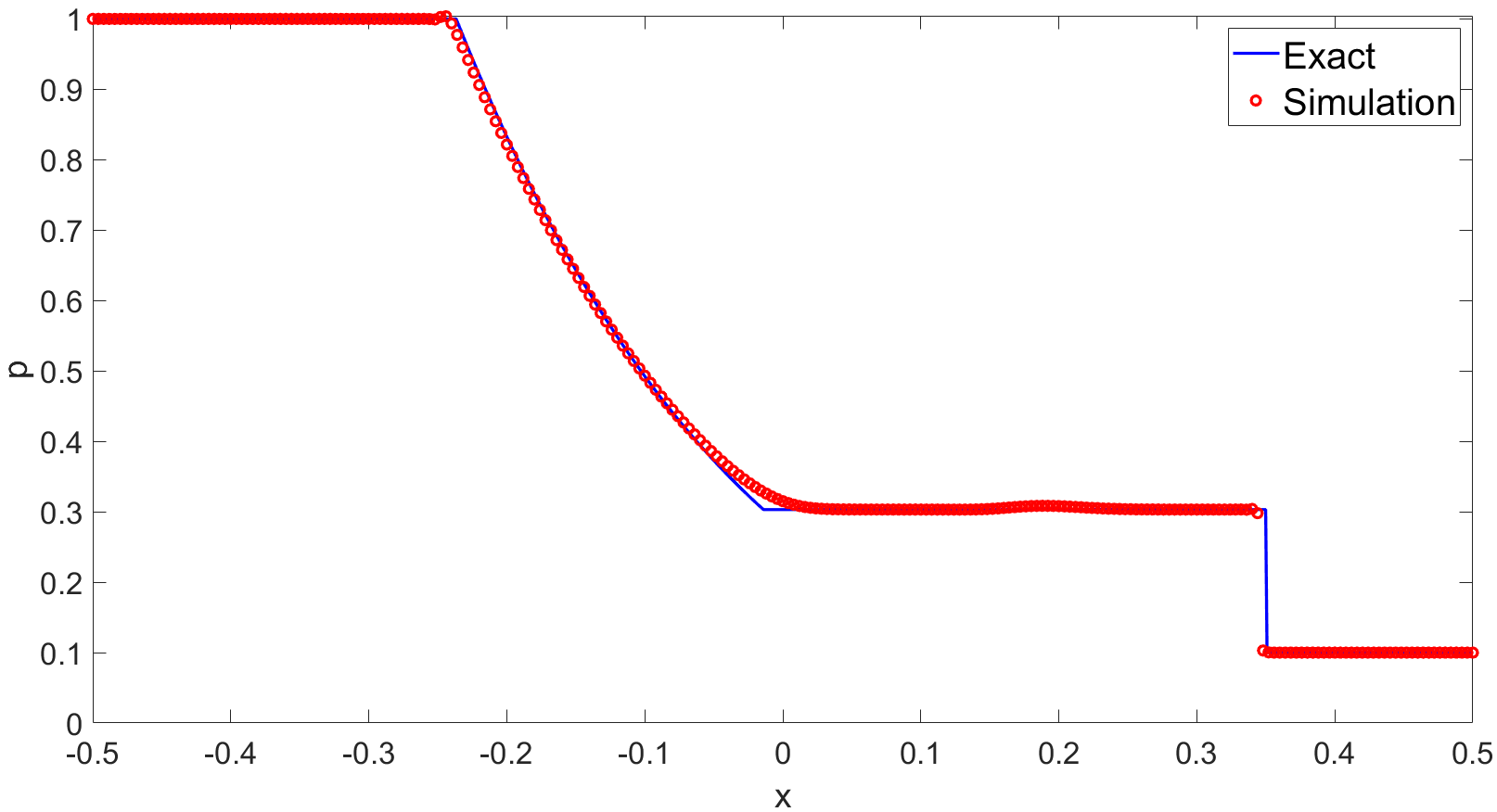} c)
	\caption{Sod shock tube problem at \(t = 0.2\), comparison with exact solution, a) density, b) velocity, c) pressure.}
	\label{fig:Sod_results}
\end{figure}

Even though the proposed method is particularly well suited for low Mach number flows, we have also tested its behaviour also in a situation in which shock waves occur. For this purpose, we have first considered the classical Sod shock tube problem for an ideal gas proposed in \cite{sod:1978} and also discussed in \cite{dumbser:2016}. It consists of a right-moving shock wave, an intermediate contact discontinuity and a left-moving rarefaction fan. In this higher Mach number regime, as done also in \cite{dumbser:2016}, for further stabilization in presence of shocks and discontinuities, the  numerical flux employed for the quantities computed explicitly in the above weak formulations is the classical Local Lax-Friedrichs flux (LLF), instead of the upwind flux, defined by setting

\begin{align*}
\lambda^{(n,1)} &= \max\left\{\left\|\mathbf{u}^{(n,1)^{+}}\right\| + \frac{c^{(n,1)^{+}}}{M},\left\|\mathbf{u}^{(n,1)^{-}}\right\| + \frac{c^{(n,1)^{-}}}{M}\right\} \\
\lambda^{(n,2)} &= \max\left\{\left\|\mathbf{u}^{(n,2)^{+}}\right\| + \frac{c^{(n,2)^{+}}}{M},\left\|\mathbf{u}^{(n,2)^{-}}\right\| + \frac{c^{(n,2)^{-}}}{M}\right\}.
\end{align*}
The presence of discontinuities requires the use of a monotonic scheme to avoid undershoots and overshoots. It is well known that using \(Q_0\) finite elements in combination with LLF and an explicit time integration method that complies with the monotonicity constraints discussed in \cite{ferracina:2004, higueras:2004, gottlieb:2001} guarantees the monotonicity of the solution. Hence, a way to obtain monotonic results is to project the numerical solution onto the \(Q_0\) subspace for each element in which a suitable jump indicator exceeds a certain threshold. Similar projections onto low order components of the solution are also used in several monotonization approaches, see e.g. \cite{dumbser:2014}.
However, since in the proposed scheme only the density is treated in a full explicit fashion, in order to avoid an excessive complication in the structure of the resulting method we choose to apply this \(Q_0\) projection strategy only for the density variable, without introducing monotonization for the velocity and the pressure. While we are aware that this is not sufficient to guarantee full monotonicity, the derivation of a fully monotonic IMEX scheme goes beyond the scope of this work and we do not investigate this issue further here. Therefore, the results in this Section are to be interpreted merely as a first stress test of the proposed scheme at higher Mach number values. In future work, we plan to investigate the behaviour of the monotonization approach proposed in \cite{orlando:2022} for full explicit schemes, when applied only to the density.  
We use a smoothness indicator based on the jump of the density across two faces. More in detail, we define for each element \(K\) the quantity

\[\eta_K = \sum_{\Gamma \in \mathcal{E}_K}\left\|\rho^{+} - \rho^{-}\right\|_{2, \Gamma}^2\] 
where \(\mathcal{E}_K\) denotes the set of all faces belonging to cell \(K\) and \(\left\|\cdot\right\|\) represents the standard \(L^2\) norm on \(\Gamma\). The chosen threshold in this case is equal to \(10^{-6}\). Table \ref{tab:Sod_setting} defines the initial conditions and position of the initial discontinuity. We consider the domain \(\Omega = \left(-0.5, 0.5\right)\) and a one-dimensional mesh  composed by \(500\) elements with a time step \(\Delta t = 1 \cdot 10^{-4}\), chosen in such a way that the maximum Courant number is  \(C \approx 0.07\), while the maximum advective Courant number is \(C_u \approx 0.06\). Figure \ref{fig:Sod_results} shows the results for the density, the velocity and the pressure at \(t = 0.2\) compared with the exact solution. One can easily notice that the shocks are located at the right position and that the  values in the  wake of the shocks are correct. The technique applied to guarantee the monotonicity introduces however an excessive amount of numerical diffusion in correspondence of contact wave, which is not sharply resolved, in contrast  to what happens for the shock and the rarefaction waves. On the other hand, if one decreases the value of the threshold, an oscillating solution is obtained. While far from optimal, these results highlight however the robustness of the proposed approach also in the higher Mach number case.

\begin{table}[H]
	\centering
	\begin{tabular}{|c|c|c|c|c|c|c|}
		\hline
		\(\rho_L\) & \(u_L\) & \(p_L\) & \(\rho_R\) & \(u_R\) & \(p_R\) & \(x_d\)  \\
		\hline
		1 & 0 & 1 & 0.125 & 0 & 0.1 & 0 \\
		\hline
	\end{tabular}
	\caption{Initial left and right states for Sod shock tube problem. \(x_d\) denotes the position of the initial discontinuity.}
	\label{tab:Sod_setting}
\end{table}

\begin{figure}[H]
	\includegraphics[width=0.29\textwidth]{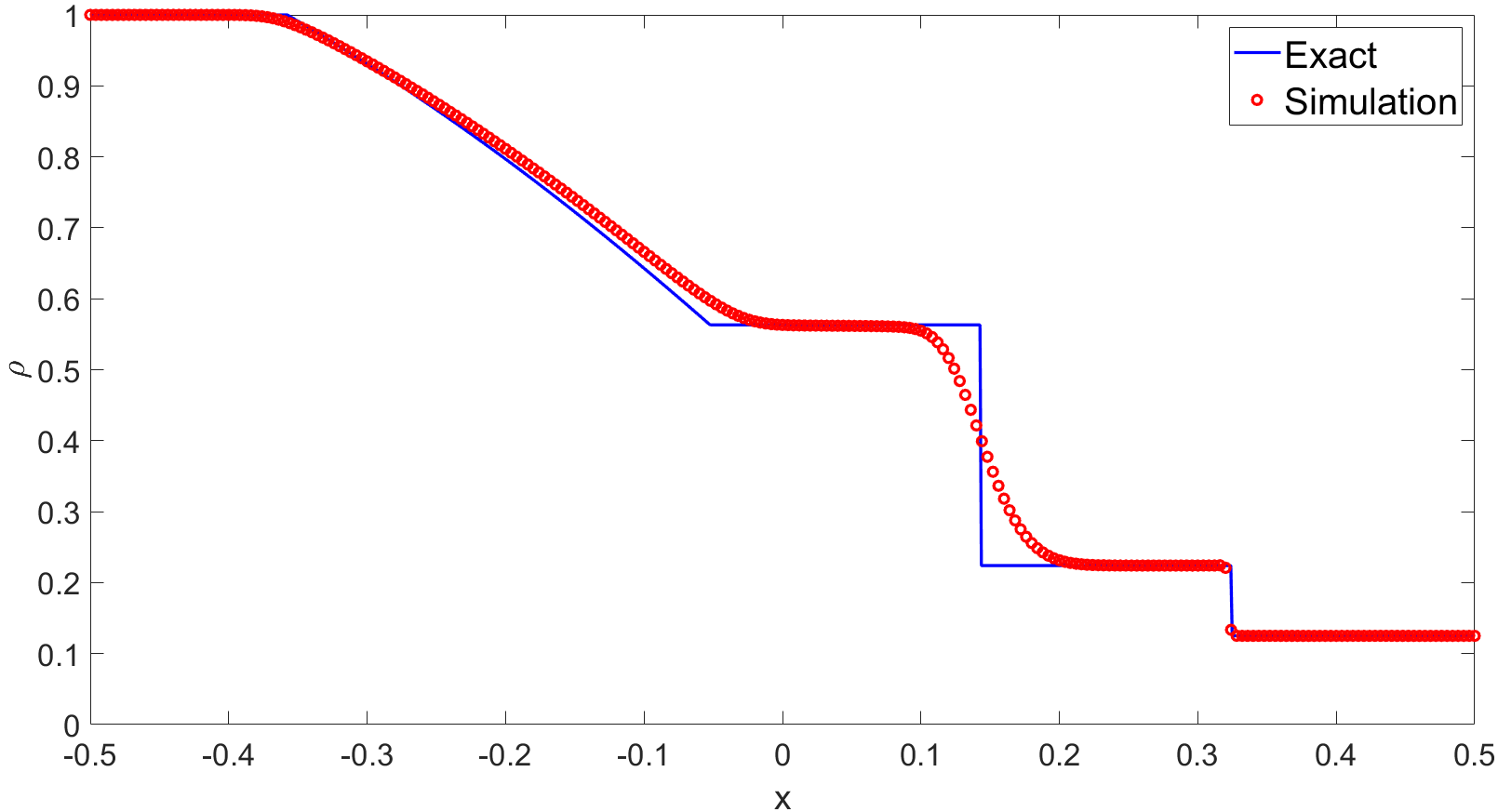} a)
	\includegraphics[width=0.29\textwidth]{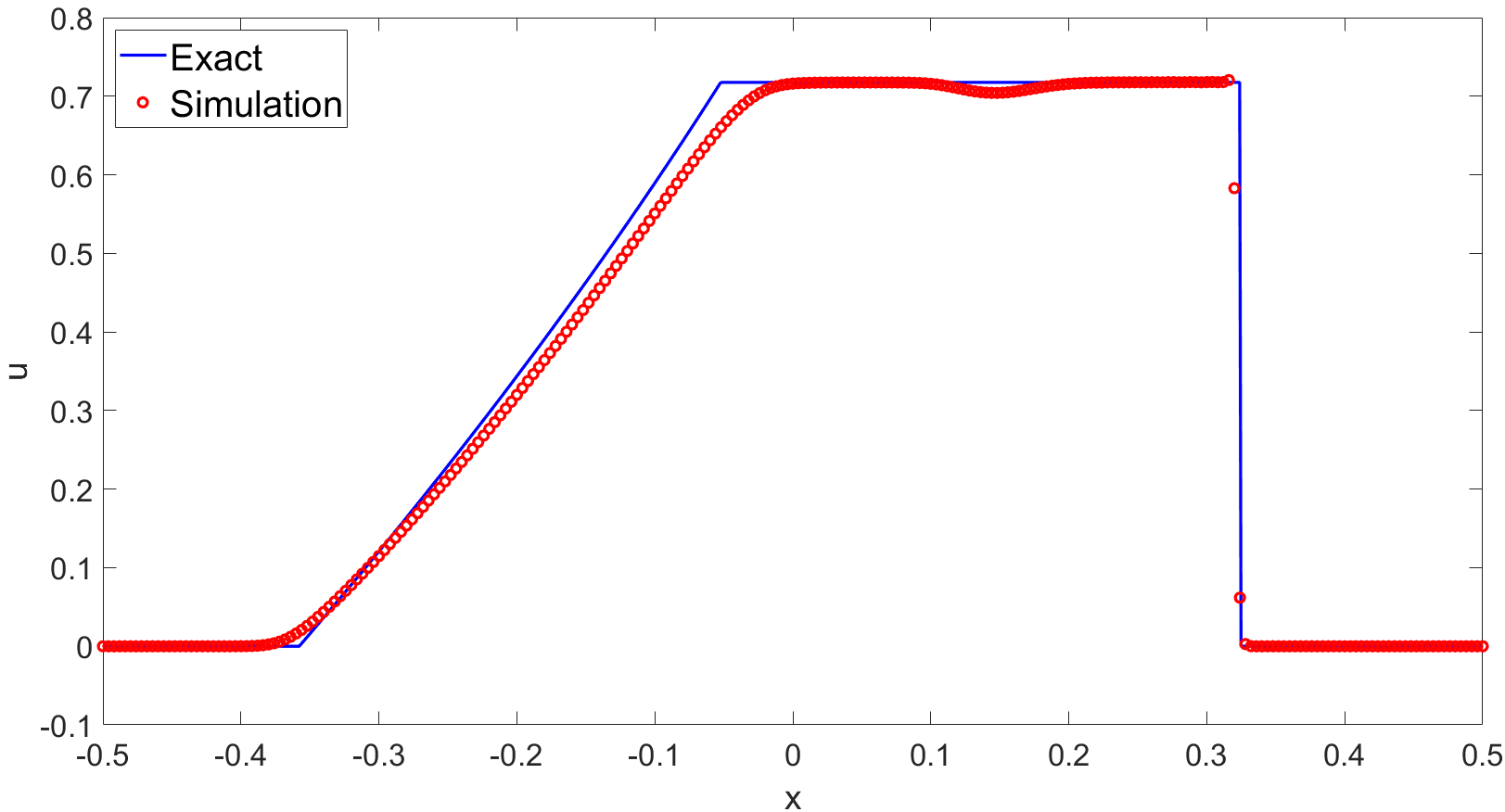} b)
	\includegraphics[width=0.29\textwidth]{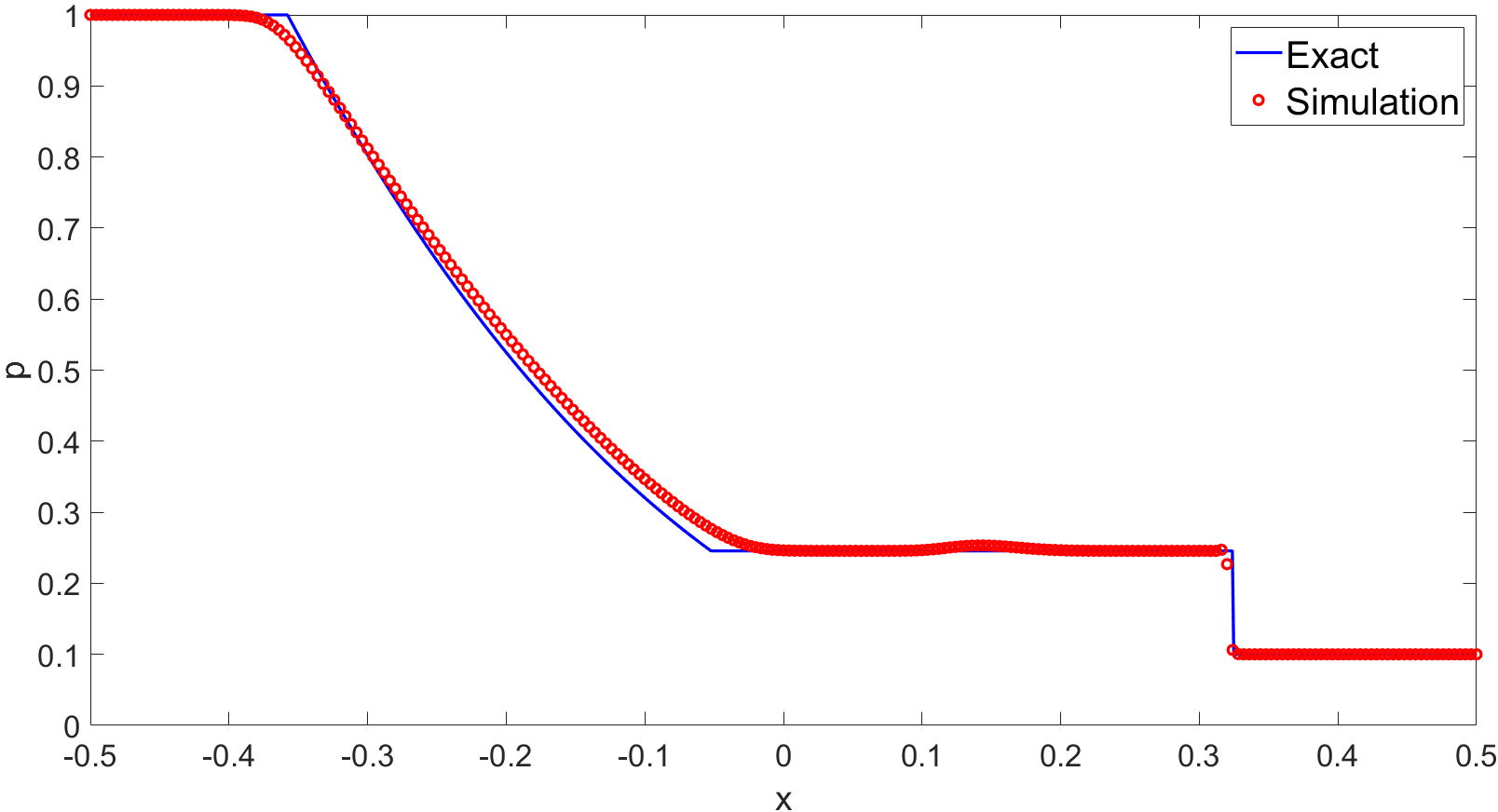} c)
	\caption{Sod shock tube problem with van der Waals EOS at \(t = 0.2\), comparison with exact solution, a) density, b) velocity, c) pressure.}
	\label{fig:Sod_results_vdW}
\end{figure}

Following \cite{dumbser:2016}, we have also considered the same problem in the case of the van der Waals EOS, taking \(\tilde a = \tilde b = 0.5\). All the other parameters are the same as in the previous case. Figure \ref{fig:Sod_results_vdW} shows the results for the density, the velocity and the pressure at \(t = 0.2\). 
With the same limitations discussed before, a good agreement between the numerical results and the exact solution is observed also in this case for all the quantities.

\begin{figure}[H]
	\includegraphics[width=0.29\textwidth]{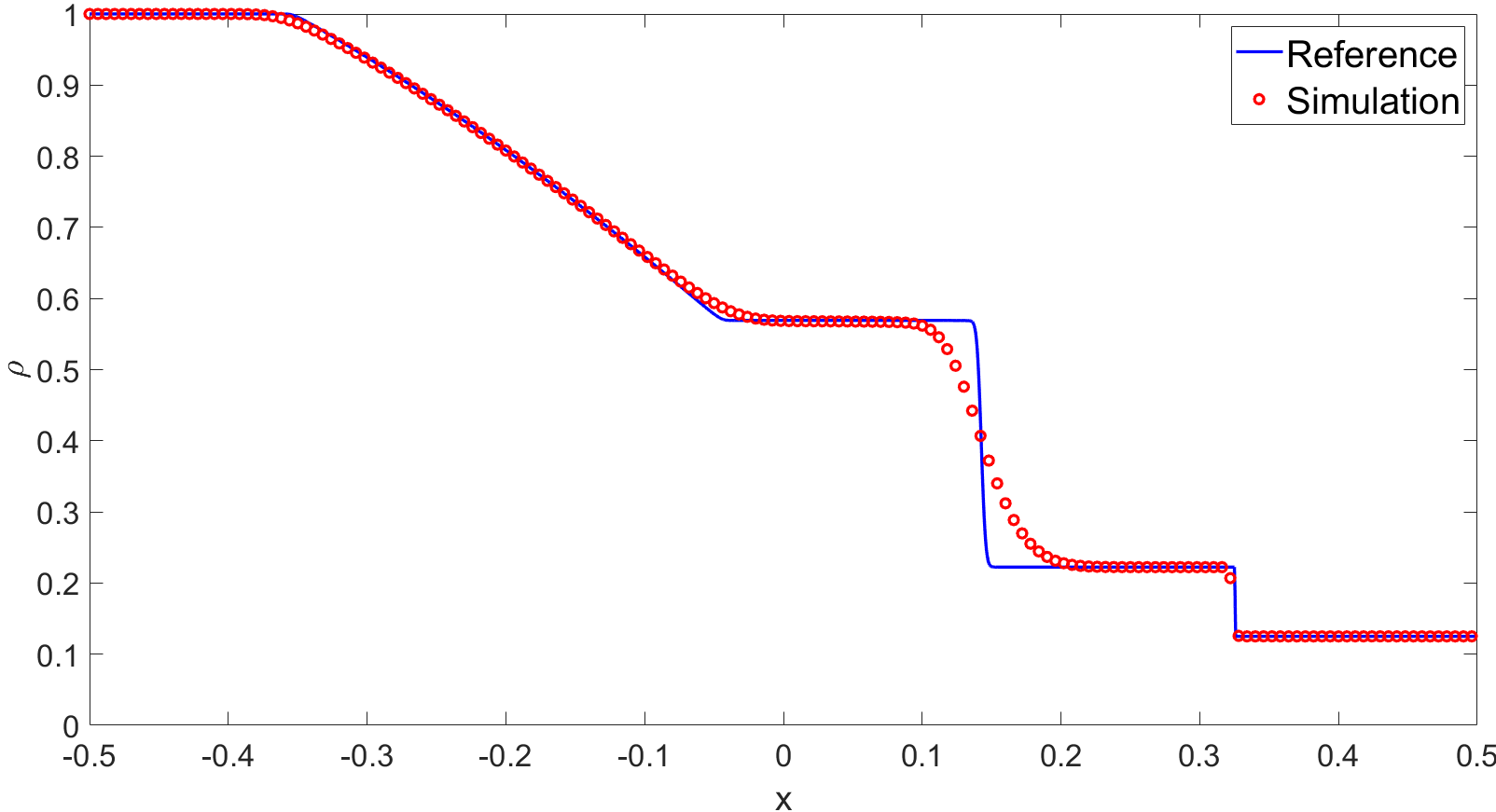} a)
	\includegraphics[width=0.29\textwidth]{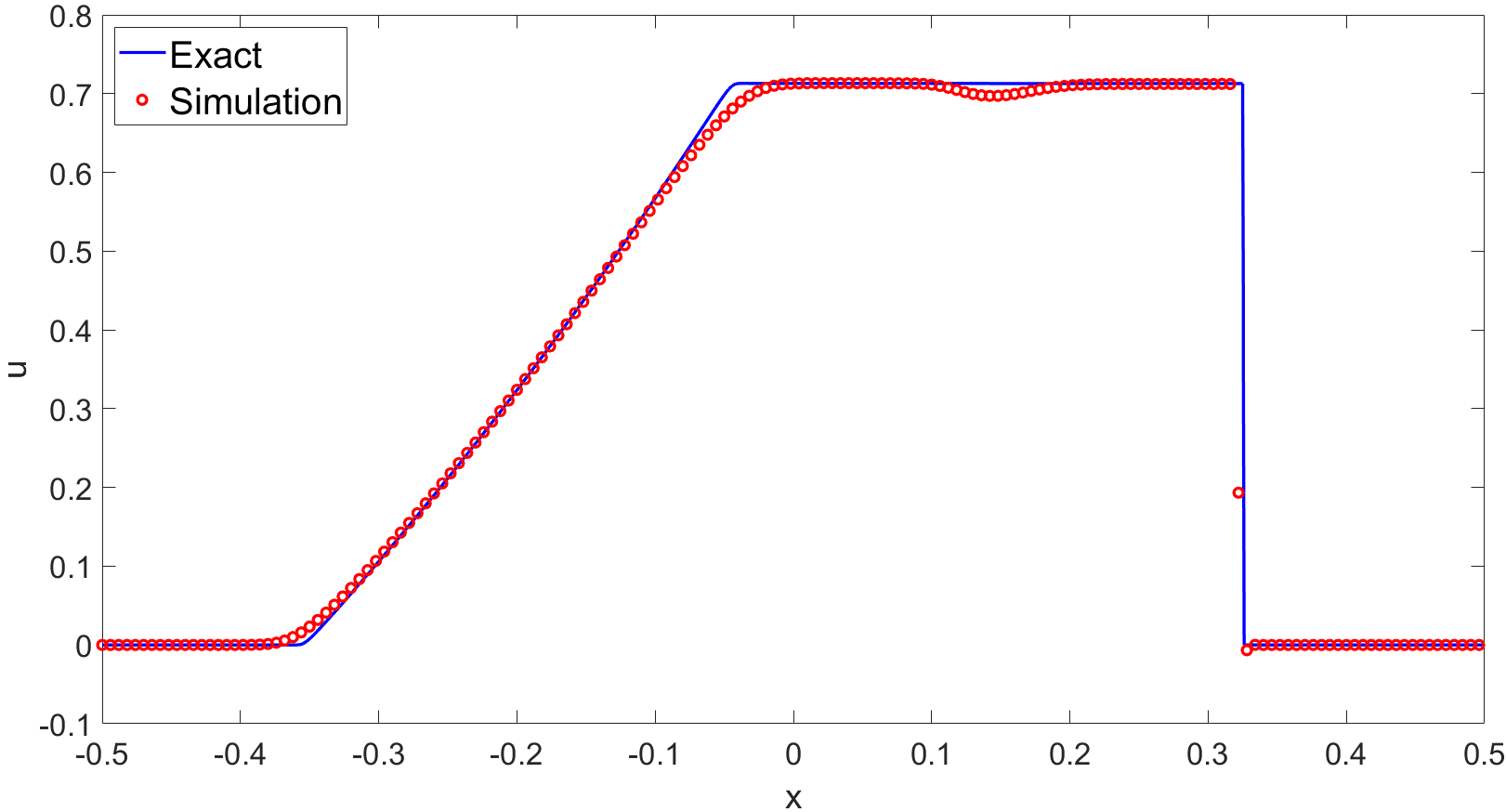} b)
	\includegraphics[width=0.29\textwidth]{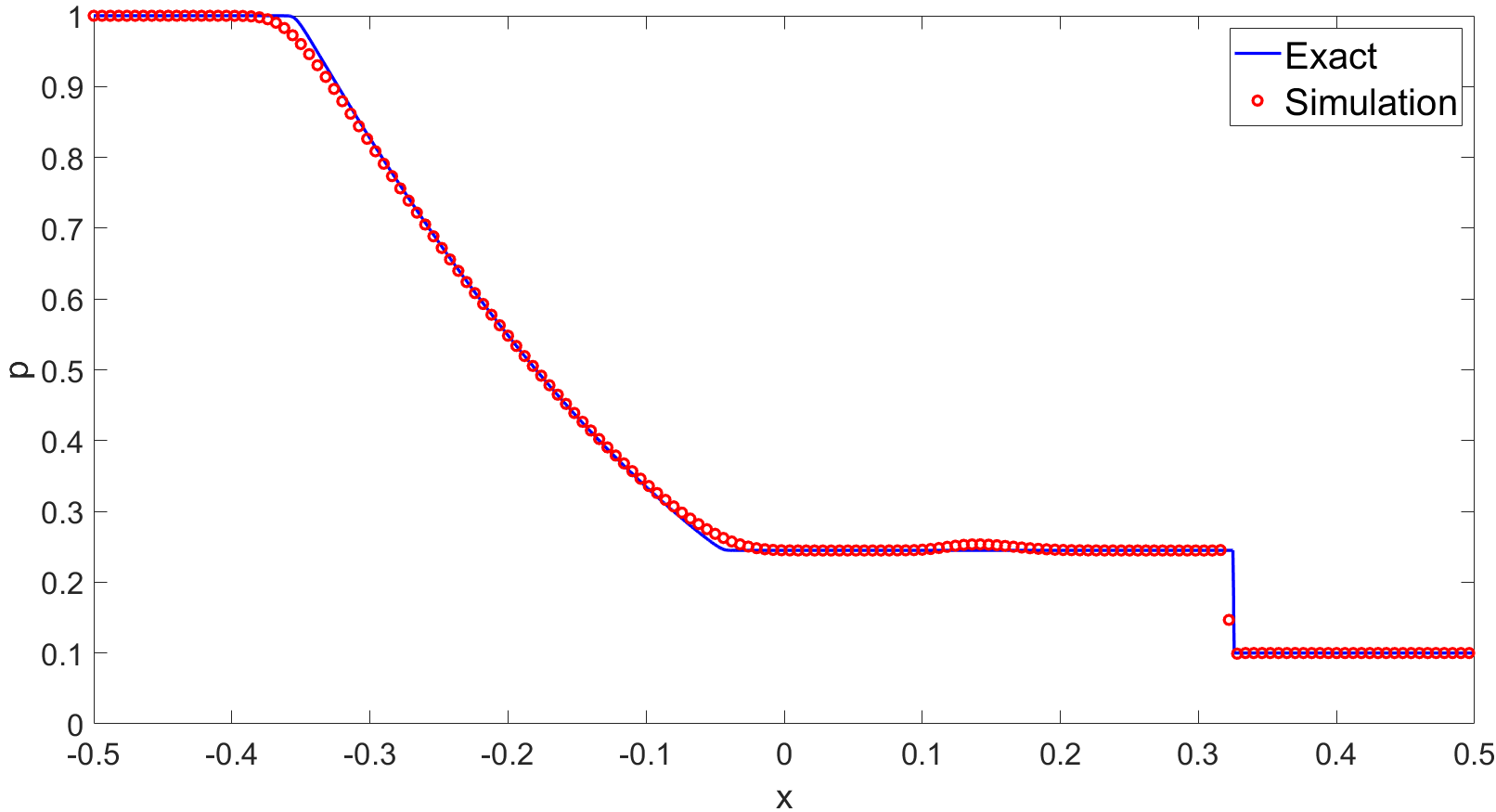} c)
	\caption{Sod shock tube problem with Peng-Robinson EOS at \(t = 0.2\), comparison with reference solution, a) density, b) velocity, c) pressure.}
	\label{fig:Sod_results_PG}
\end{figure}

Finally, we have analyzed the same test case for the Peng-Robinson EOS with \(\tilde a = \tilde b = 0.5\). In this case no analytic solution is available and a reference solution is computed using the third order optimal SSP Runge-Kutta method derived in \cite{gottlieb:1998} in combination with \(Q_{0}\) finite elements and LLF as numerical flux on a mesh with \(16000\) elements.  
Notice that, in the case of the van der Waals EOS,  this procedure was found to provide a solution  overlapping with the exact one. The good agreement between our numerical results and the reference ones is established also for this equation of state, as evident from Figure \ref{fig:Sod_results_PG}, again with significant smoothing
of the contact discontinuity. 

\subsection{2D Lid-driven cavity}
\label{ssec:2dlid}

\begin{figure}[H]
	\includegraphics[width=0.45\textwidth]{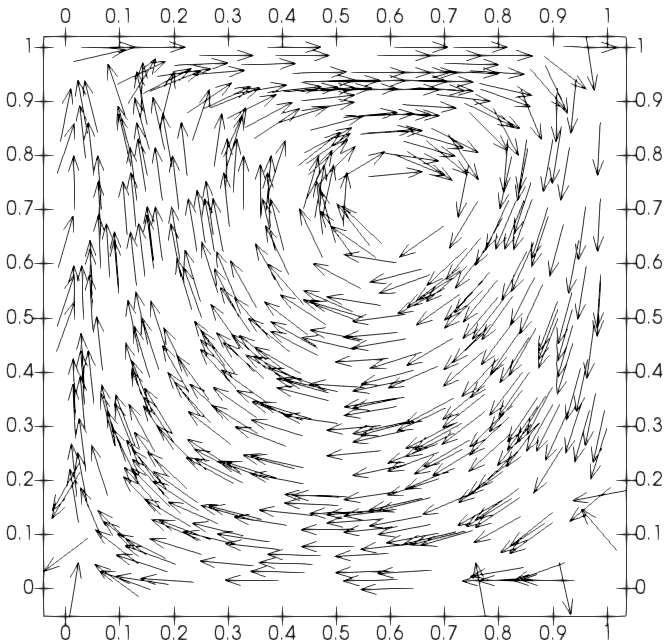} a)
	\includegraphics[width=0.45\textwidth]{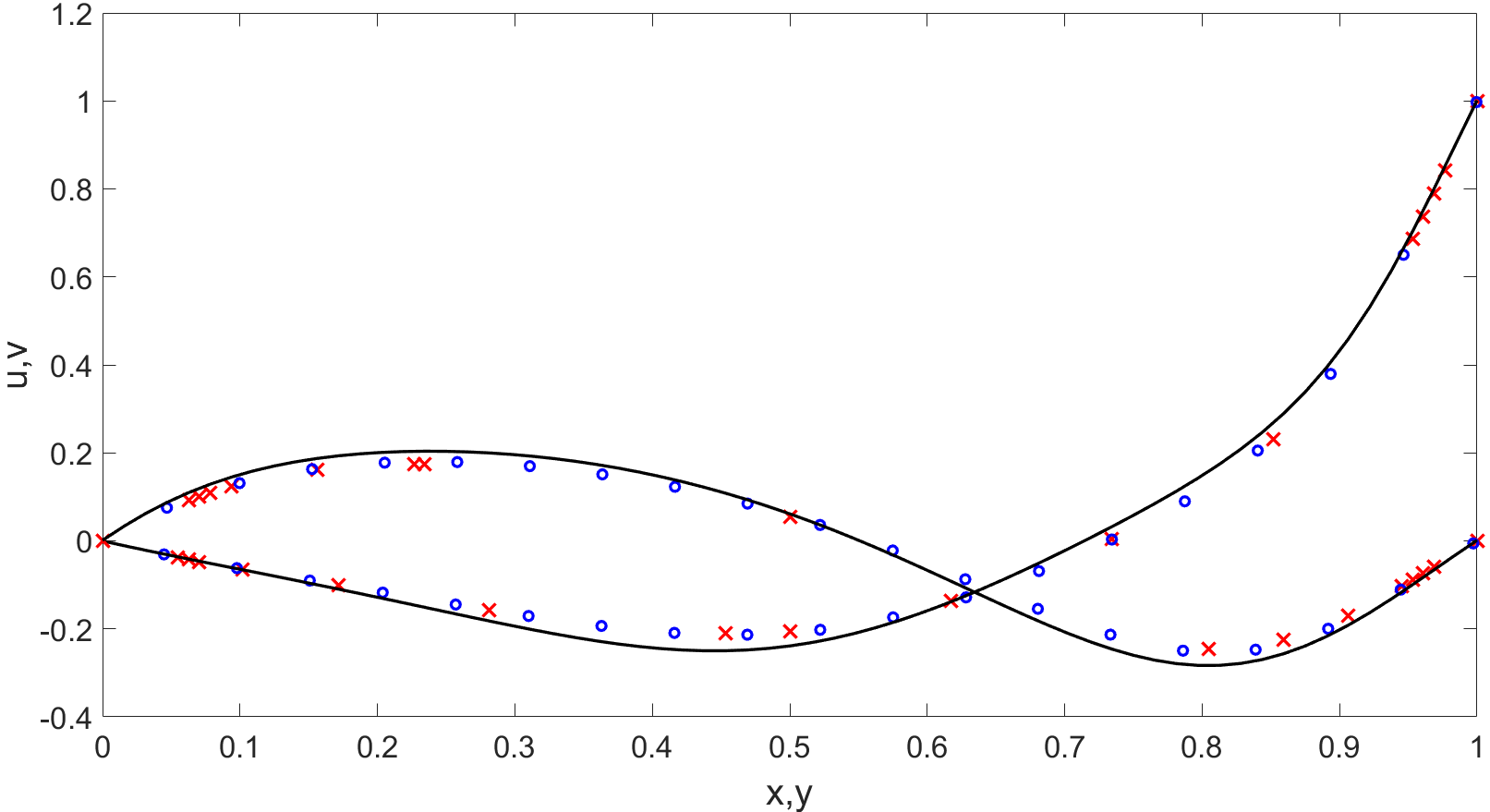} b)
	\caption{Computational results for the 2D lid-driven cavity, a) streamlines, b) comparison with the solutions in \cite{ghia:1982} and in \cite{tavelli:2017}. Blue dots denote the results in \cite{ghia:1982}, red crosses the results in \cite{tavelli:2017} and the  black line our numerical results.}
	\label{fig:Lid_Driven_Cavity_Q1}
\end{figure}

We consider now the classical 2D lid-driven cavity test case. The computational domain is the box \(\Omega = \left(0,1\right) \times \left(0,1\right)\) which is initialized with a density \(\rho = 1\) and a velocity \(\mathbf{u} = \mathbf{0}\). The flow is driven by the upper boundary, whose velocity is set to \(\mathbf{u} = \left(1,0\right)^{T}\), while on the other three boundaries a no-slip condition is imposed. We set \(Re = 100\) and \(M^2 = 10^{-5}\). The advantage of the proposed scheme is that the allowed time step is more than 100 times larger than that of a fully explicit scheme. Indeed, the time-step chosen is such that the maximum advective Courant number \(C_u\) is around \(0.12\), while the maximum Courant number \(C\) is around \(49\). The streamlines are shown in Figure \ref{fig:Lid_Driven_Cavity_Q1} and highlight the formation of the main recirculation pattern. A comparison of the horizontal component of the velocity along the vertical middle line and of the vertical component of the velocity along the horizontal middle line with the reference solutions in \cite{ghia:1982, tavelli:2017} is also presented. 

We note a reasonable agreement between the different solutions, even though there is a still visible discrepancy between our results and the reference ones. Since the solution in \cite{tavelli:2017} is obtained using third degree polynomials, in order to further improve the results, we consider also the case \(r = 2\). For this higher order approximation we note that our results fit very well both the reference solutions. 

\begin{figure}[H]
	\includegraphics[width=0.45\textwidth]{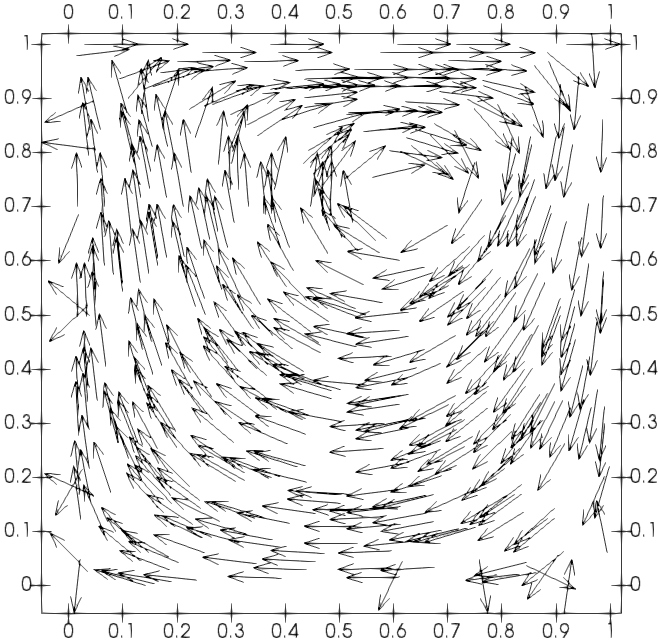} a)
	\includegraphics[width=0.45\textwidth]{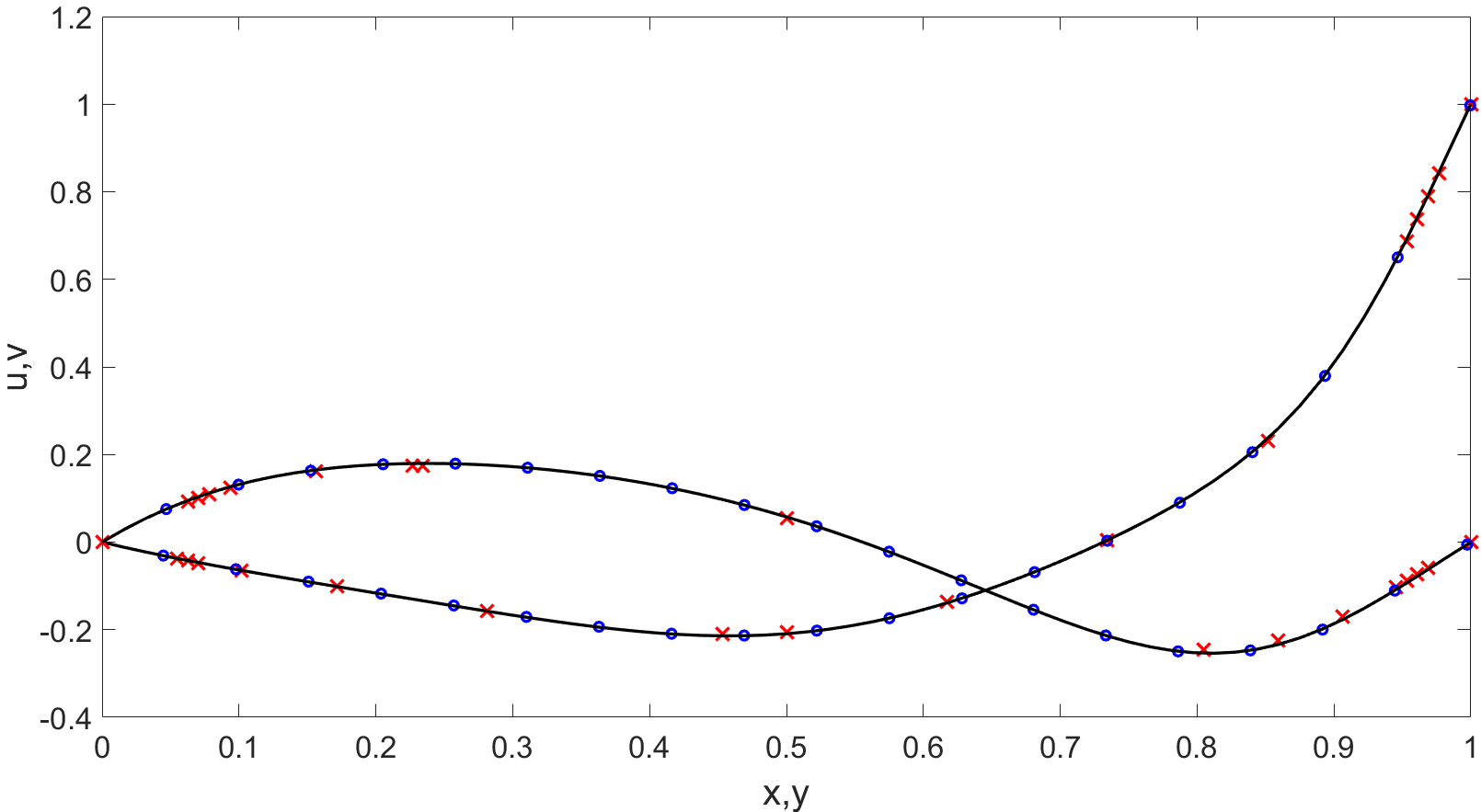} b)
	\caption{Computational results for the 2D lid-driven cavity with \(k = 2\), a) streamlines, b) comparison with the solutions in \cite{ghia:1982} and in \cite{tavelli:2017}. Blue dots denote the results in \cite{ghia:1982}, red crosses the results in \cite{tavelli:2017} and black line our numerical results.}
	\label{fig:Lid_Driven_Cavity_Q2}
\end{figure}

\begin{figure}[H]
	\includegraphics[width=0.45\textwidth]{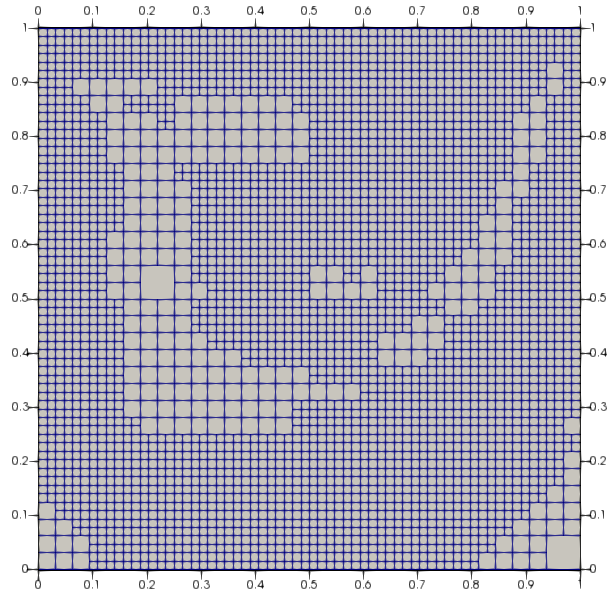} a)
	\includegraphics[width=0.45\textwidth]{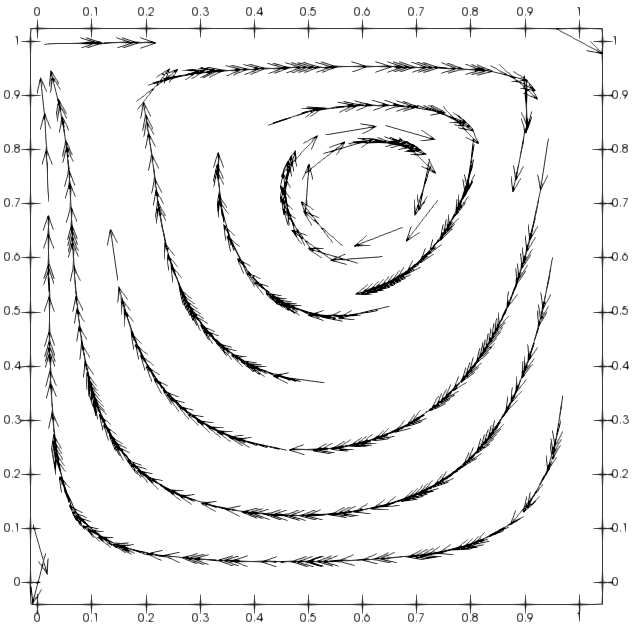} b)
	\caption{Adaptive simulation for the 2D lid-driven cavity, a) mesh at steady state, b) streamlines.}
	\label{fig:Lid_Driven_Cavity_adaptive}
\end{figure}

\begin{figure}[H]
	\includegraphics[width=0.45\textwidth]{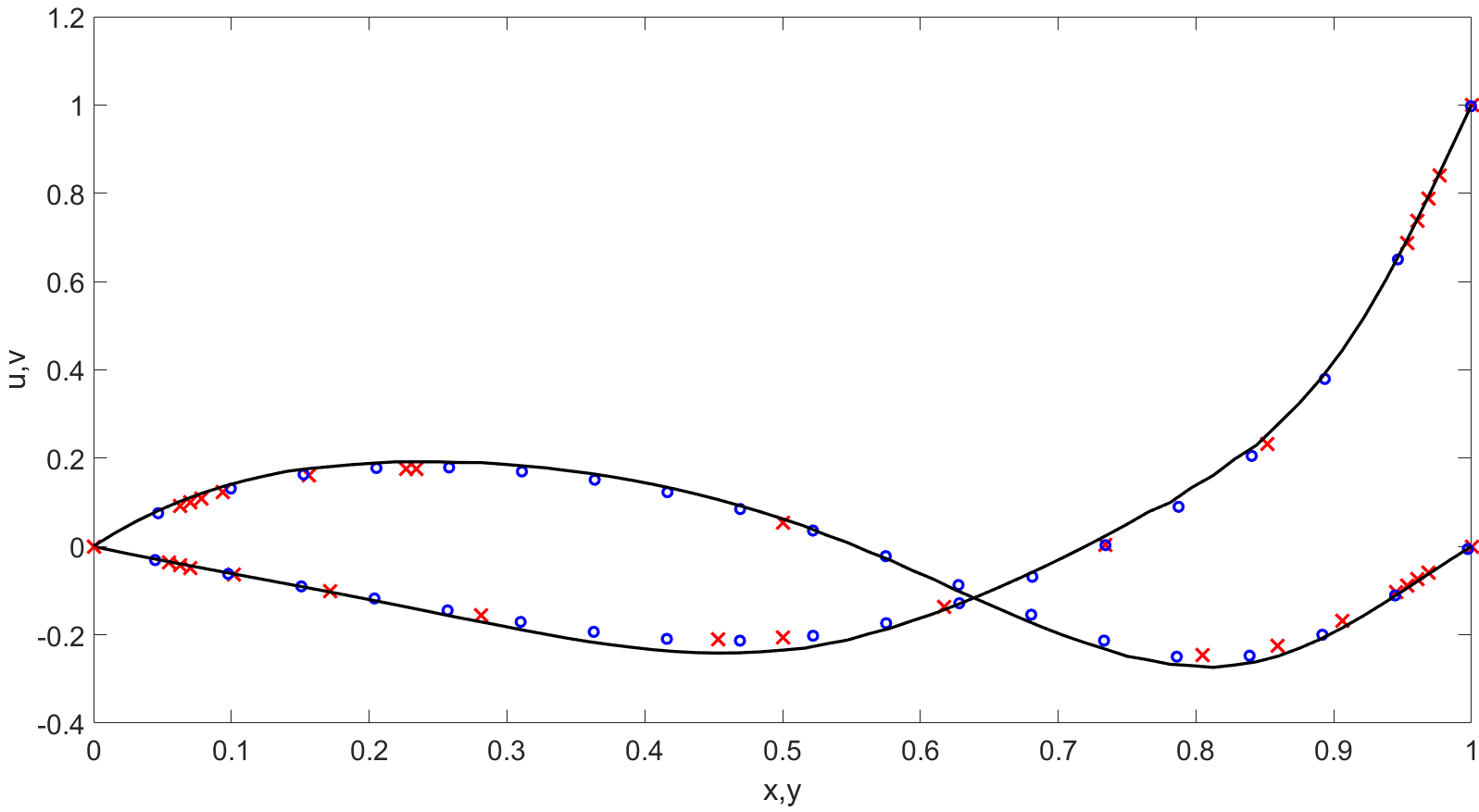} a)
	\includegraphics[width=0.45\textwidth]{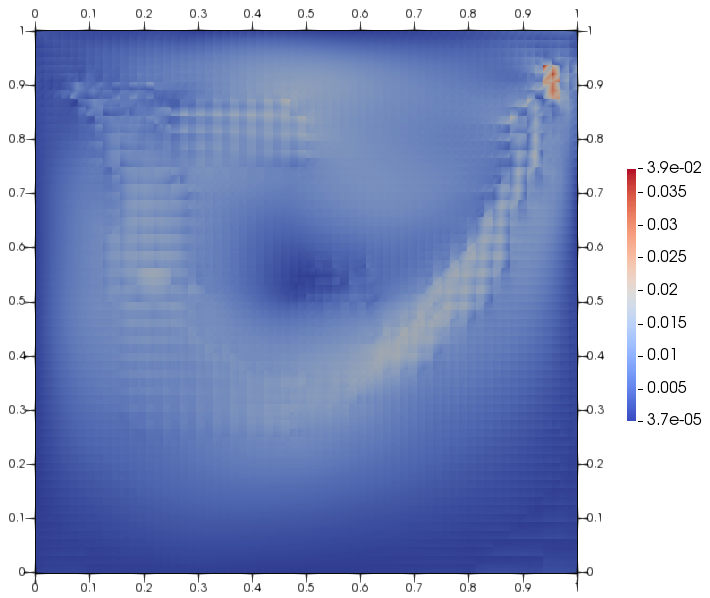} b)
	\caption{Adaptive simulation for the 2D lid-driven cavity, a) comparison with the solutions in \cite{ghia:1982} and in \cite{tavelli:2017}. Blue dots denote the results in \cite{ghia:1982}, red crosses the results in \cite{tavelli:2017} and black line our numerical results, b) difference for velocity magnitude between the fixed grid simulation and the adaptive simulation (interpolated to the fixed grid).}
	\label{fig:Lid_Driven_Cavity_adaptive_bis}
\end{figure}

We have also tested the $h-$adaptive version of the same algorithm, using a refinement criterion based on the vorticity. More specifically, we define 

\[\eta_K = \text{diam}\left(K\right)^{2} \left\|\nabla\times\mathbf{u}\right\|^{2}_{2,K}\]
as local indicator. We start from a uniform Cartesian mesh with \(16\) elements along each direction. We allowed refinement for \(5\%\) of the elements with largest indicator values  and coarsening for \(30\%\) of the elements with the smallest indicator values. The minimum element diameter allowed is \(\mathcal{H} = \frac{1}{64}\), while the maximum element diameter is \(\mathcal{H} = \frac{1}{16}\). Figure \ref{fig:Lid_Driven_Cavity_adaptive} reports the computational mesh at steady state and the computed streamlines. One can easily notice that the local refinement criterion is able to  enhance automatically the resolution in the zones where vortices appear, as well as along the top boundary of the domain. For a more quantitative view, in Figure \ref{fig:Lid_Driven_Cavity_adaptive_bis} we compare again the components of the velocity along the middle lines. Moreover, the absolute difference between the velocities of the fixed mesh and those of adaptive simulations is plotted over the whole domain, showing that no substantial loss of accuracy has occurred with a reduction of around 25\% of the required computational time. 

\subsection{Cold bubble}
\label{ssec:cold_bubble}

In this Section, we consider a test case proposed in \cite{restelli:2007, restelli:2009} for an ideal gas in which the gravity force is active. The computational domain is the rectangle \(\left(0,1000\right) \times \left(0,2000\right)\) and the initial condition is represented by a thermal anomaly introduced in an isentropic background atmosphere with constant potential temperature \(\theta_0 = 303\). The perturbation of potential temperature \(\theta^{'}\) defines the initial datum and it is given by

\begin{equation}
\theta^{'} = 
\begin{cases}
A & \text{if } \tilde r \le r_0 \\
A\exp\left(-\frac{\left(\tilde r - r_0\right)^2}{\sigma^2}\right) & \text{if } \tilde r > r_0,
\end{cases}
\end{equation}
with \(\tilde r^2 = \left(x - x_0\right)^2 + \left(y - y_0\right)^2\) and \(x_0 = 500\), \(y_0 = 1250\), \(r_0 = 50\), \(\sigma = 100\) and \(A = -15\). Moreover, we set \(Fr^2 = \frac{1}{9.81}, M^2 = 10^{-5}\) and \(T_{f} = 50\). The expression of the Exner pressure is given by 

\[\Pi = 1 - \frac{M^2}{Fr^2}\frac{y}{\tilde c_{p}\theta},\]
with \(y\) denoting the vertical coordinate and \(\tilde c_{p} = \frac{\gamma}{\gamma - 1}\tilde R_{g} = 1.0045 \cdot 10^{-2}\) denoting the non-dimensional specific heat at constant pressure. Notice that these values are obtained by considering \(\cal R = \SI{1}{\kilogram\per\meter\cubed}, \cal T \cal T = \SI{1}{\kelvin}\) and \({\cal P} = 10^{5} \hspace{0.05cm} \SI{}{\pascal}\). Moreover, it is to be remarked that, unlike in \cite{restelli:2007}, no artificial viscosity has been added to stabilize the computation. Wall boundary conditions are imposed at all the boundaries. The time step is taken to be \(\Delta t = 0.08\), corresponding to a maximum Courant number \(C \approx 5.6\) and a maximum advective Courant number \(C_u \approx 0.18\). 

\begin{figure}[H]
	\centering
	\begin{subfigure}{0.5\textwidth}
		\includegraphics[width=0.9\textwidth]{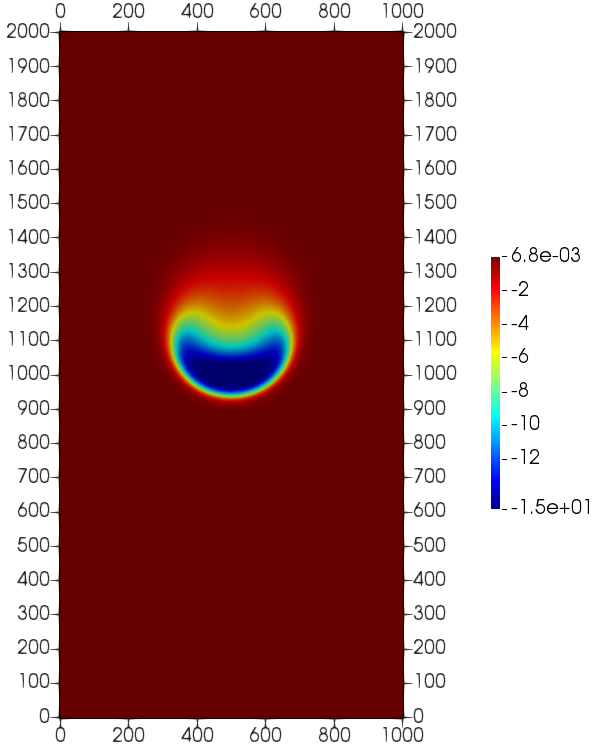} a)
	\end{subfigure}\hspace*{\fill}
	\begin{subfigure}{0.5\textwidth}
		\includegraphics[width=0.9\textwidth]{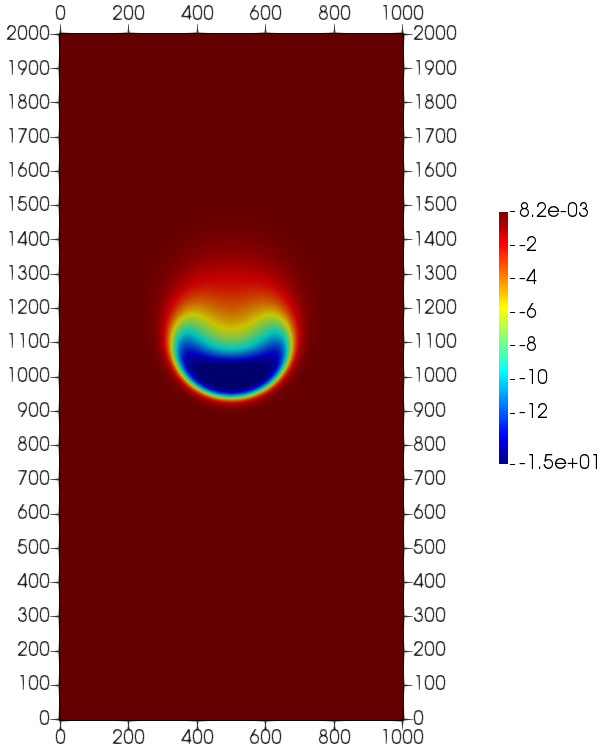} b)
	\end{subfigure}
	\caption{Cold bubble test case, results at \(t = T_{f}\), a) contour plot of potential temperature perturbation for the reference explicit simulation, b) contour plot of the potential temperature perturbation for the simulation with IMEX scheme.} 
	\label{fig:Cold_bubble_ideal}
\end{figure}

\begin{figure}[H]
	\centering
	\includegraphics[width = 0.8\textwidth]{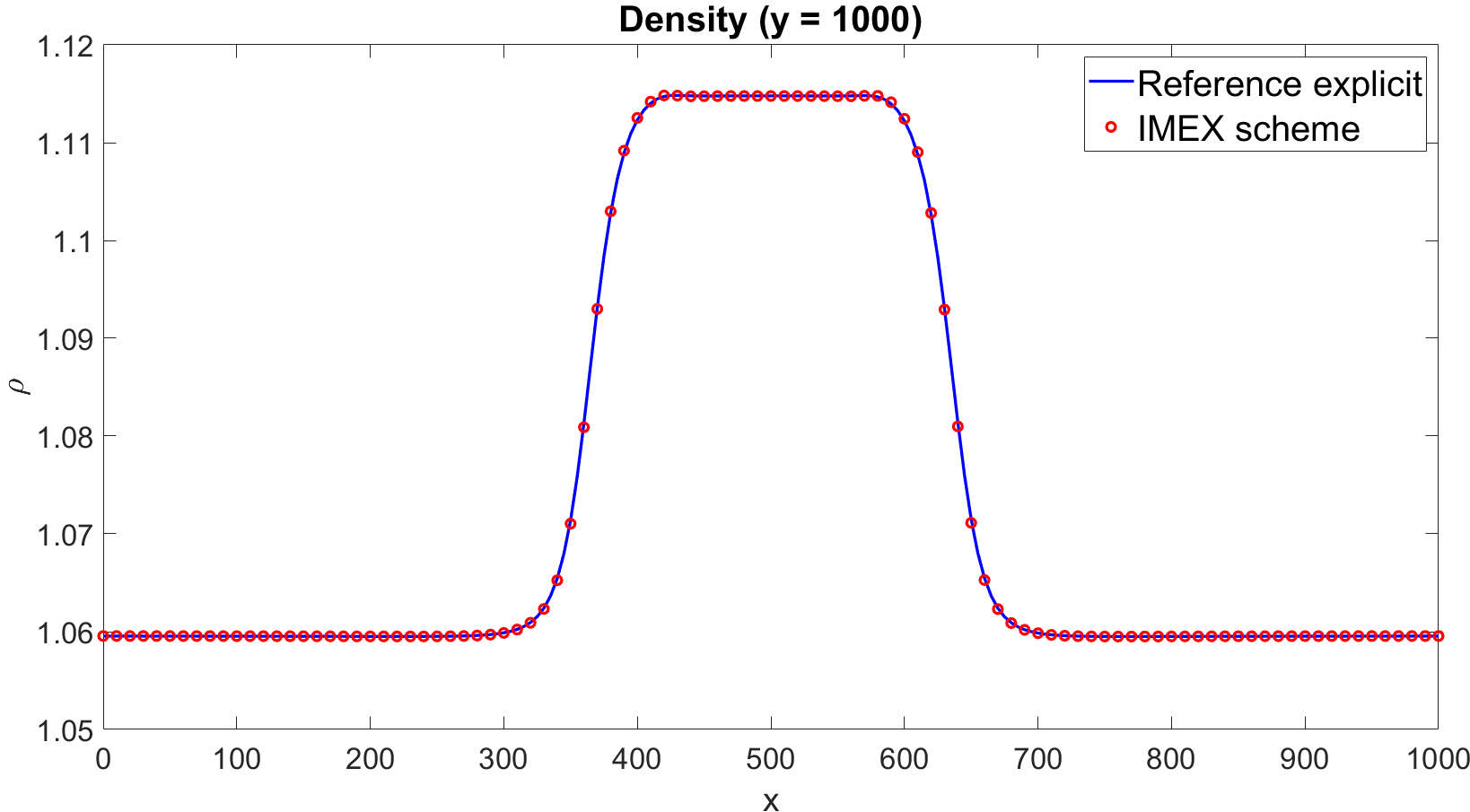}	
	\caption{Cold bubble test case, results at \(t = T_{f}\), density profile at \(y = 1000\). The continuous blue line represents the results for the reference explicit simulation, whereas the red dots denote the results for the IMEX scheme.} 
	\label{fig:Cold_bubble_ideal_rho}
\end{figure}

For the purpose of a quantitative comparison, a reference solution is computed with an explicit time discretization given by the optimal third order SSP scheme mentioned in Section \ref{ssec:Sod}. Figure \ref{fig:Cold_bubble_ideal} shows the contour plot of the potential temperature perturbation at \(t = T_{f}\) and one can easily notice that we are able to recover correctly the shape of the reference solution. For a more quantitative point of view, the profile of the density at \(y = 1000\) is reported in Figure \ref{fig:Cold_bubble_ideal_rho} and a good agreement between the reference results and those  obtained with the IMEX scheme is established. The IMEX scheme allows to employ a time step \(40\) times larger compared to the fully explicit scheme with a computational saving  of around 90\%. Three fixed-point iterations were required on average for each IMEX stage.  

\begin{figure}[H]
	\centering
	\begin{subfigure}{0.5\textwidth}
		\includegraphics[width=0.9\textwidth]{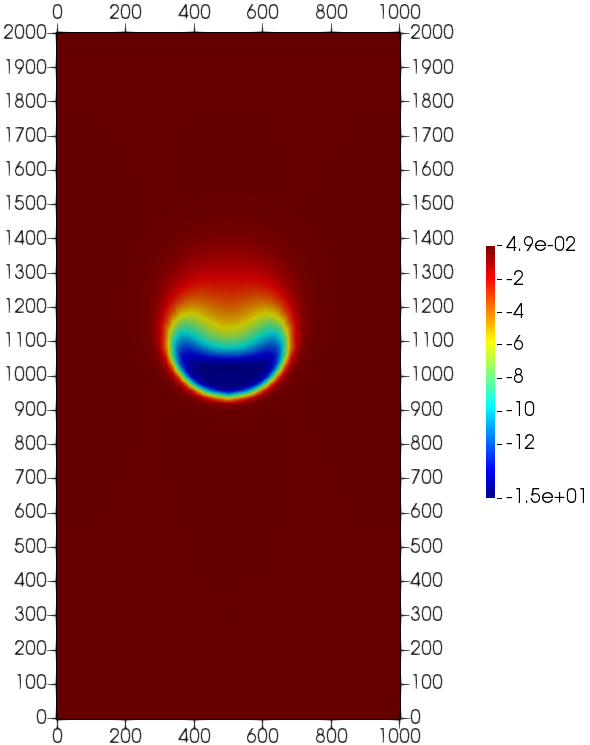} a)
	\end{subfigure}\hspace*{\fill}
	\begin{subfigure}{0.5\textwidth}
		\includegraphics[width=0.7\textwidth]{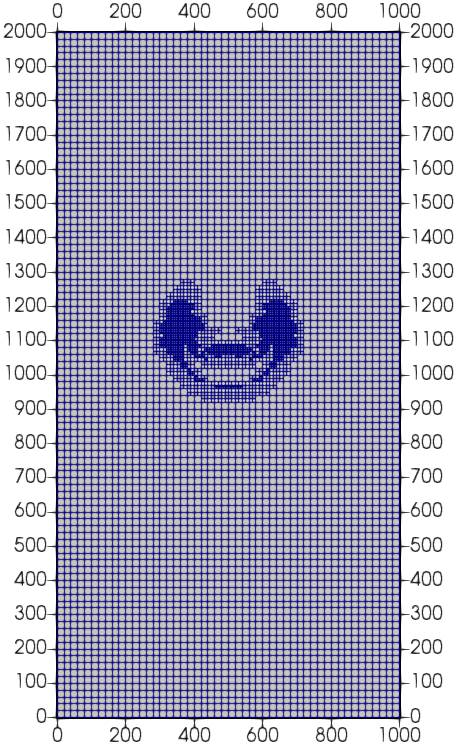} b)
	\end{subfigure}
	\caption{Cold bubble test case, adaptive simulation, results at \(t = T_{f}\), a) contour plot of potential temperature perturbation, b) computational grid.} 
	\label{fig:Cold_bubble_ideal_adaptive}
\end{figure}

\begin{figure}[H]
	\centering
	\includegraphics[width = 0.8\textwidth]{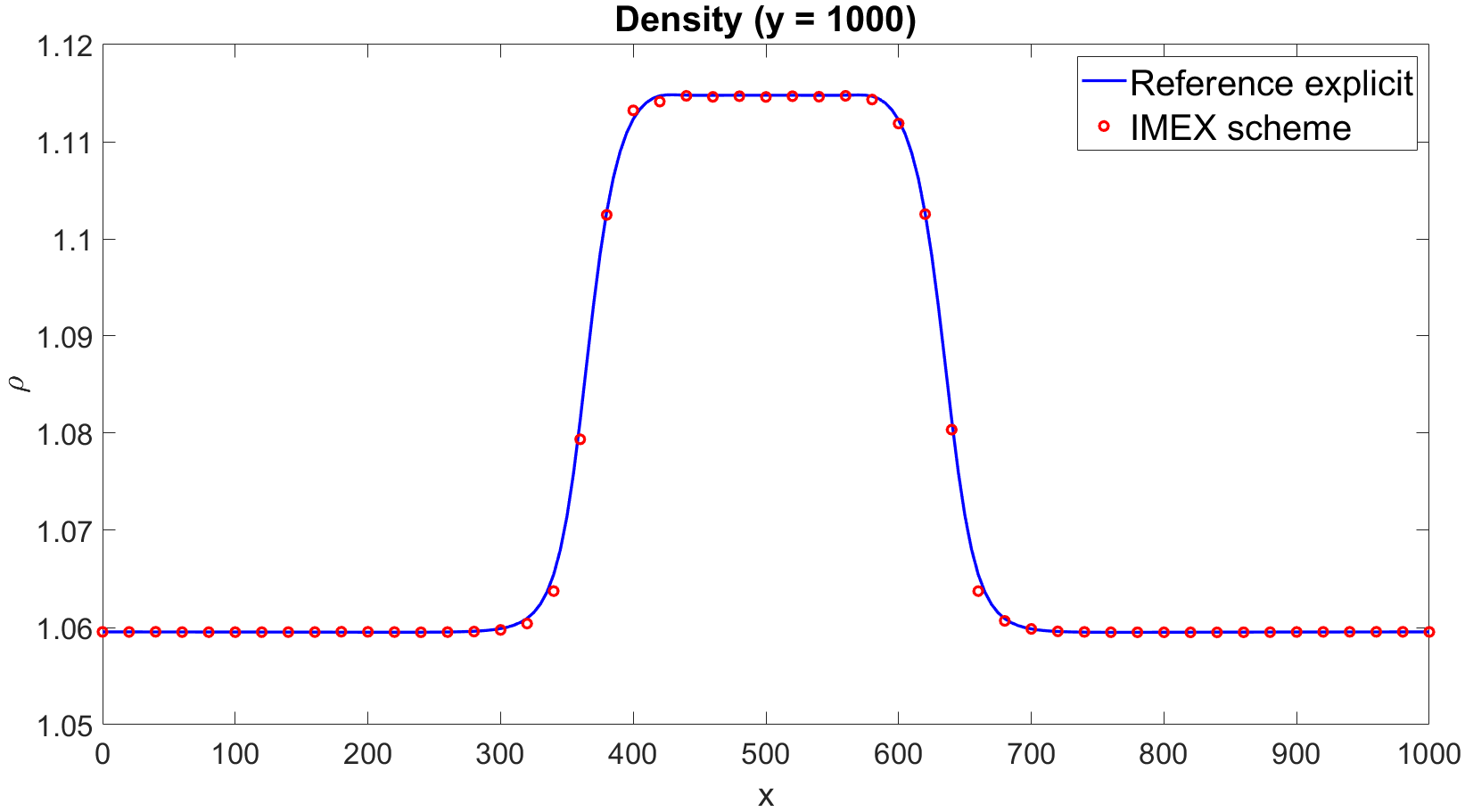}	
	\caption{Cold bubble test case, adaptive simulation, results at \(t = T_{f}\), density profile at \(y = 1000\). The continuous blue line represents the results for the reference explicit simulation, whereas the red dots denote the results for the IMEX scheme.} 
	\label{fig:Cold_bubble_ideal_rho_adaptive}
\end{figure}

In order to further enhance the computational efficiency, we employ again the code $h-$adaptivity capabilities. As mentioned in Section \ref{sec:EOS}, we use as refinement indicator the gradient of the potential temperature, since this quantity allows to identify the cold bubble. More specifically, we set

\begin{equation}
\eta_K = \max_{i \in \mathcal{N}_{K}} \left|\nabla\theta\right|_{i} 
\end{equation}
as local indicator, where \(\mathcal{N}_{K}\) is defined as in \eqref{eq:vortex_indicator}, and we allow to refine when \(\eta_K\) exceeds \(10^{-1}\) and to coarsen below \(6 \cdot 10^{-2}\). The initial computational grid is composed by \(50\times 100\) elements and we allow up to two local refinements only, so as to keep the advective Courant number under control and to recover the same maximum resolution employed for the non adaptive mesh simulation. As one can easily notice from Figure \ref{fig:Cold_bubble_ideal_adaptive}, the refinement criterion is able to track the bubble and the one-dimensional density profile  at \(y = 1000\)  in Figure \ref{fig:Cold_bubble_ideal_rho_adaptive} is correctly reproduced. The final mesh consists of 6914 elements instead of the 80000 elements of the full resolution mesh and a further 50\% reduction in computational time is achieved.

\begin{figure}[H]
	\centering
	\includegraphics[width = 0.8\textwidth]{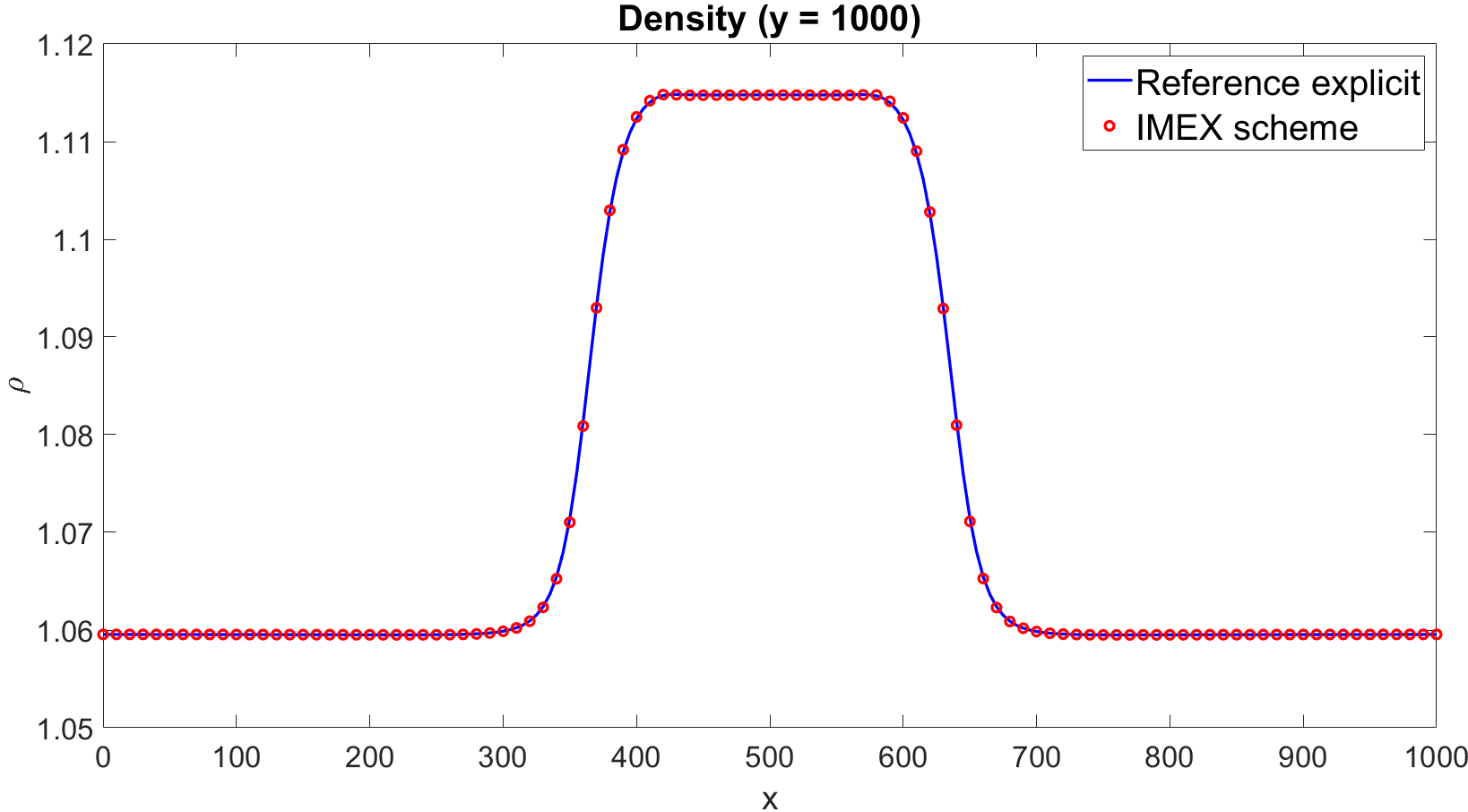}	
	\caption{Cold bubble test case, van der Waals EOS with \(\tilde a = 5 \cdot 10^{-9}\) and \(\tilde b = 5 \cdot 10^{-4}\), results at \(t = T_{f}\), density profile at \(y = 1000\). The continuous blue line represents the results for the reference explicit simulation, whereas the red dots denote the results for the IMEX scheme.} 
	\label{fig:Cold_bubble_vdW_z1_rho}
\end{figure}

We repeat now the same test using non-ideal equations of state. We first consider the van der Waals equation with a constant \(\tilde c_{v}\) given by \(\tilde c_{v} = \frac{\tilde R_{g}}{\gamma - 1} = 7.175 \cdot 10^{-3}, \tilde a = 5 \cdot 10^{-9}\) and \(\tilde b = 5 \cdot 10^{-4}\), so that the same specific heat at constant volume with respect to the ideal gas case is obtained and \(z \approx 1\). The fluid is initialized using the same pressure and the same density values as in the ideal gas case. Notice that \(\frac{d \tilde a}{dT} = \frac{d \tilde c_{v}}{dT} = 0\) and so it is not necessary to explicitly compute the temperature for \eqref{eq:nlineq1_cubic} and \eqref{eq:nlineq2_cubic}. We expect a behaviour similar to the ideal gas one and this is confirmed by the density profile reported in Figure \ref{fig:Cold_bubble_vdW_z1_rho}. 

We then consider the case with \(\tilde a = 1.6 \cdot 10^{-1}\) and \(\tilde b = 5 \cdot 10^{-4}\), which yield an average compressibility factor \(z \approx 0.83\). In this case, we expect effects due to conditions far from the ideal ones and so we first compute a reference solution with the explicit time discretization. The time step for the IMEX simulation is kept equal to \(\Delta t = 0.08\), yielding a maximum Courant number \(C \approx 5.3\) and a maximum advective Courant number \(C_{u} \approx 0.19\). Figure \ref{fig:Cold_bubble_vdW_z0,83} shows the contour plot for \(\beta\) at \(t = T_{f}\) for both the reference explicit and the IMEX simulations. The expected behaviour is retrieved and a good agreement  with the reference results is established. Also in this case, a computational saving of around 90\% with respect to the explicit simulation is obtained thanks to the IMEX scheme. Figure \ref{fig:Cold_bubble_vdW_z0,83_rho} reports the profile of the density for \(y = 1000\) at \(t = T_{f}\). One can notice the very good agreement  between the IMEX results and the reference ones. Furthermore, a clear discrepancy with respect to the ideal gas can be observed. The higher density values are due to the large value of \(\tilde{a}\), which means that strong forces of attraction between the gas particles are present \cite{nederstigt:2017}. 

\begin{figure}[H]
	\begin{subfigure}{0.5\textwidth}
		\includegraphics[width=0.9\textwidth]{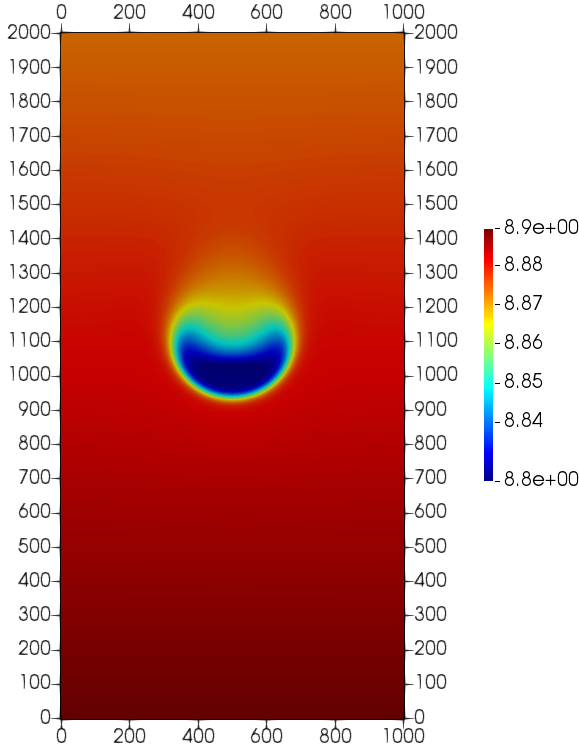} a)
	\end{subfigure}\hspace*{\fill}
	\begin{subfigure}{0.5\textwidth}
		\includegraphics[width=0.9\textwidth]{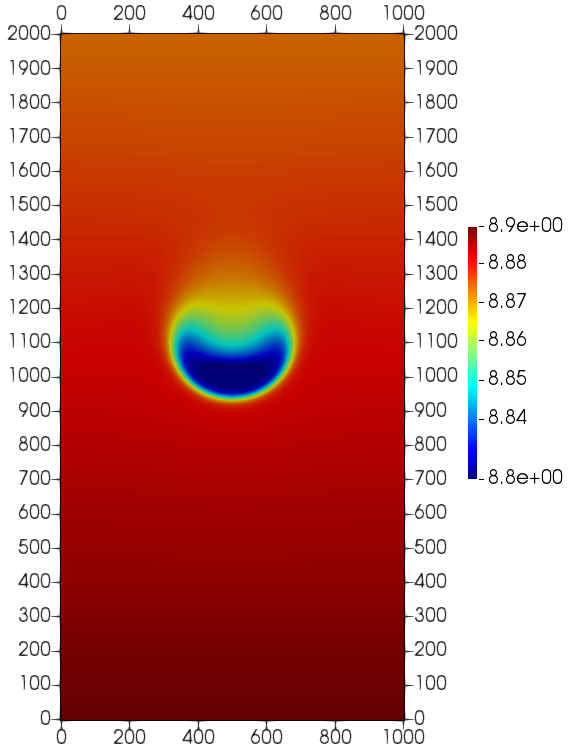} b)
	\end{subfigure}
	\caption{Cold bubble test case, van der Waals EOS with \(\tilde a = 1.6 \cdot 10^{-1}\) and \(\tilde b = 5 \cdot 10^{-4}\), results at \(t = T_{f}\), a) contour plot of \(\beta\) for the reference explicit simulation, b) contour plot of \(\beta\) for the simulation with IMEX scheme.} 
	\label{fig:Cold_bubble_vdW_z0,83}
\end{figure}

\begin{figure}[H]
	\centering
	\includegraphics[width = 0.75\textwidth]{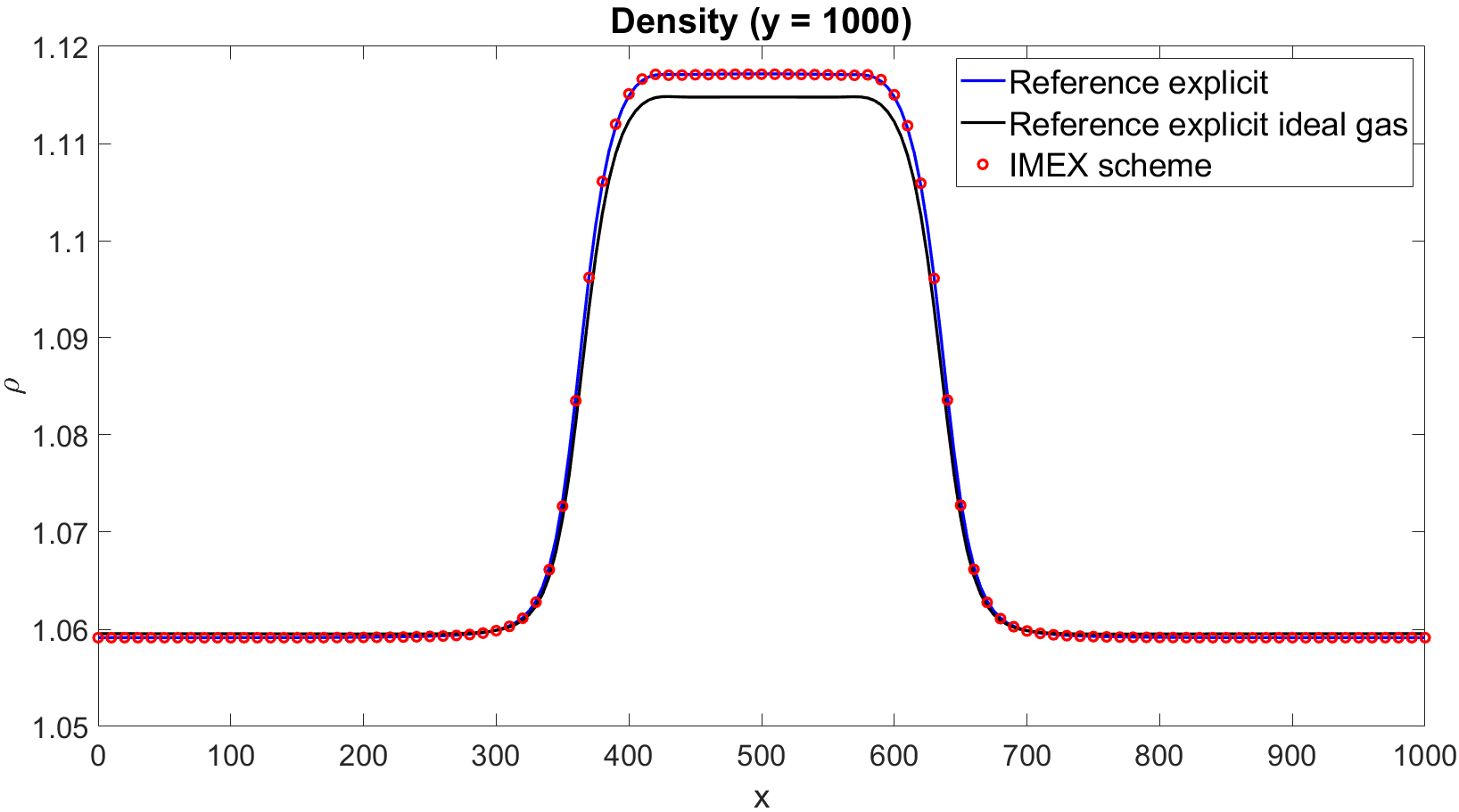}	
	\caption{Cold bubble test case, van der Waals EOS with \(\tilde a = 1.6 \cdot 10^{-1}\) and \(\tilde b = 5 \cdot 10^{-4}\), results at \(t = T_{f}\), density profile at \(y = 1000\). The continuous blue line represents the results for the full explicit simulation, the continuous black line reports the results for the reference explicit simulation with an ideal gas, whereas the red dots denote the results for the IMEX scheme.} 
	\label{fig:Cold_bubble_vdW_z0,83_rho}
\end{figure}

Concerning the adaptive simulations, since, as proven in Section \ref{sec:EOS}, the quantity \(\beta = \log(T) - 2\frac{\tilde R_{g}}{\tilde c_{v}}\atanh\left(2\rho \tilde b - 1\right)\) is constant in an isentropic process with \(\frac{d\tilde a}{dT} = \frac{d \tilde c_{v}}{dT} = 0\), we define the local refinement indicator for each element as

\begin{equation}
\eta_K = \max_{i \in \mathcal{N}_K} \left|\nabla\beta\right|_{i} .
\end{equation}
We allow to refine when \(\eta_{K}\) exceeds \(4 \cdot 10^{-4}\) and to coarsen when the indicator is below \(2 \cdot 10^{-4}\). The initial mesh is composed by \(50 \times 100\) elements and we allow up to four local refinements. For this reason, in order to keep under control the advective Courant number, we need to reduce the time step \(\Delta t = 0.02\), so as to obtain a maximum acoustic Courant number \(C \approx 5.3\) and a maximum advective Courant number \(C_{u} \approx 0.18\). Figure \ref{fig:Cold_bubble_vdW_z0,83_adaptive} confirms that \(\beta\) is an appropriate quantity to track the bubble and the one-dimensional density profile in Figure \ref{fig:Cold_bubble_vdW_z0,83_adaptive_rho} shows that no significant loss in accuracy occurs. The final mesh consists of 9086 elements. 

\begin{figure}[H]
	\centering
	\begin{subfigure}{0.5\textwidth}
		\includegraphics[width=0.9\textwidth]{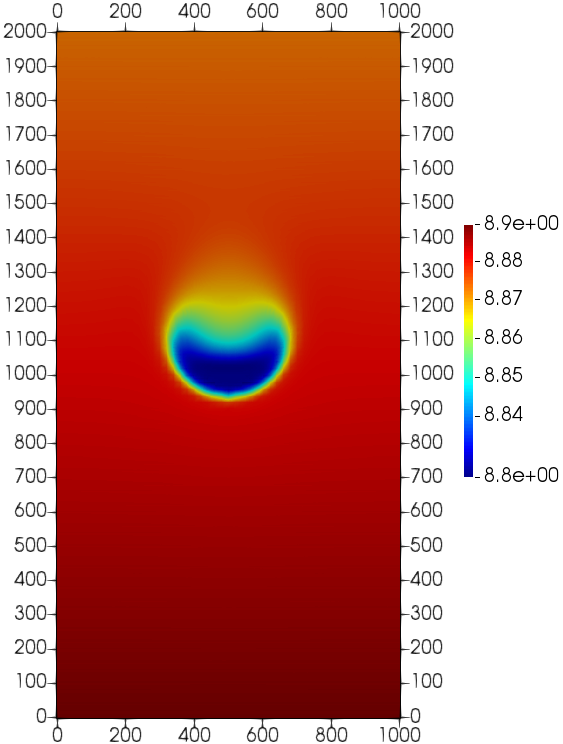} a)
	\end{subfigure}\hspace*{\fill}
	\begin{subfigure}{0.5\textwidth}
		\includegraphics[width=0.7\textwidth]{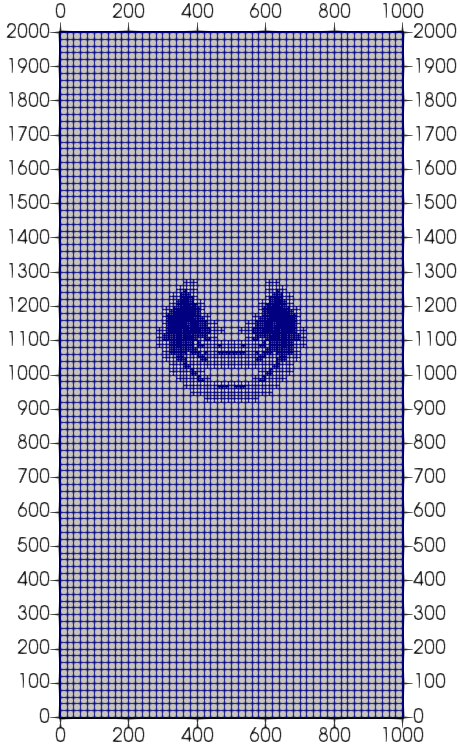} b)
	\end{subfigure}
	\caption{Cold bubble test case, van der Waals EOS with \(\tilde a = 1.6 \cdot 10^{-1}\) and \(\tilde b = 5 \cdot 10^{-4}\), adaptive simulation, results at \(t = T_{f}\), a) contour plot of \(\beta\), b) computational grid.} 
	\label{fig:Cold_bubble_vdW_z0,83_adaptive}
\end{figure}

\begin{figure}[H]
	\centering
	\includegraphics[width = 0.8\textwidth]{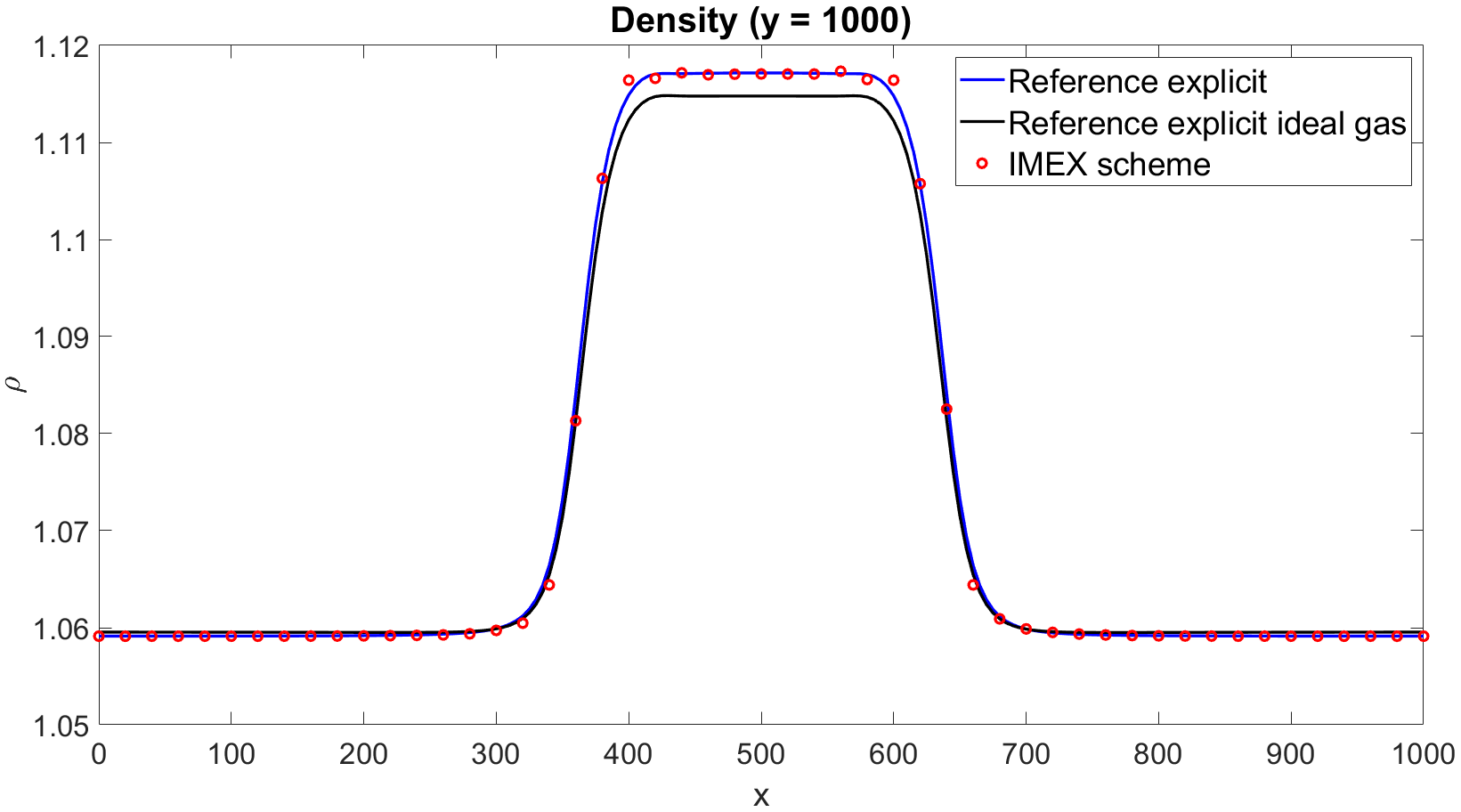}
	\caption{Cold bubble test case, van der Waals EOS with \(\tilde a = 1.6 \cdot 10^{-1}\) and \(\tilde b = 5 \cdot 10^{-4}\), adaptive simulation, results at \(t = T_{f}\), density profile at \(y = 1000\). The continuous blue line represents the results for the reference explicit simulation, the continuous black line reports the results for the reference explicit simulation with an ideal gas, whereas the red dots denote the results for the IMEX scheme in the non-ideal case.} 
	\label{fig:Cold_bubble_vdW_z0,83_adaptive_rho}
\end{figure}

\begin{figure}[H]
	\centering
	\includegraphics[width = 0.75\textwidth]{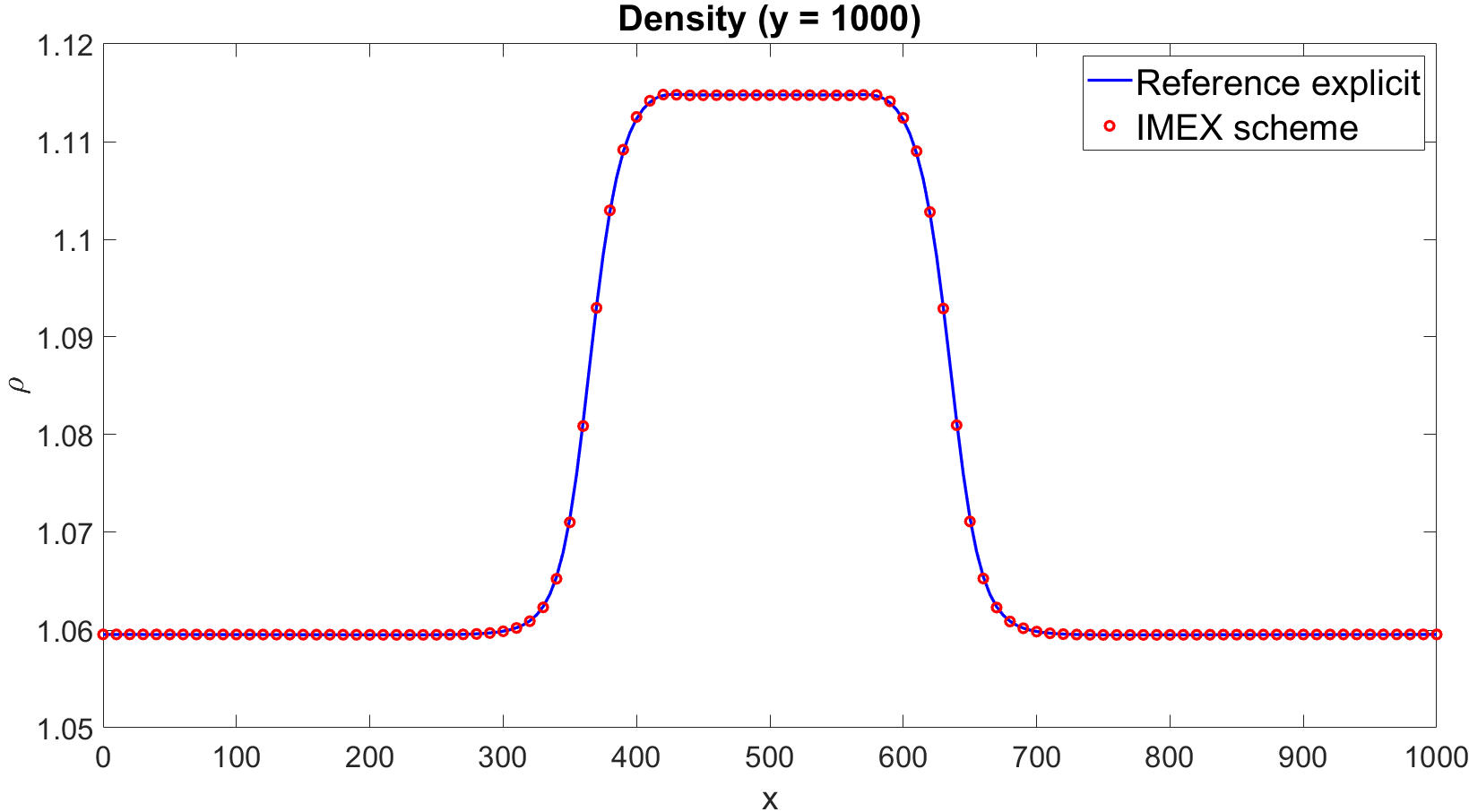}	
	\caption{Cold bubble test case, Peng-Robinson EOS with \(\tilde a = 5 \cdot 10^{-9}\) and \(\tilde b = 5 \cdot 10^{-4}\), results at \(t = T_{f}\), density profile at \(y = 1000\). The continuous blue line represents the results for the reference explicit simulation, whereas the red dots denote the results for the IMEX scheme.} 
	\label{fig:Cold_bubble_PR_z1_rho}
\end{figure}

The same analysis is carried out using the Peng-Robinson EOS. Hence, we first consider \(\tilde R_{g} = 2.87 \cdot 10^{-3}\), \(\tilde c_{v} = 7.175 \cdot 10^{-3}\), \(\tilde a = 5 \cdot 10^{-9}\) and \(\tilde b = 5 \cdot 10^{-4}\), so that \(z \approx 1\). The density profile reported in Figure \ref{fig:Cold_bubble_PR_z1_rho} highlights, as expected, a behaviour entirely analogous   to that of the ideal gas. 

Next, we take \(\tilde a = 1.6 \cdot 10^{-1}\) and \(\tilde b = 5 \cdot 10^{-4}\), so that \(z \approx 0.83\), and we perform both uniform mesh and adaptive simulations, using the same parameters employed for the van der Waals EOS. The results are compared with a reference solution computed with the explicit method. Figure \ref{fig:Cold_bubble_PR_z0,83} shows similar contour plots for all the configurations as well as for the adaptive mesh at \(t = T_{f}\), which consists of \(9137\) elements and it is able to track the bubble correctly. Figure \ref{fig:Cold_bubble_PR_z0,83_rho} reports the comparison for the one-dimensional profile of the density at \(y = 1000\) and the same considerations of the van der Waals EOS are still valid. 
We want to test in this case the refinement indicator based on \eqref{eq:beta_1}. More specifically, we set

\begin{equation}\label{eq:indicator_1}
\eta_{K} = \max_{i \in \mathcal{N}_K} \left|\nabla\left(\frac{p}{\rho^{\gamma_{p\rho}}}\right)\right|_{i}
\end{equation} 
and we allow to refine in case \(\eta_{K}\) is above \(4 \cdot 10^{-4}\) and to coarsen below \(2 \cdot 10^{-4}\) with the same remeshing procedure adopted so far for non-ideal gases. Figure \ref{fig:Cold_bubble_PR_z0,83_adaptive_bis} shows the contour plot of \eqref{eq:beta_1} and the computational mesh at \(t = T_{f}\). The mesh consists of \(8294\) elements and one can easily notice that more resolution is added only in correspondence of the bubble. 

\begin{figure}[H]
	\centering
	\begin{subfigure}{0.5\textwidth}
		\includegraphics[width=0.9\textwidth]{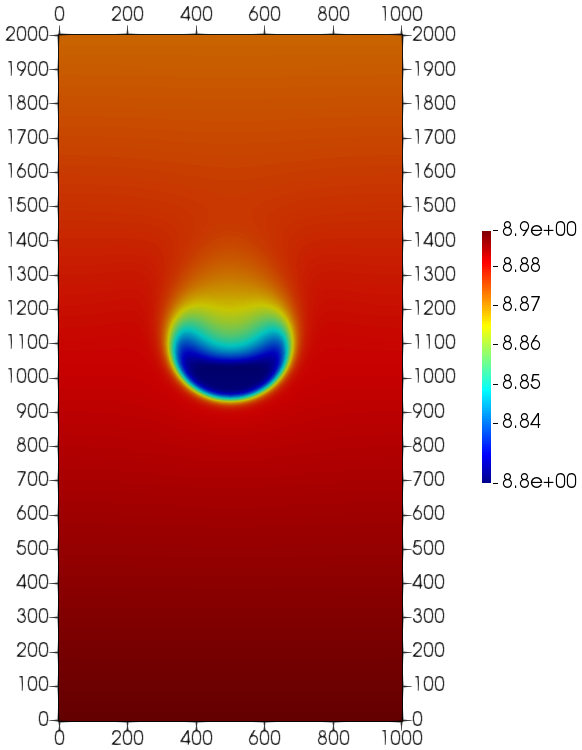} a)
	\end{subfigure}\hspace*{\fill}
	\begin{subfigure}{0.5\textwidth}
		\includegraphics[width=0.9\textwidth]{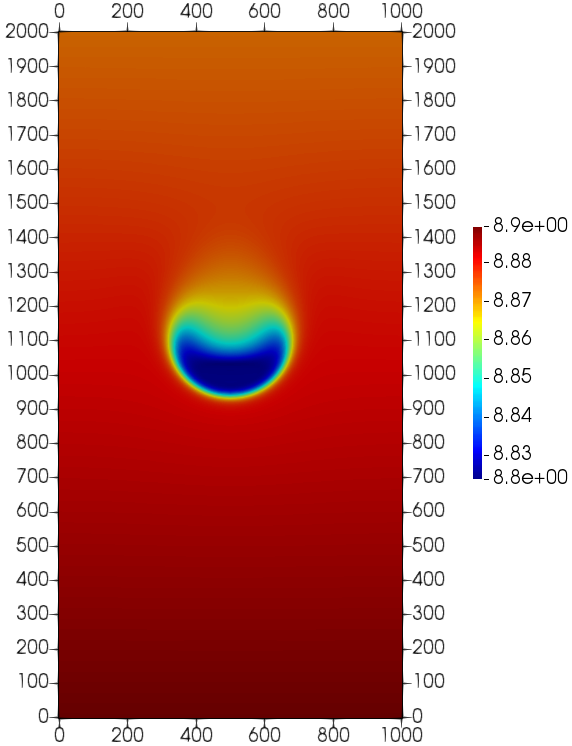} b)
	\end{subfigure}
	\begin{subfigure}{0.5\textwidth}
		\includegraphics[width=0.9\textwidth]{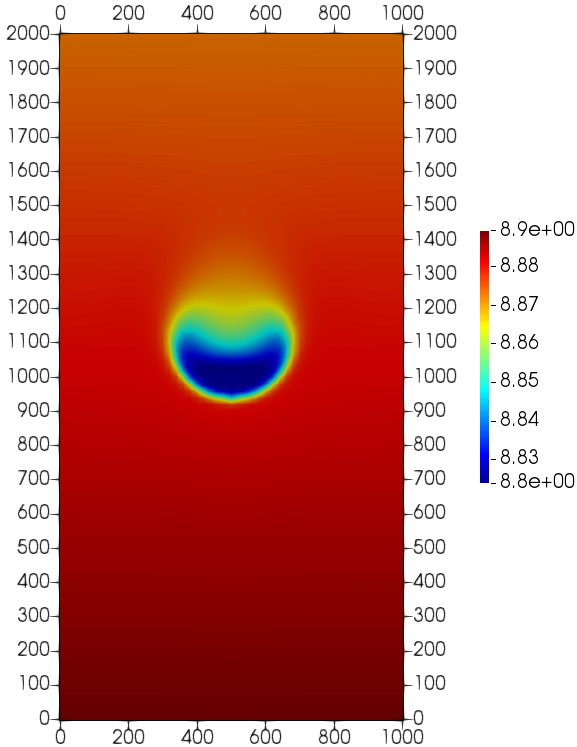} a)
	\end{subfigure}\hspace*{\fill}
	\begin{subfigure}{0.5\textwidth}
		\includegraphics[width=0.7\textwidth]{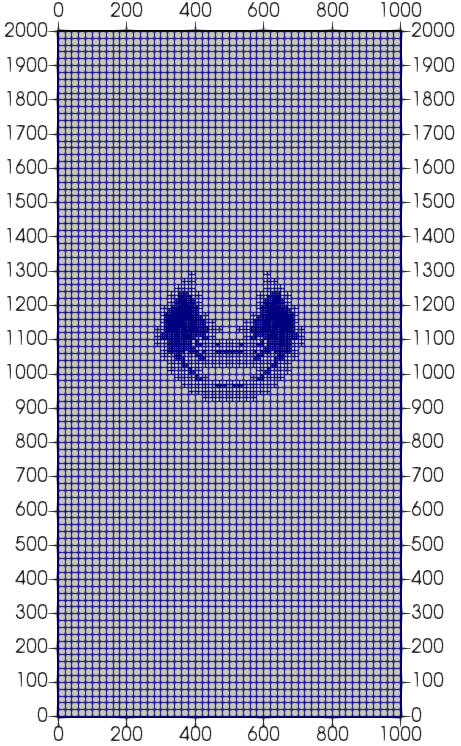} b)
	\end{subfigure}
	\caption{Cold bubble test case, Peng-Robinson EOS with \(\tilde a = 1.6 \cdot 10^{-1}\) and \(\tilde b = 5 \cdot 10^{-4}\), results at \(t = T_{f}\), a) contour plot of \(\beta\) for the reference explicit simulation, b) contour plot of \(\beta\) for the constant mesh simulation with IMEX scheme, c) contour plot of \(\beta\) for adaptive simulation with IMEX scheme, d) adaptive mesh.} 
	\label{fig:Cold_bubble_PR_z0,83}
\end{figure}

\begin{figure}[H]
	\centering
	\includegraphics[width = 0.8\textwidth]{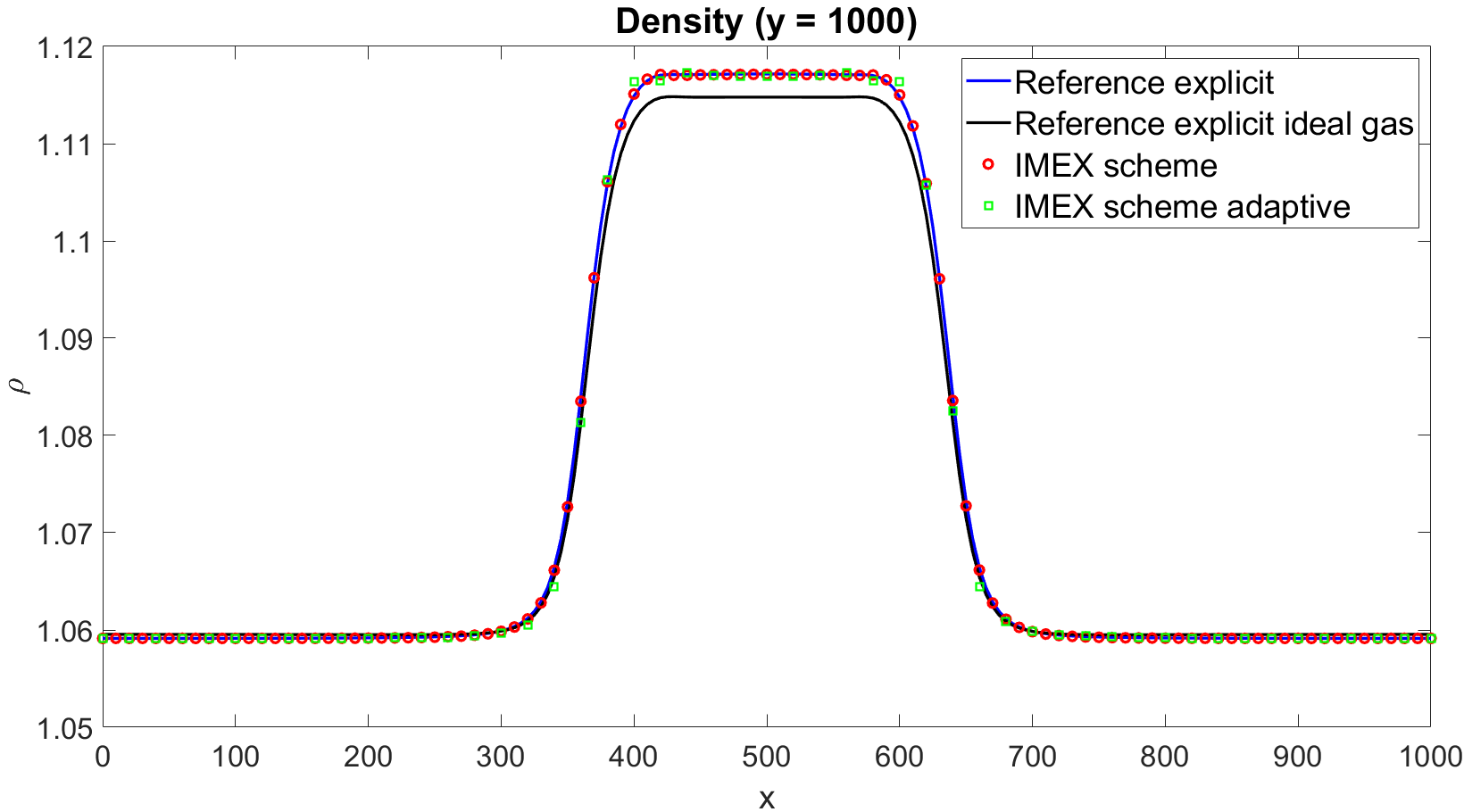}	
	\caption{Cold bubble test case, Peng-Robinson EOS with \(\tilde a = 1.6 \cdot 10^{-1}\) and \(\tilde b = 5 \cdot 10^{-4}\), adaptive simulation, results at \(t = T_{f}\), density profile at \(y = 1000\). The continuous blue line represents the results for the reference explicit simulation, the continuous black line reports the results for the reference explicit simulation with an ideal gas, the red dots denote the results for the IMEX scheme, whereas the green squares represent the results for the adaptive simulation with IMEX scheme.} 
	\label{fig:Cold_bubble_PR_z0,83_rho}
\end{figure}

\begin{figure}[H]
	\centering
	\begin{subfigure}{0.47\textwidth}
		\includegraphics[width=0.8\textwidth]{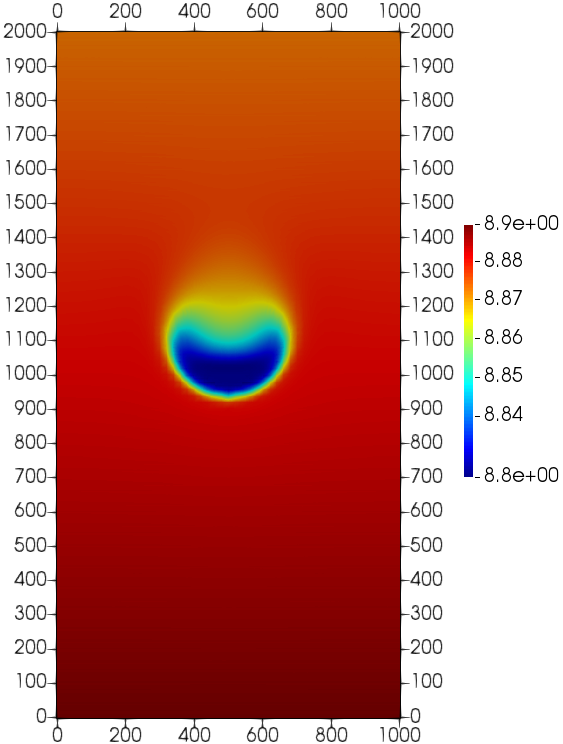} a)
	\end{subfigure}\hspace*{\fill}
	\begin{subfigure}{0.47\textwidth}
		\includegraphics[width=0.65\textwidth]{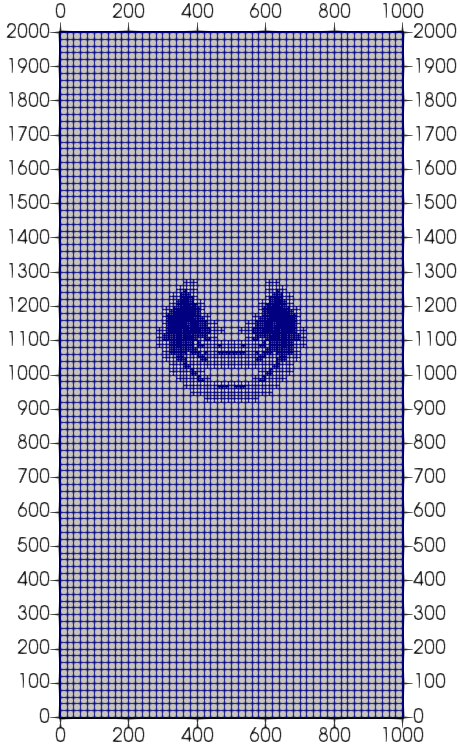} b)
	\end{subfigure}
	\caption{Cold bubble test case, Peng-Robinson EOS with \(\tilde a = 1.6 \cdot 10^{-1}\) and \(\tilde b = 5 \cdot 10^{-4}\), adaptive simulation with criterion \eqref{eq:indicator_1}, results at \(t = T_{f}\), a) contour plot of \eqref{eq:beta_1}, b) adaptive mesh} 
	\label{fig:Cold_bubble_PR_z0,83_adaptive_bis}
\end{figure}

\subsection{Warm bubble}
\label{ssec:warm_bubble}

\begin{figure}[H]
	\begin{subfigure}{0.45\textwidth}
		\centering
		\includegraphics[width=\textwidth]{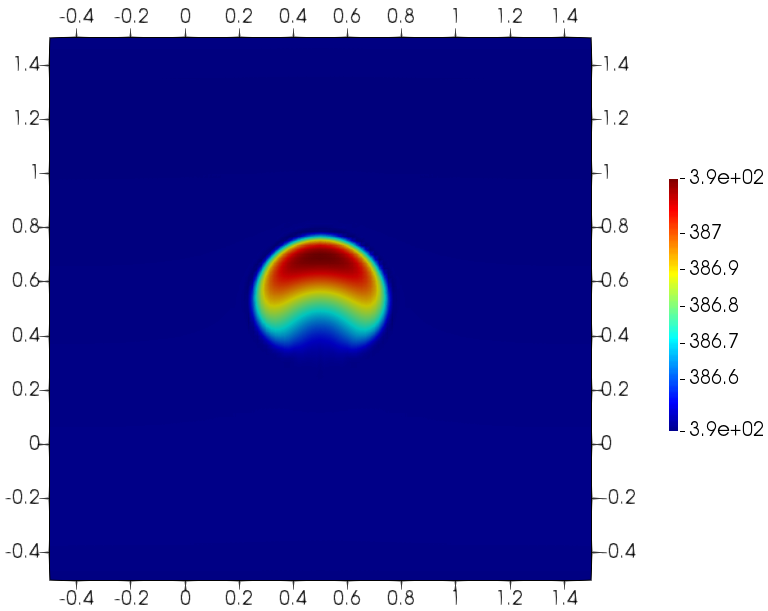}
	\end{subfigure}\hspace*{\fill}
	\begin{subfigure}{0.45\textwidth}
		\centering
		\includegraphics[width=\textwidth]{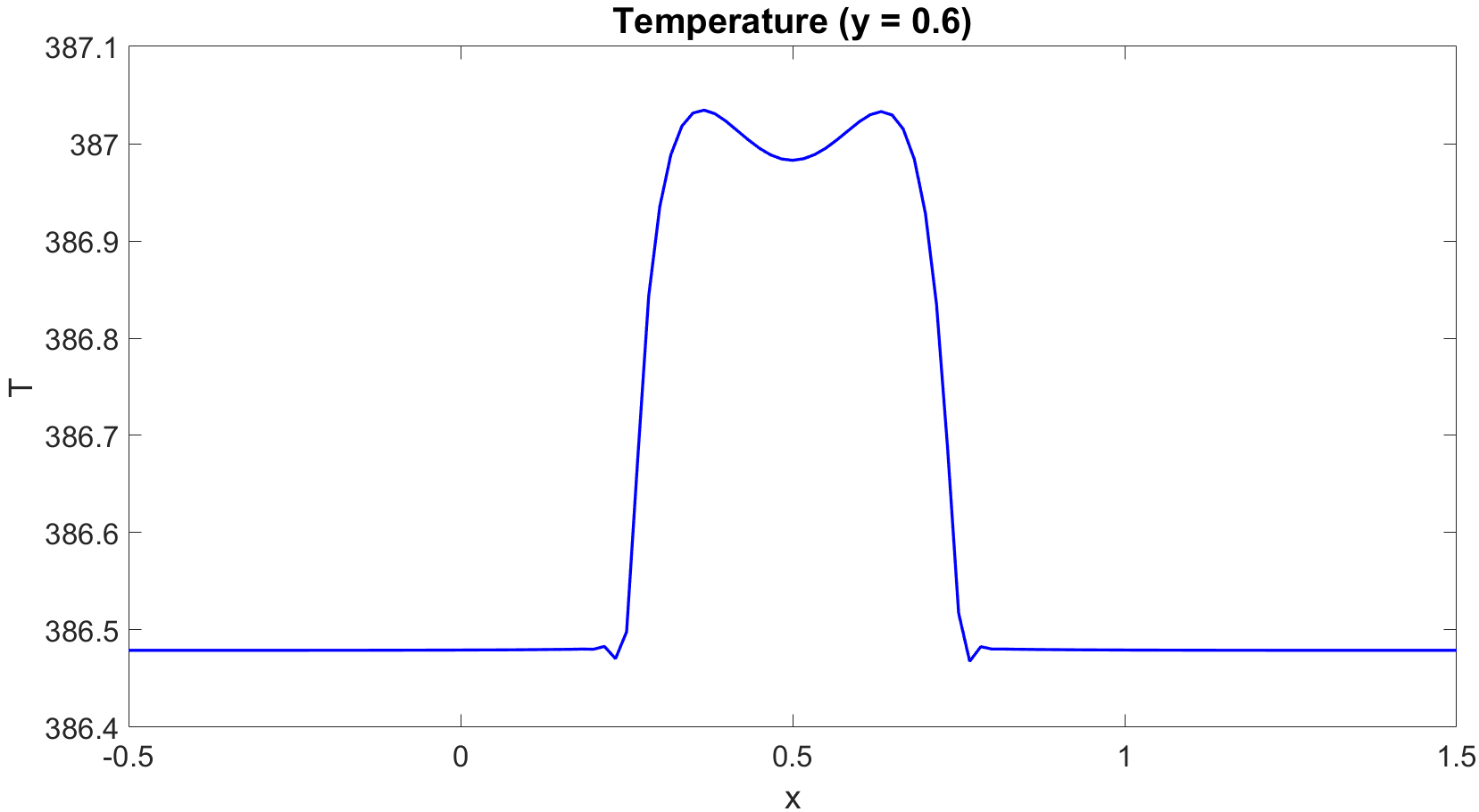}
	\end{subfigure}	
	\begin{subfigure}{0.45\textwidth}
		\centering
		\includegraphics[width=\textwidth]{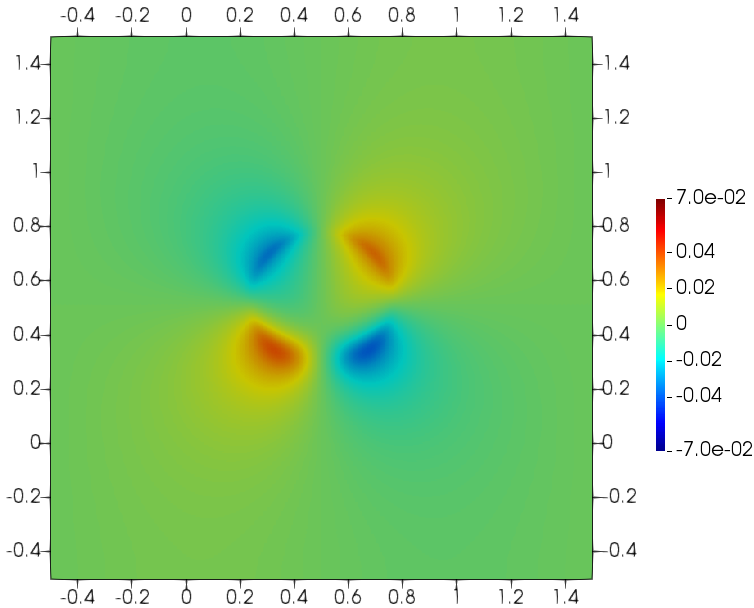}
	\end{subfigure}\hspace*{\fill}
	\begin{subfigure}{0.45\textwidth}
		\centering
		\includegraphics[width=\textwidth]{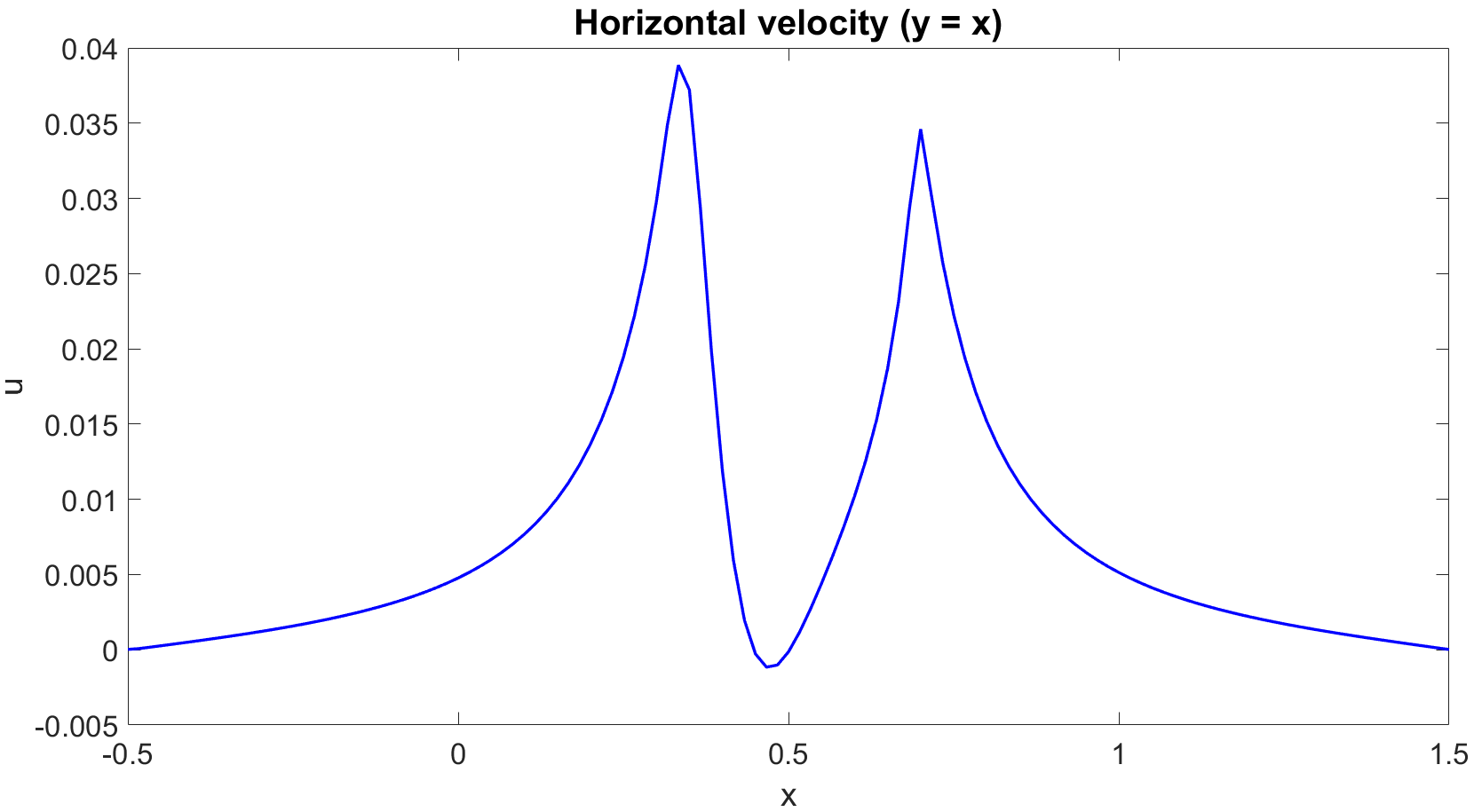}
	\end{subfigure}	
	\begin{subfigure}{0.45\textwidth}
		\centering
		\includegraphics[width=\textwidth]{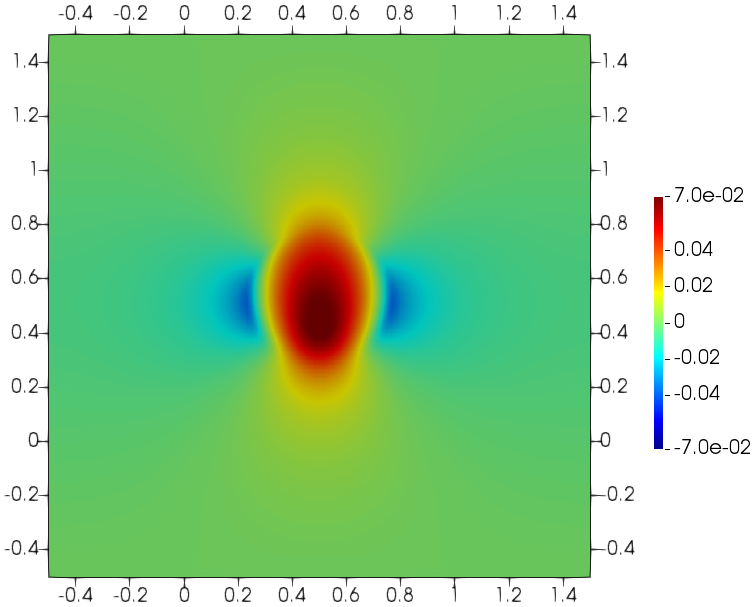}
	\end{subfigure}\hspace*{\fill}
	\begin{subfigure}{0.45\textwidth}
		\centering
		\includegraphics[width=\textwidth]{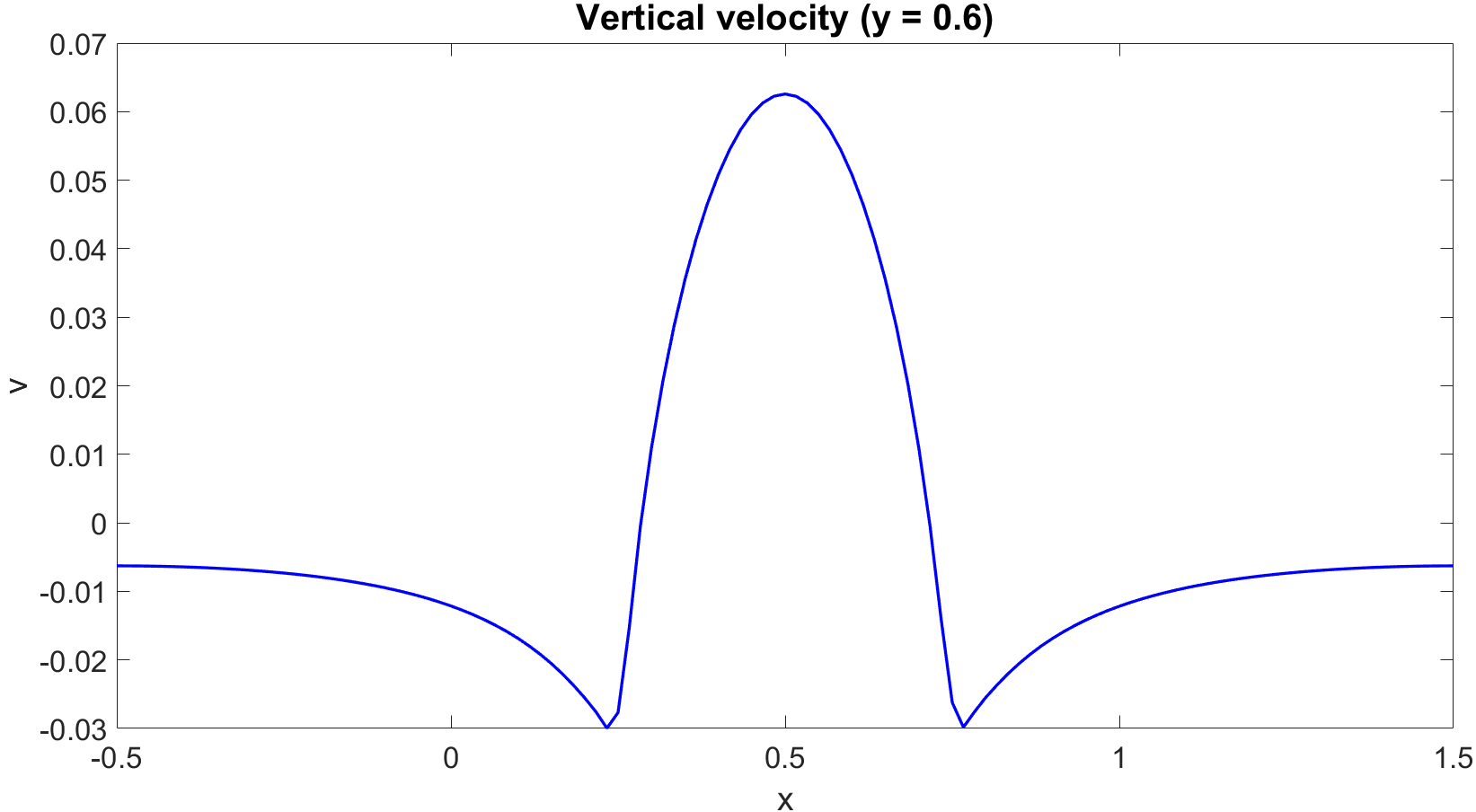}
	\end{subfigure}	
	\caption{Warm bubble test case, results at \(t = \SI{10}{\second}\). From bottom to top: temperature, horizontal velocity and vertical velocity.}
	\label{fig:Warm_bubble_t10}
\end{figure}

\begin{figure}[H]
	\begin{subfigure}{0.45\textwidth}
		\centering
		\includegraphics[width=\textwidth]{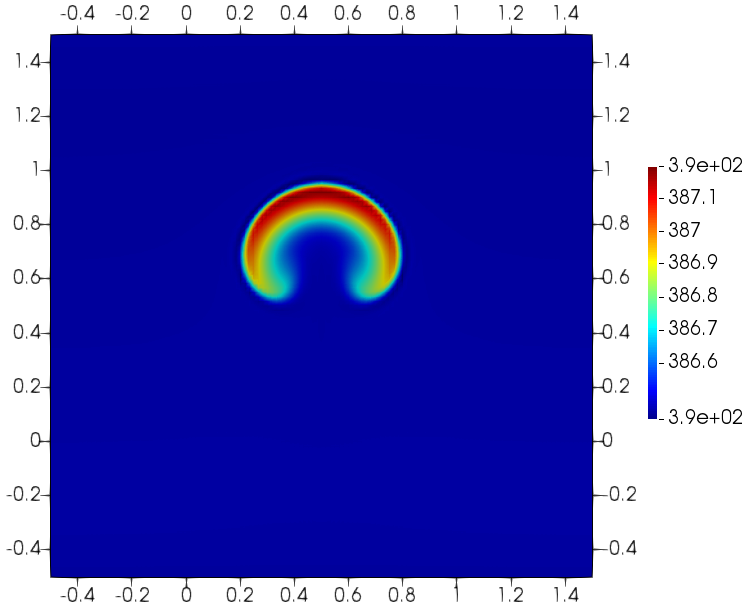}
	\end{subfigure}\hspace*{\fill}
	\begin{subfigure}{0.45\textwidth}
		\centering
		\includegraphics[width=\textwidth]{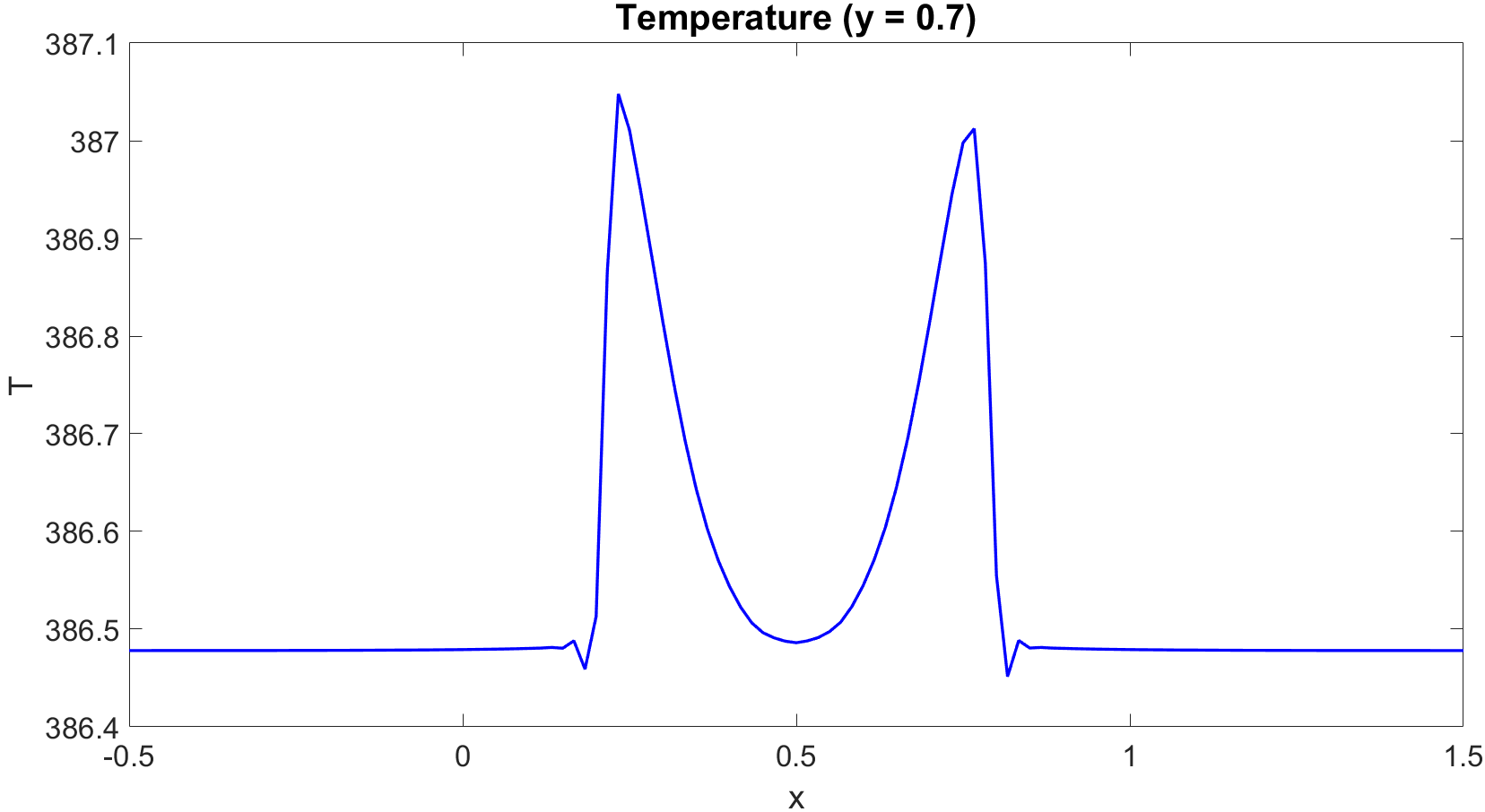}
	\end{subfigure}	
	\begin{subfigure}{0.45\textwidth}
		\centering
		\includegraphics[width=\textwidth]{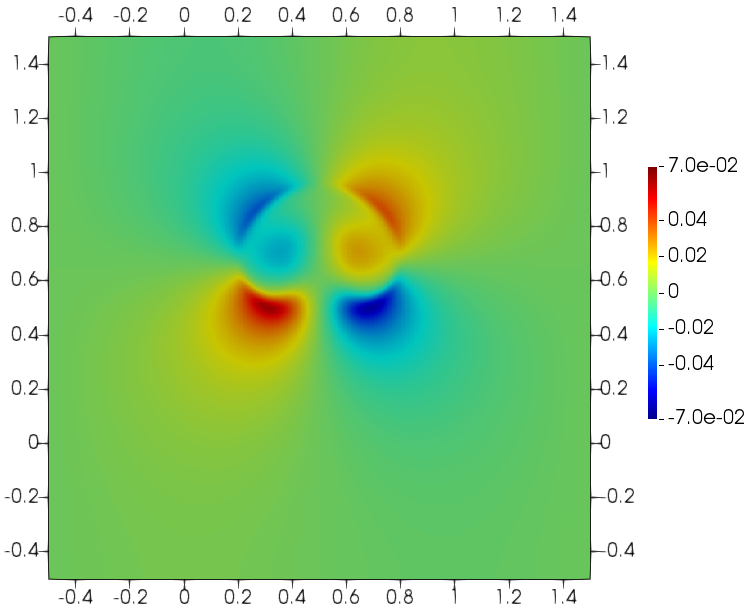}
	\end{subfigure}\hspace*{\fill}
	\begin{subfigure}{0.45\textwidth}
		\centering
		\includegraphics[width=\textwidth]{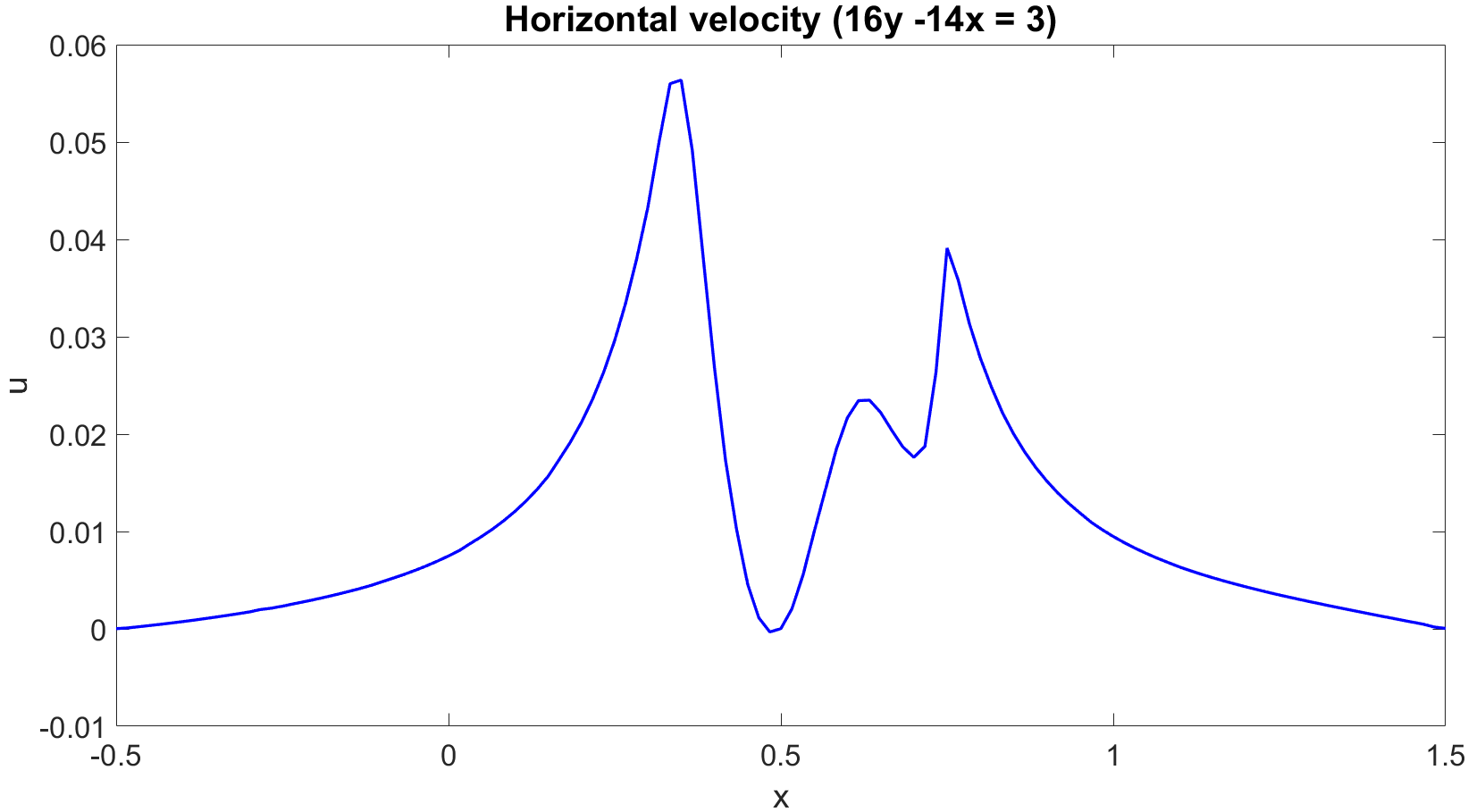}
	\end{subfigure}	
	\begin{subfigure}{0.45\textwidth}
		\centering
		\includegraphics[width=\textwidth]{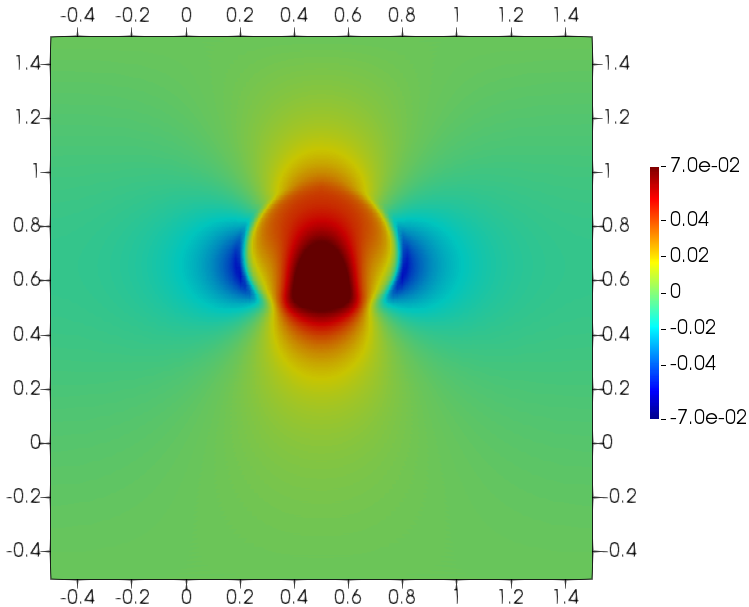}
	\end{subfigure}\hspace*{\fill}
	\begin{subfigure}{0.45\textwidth}
		\centering
		\includegraphics[width=\textwidth]{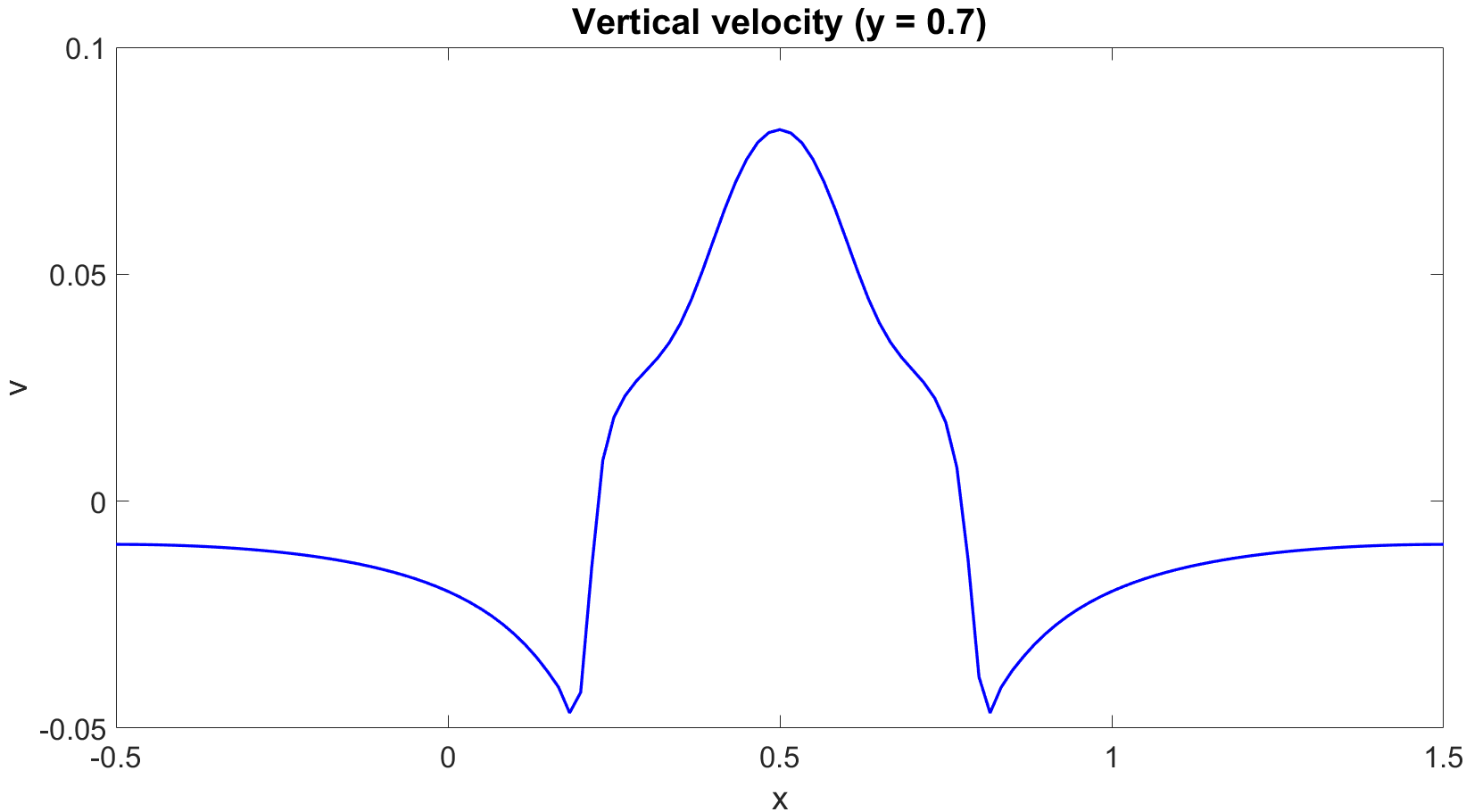}
	\end{subfigure}	
	\caption{Warm bubble test case, results at \(t = \SI{15}{\second}\). From bottom to top: temperature, horizontal velocity and vertical velocity.}
	\label{fig:Warm_bubble_t15}
\end{figure}

\begin{figure}[H]
	\begin{subfigure}{0.45\textwidth}
		\centering
		\includegraphics[width=\textwidth]{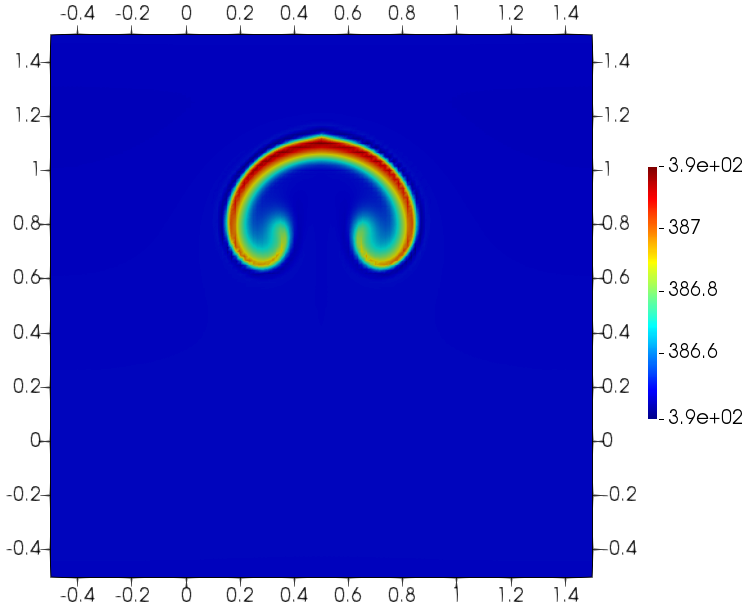}
	\end{subfigure}\hspace*{\fill}
	\begin{subfigure}{0.45\textwidth}
		\centering
		\includegraphics[width=\textwidth]{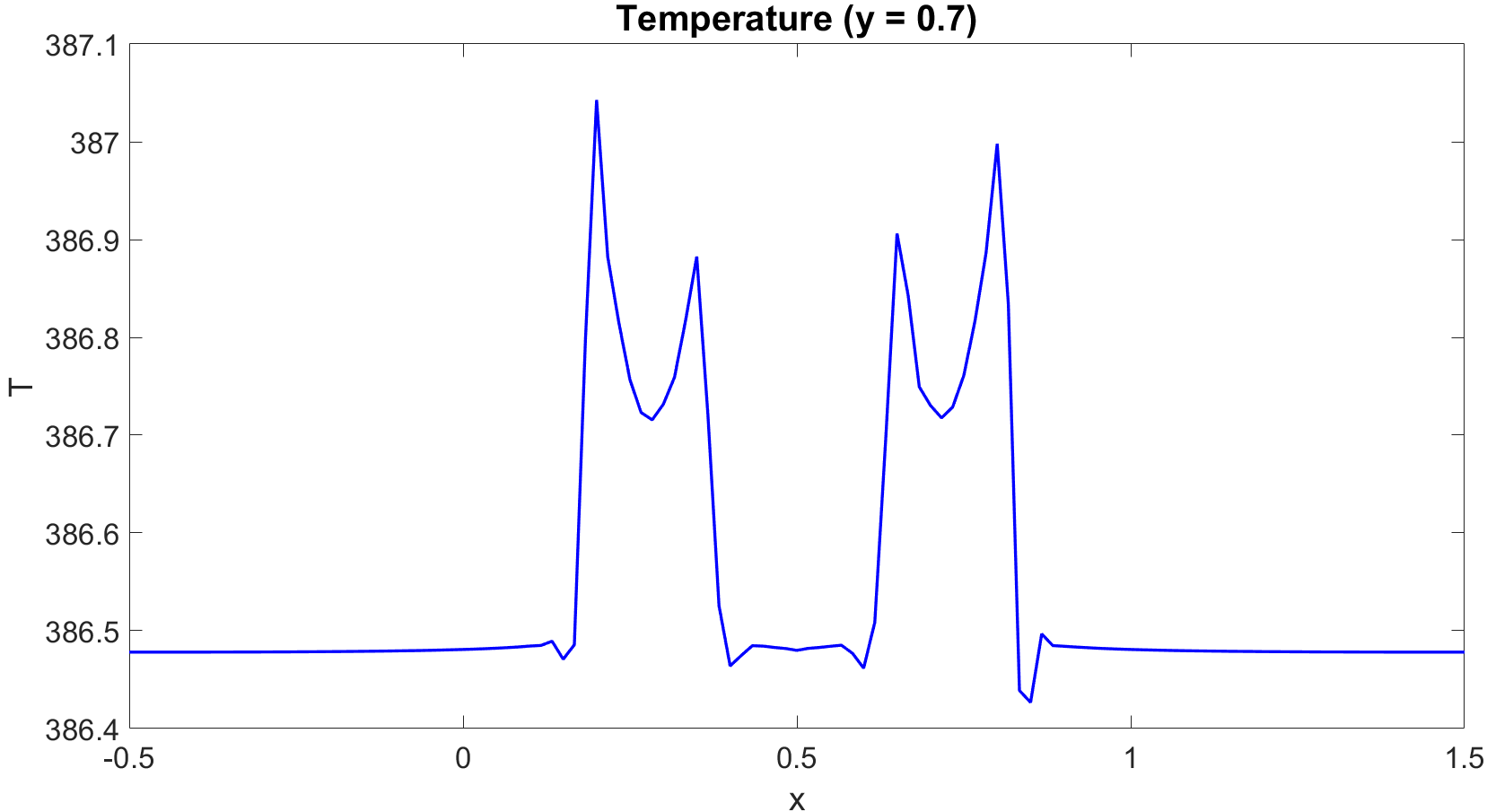}
	\end{subfigure}	
	\begin{subfigure}{0.45\textwidth}
		\centering
		\includegraphics[width=\textwidth]{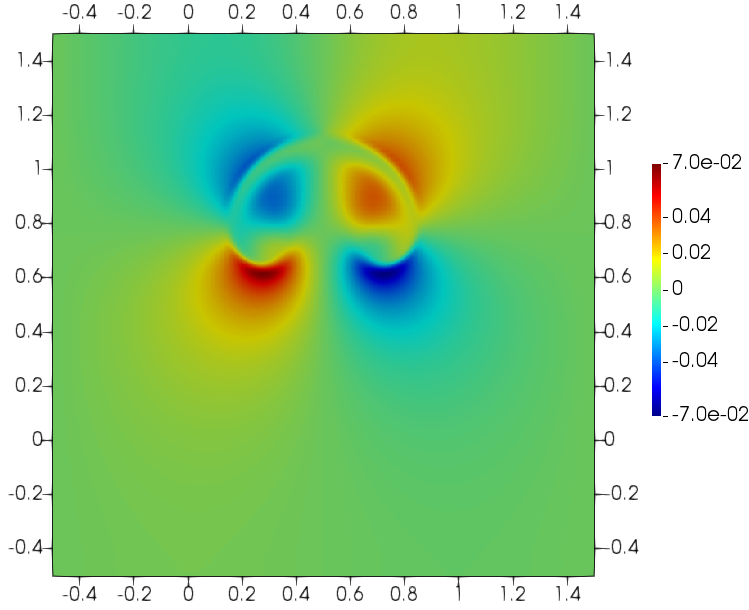}
	\end{subfigure}\hspace*{\fill}
	\begin{subfigure}{0.45\textwidth}
		\centering
		\includegraphics[width=\textwidth]{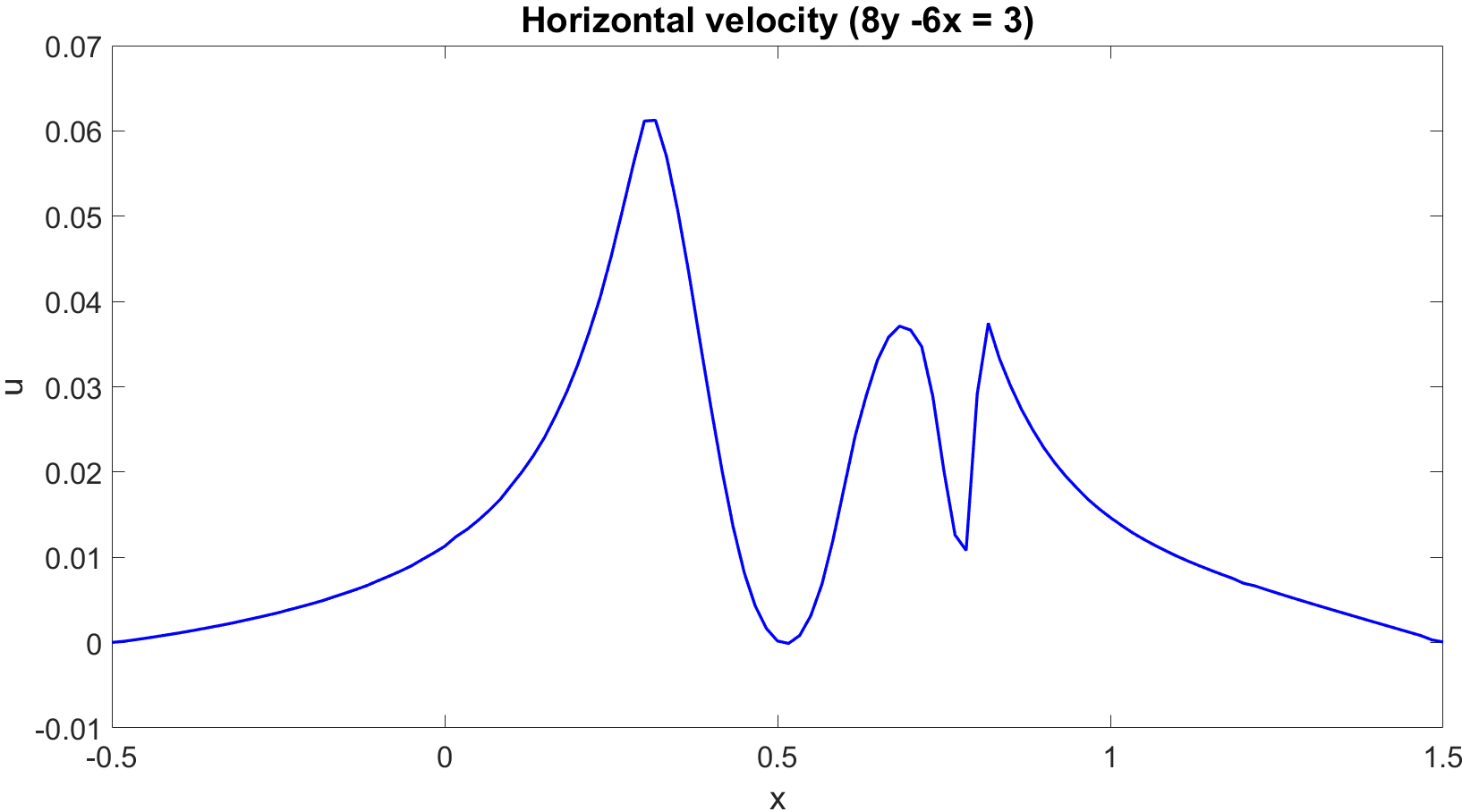}
	\end{subfigure}	
	\begin{subfigure}{0.45\textwidth}
		\centering
		\includegraphics[width=\textwidth]{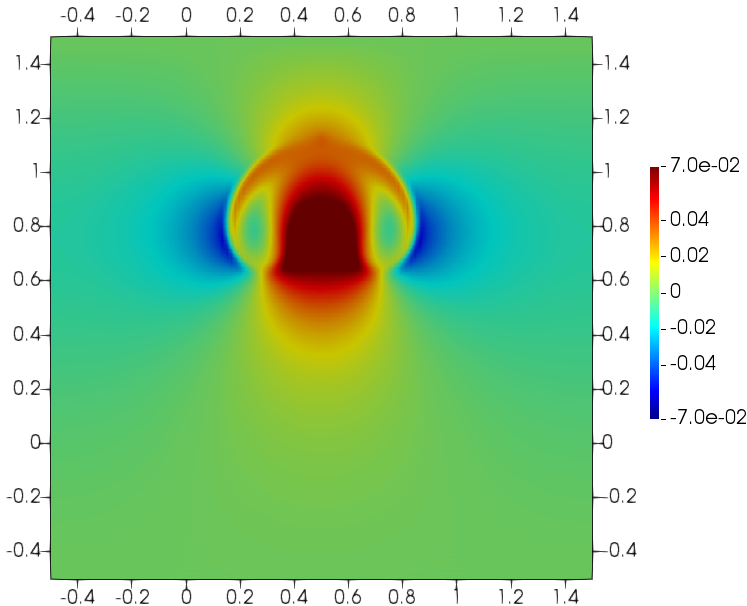}
	\end{subfigure}\hspace*{\fill}
	\begin{subfigure}{0.45\textwidth}
		\centering
		\includegraphics[width=\textwidth]{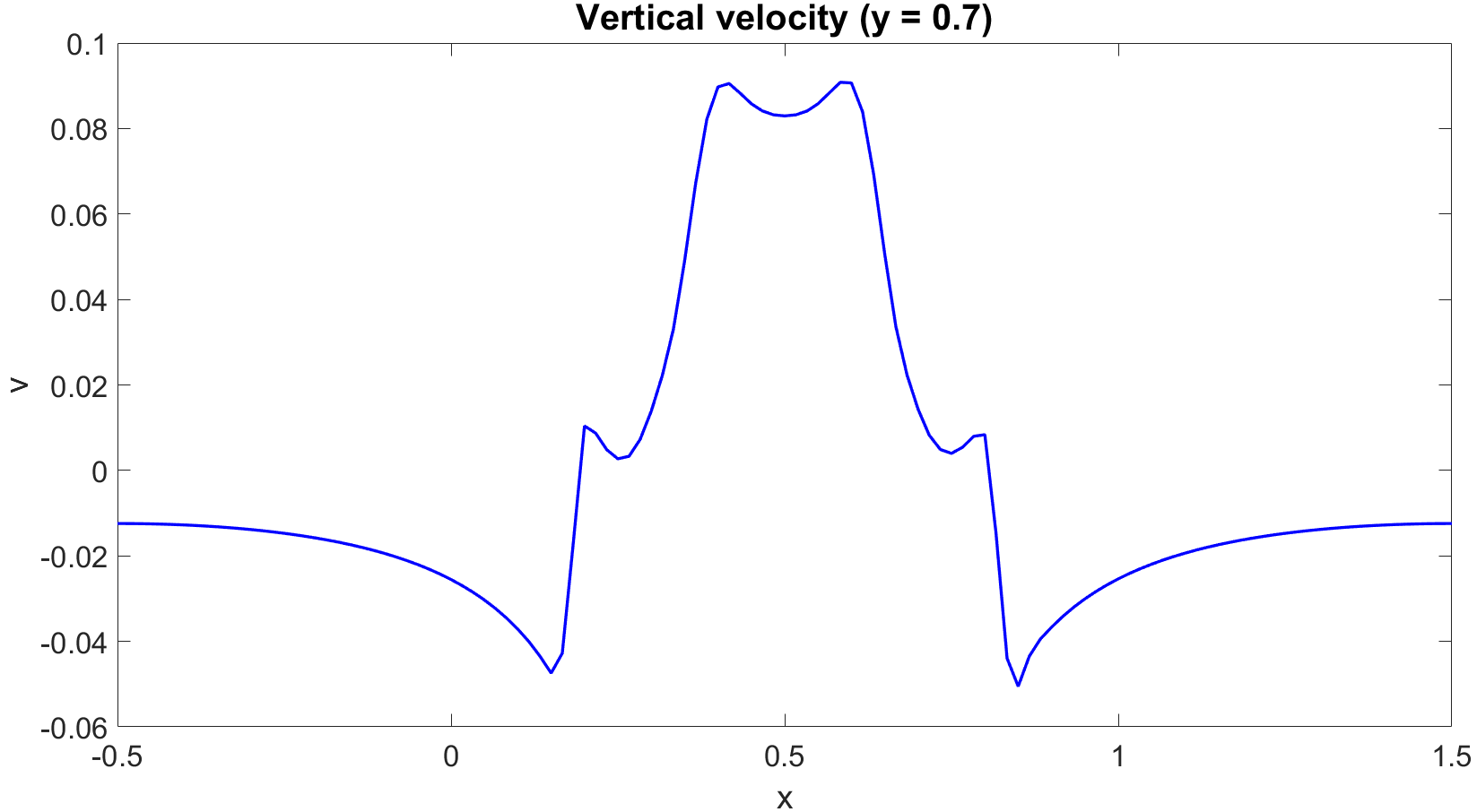}
	\end{subfigure}	
	\caption{Warm bubble test case, results at \(t = \SI{20}{\second}\). From bottom to top: temperature, horizontal velocity and vertical velocity.}
	\label{fig:Warm_bubble_t20}
\end{figure}

In order to test the method in presence of heat conduction, we now consider for an ideal gas the test case of a rising warm bubble proposed in \cite{busto:2020}. The domain is the square box \(\Omega = \left(-0.5, 1.5\right) \times \left(-0.5, 1.5\right)\) with periodic boundary conditions on the lateral boundaries and wall boundary conditions on the top and on the bottom of the domain. The initial temperature corresponds to a truncated Gaussian profile

\begin{equation}\label{eq:warm_bubble_T0_ideal}
T(\mathbf{x},0) = \begin{cases}
386.48 \qquad &\text{if } \tilde r > r_0 \\
\frac{\tilde p_0}{\tilde R_{g,air} \left(1 - 0.1e^{\frac{\tilde r^2}{\sigma^2}}\right)} \qquad &\text{if } \tilde r \le r_0,
\end{cases}
\end{equation}
where \(\tilde r^2 = \left(x - x_0\right)^2 + \left(y - y_0\right)^2\) is the distance from the center with coordinates \(x_0 = 0.5\) and \(y_0 = 0.35\), \(r_0 = 0.25\) is the radius and \(\sigma = 2\). In this Section, we consider unitary reference values for density, pressure and temperature and, therefore, we set \(\tilde p_0 = 10^5\) and \(\tilde R_{g,air} = 287\). Moreover, following \cite{busto:2020}, we consider:

\[Re = 804.9 \qquad Pr = 0.71 \qquad Fr \approx 0.004 \qquad M \approx 0.01.\]
The grid is composed by 120 elements along each direction and the time step is such that the maximum Courant number \(C \approx 118\) and the maximum value of advective Courant number \(C_u\) is around \(0.03\). Figures \ref{fig:Warm_bubble_t10}, \ref{fig:Warm_bubble_t15} and \ref{fig:Warm_bubble_t20} show the results at \(t = 10, 15, 20 \hspace{0.05cm} \SI{}{\second}\) both in terms of contours and plots along the same specific sections along \(x-\)axis chosen in \cite{busto:2020}. All the results are in good agreement with the reference ones and we are able to recover the development of the expected Kelvin-Helmholtz instability. 

\begin{figure}[H]
	\begin{subfigure}{0.5\textwidth}
		\centering
		\includegraphics[width=0.9\textwidth]{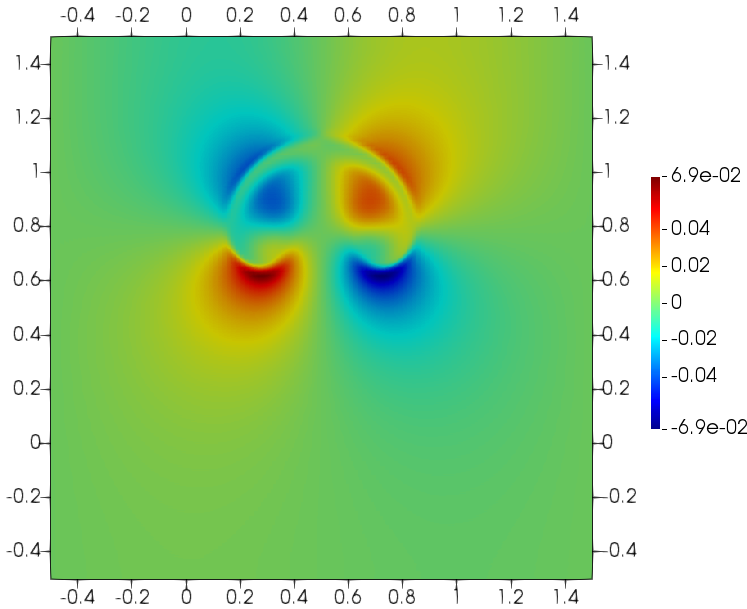} a)
	\end{subfigure}\hspace*{\fill}
	\begin{subfigure}{0.5\textwidth}
		\centering
		\includegraphics[width=0.9\textwidth]{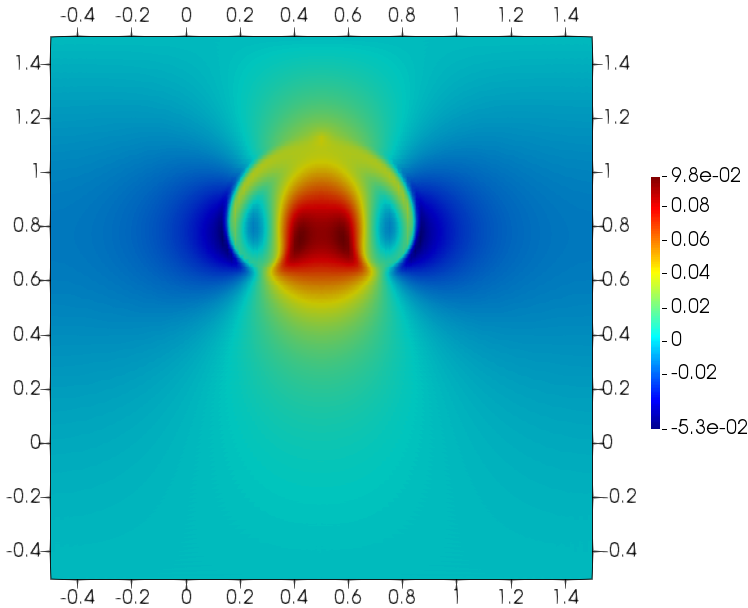} b)
	\end{subfigure}
	\begin{subfigure}{0.5\textwidth}
		\centering
		\includegraphics[width=0.9\textwidth]{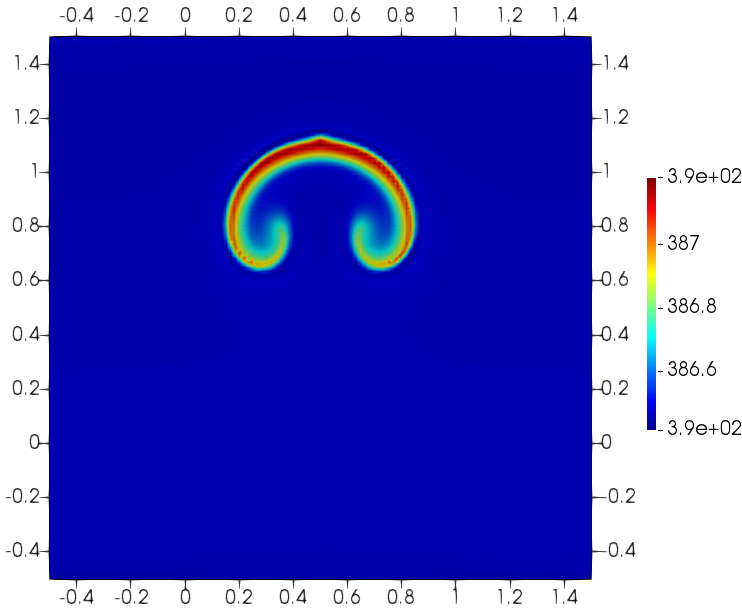} c)
	\end{subfigure}
	\caption{Warm bubble test case for \(N_{2}O\) with Peng-Robinson EOS, results at \(t = \SI{20}{\second}\), a) horizontal velocity, b) vertical velocity, c) temperature.}
	\label{fig:Warm_bubble_N2O_t20_PG}
\end{figure}

\begin{figure}[H]
	\centering
	\includegraphics[width=0.9\textwidth]{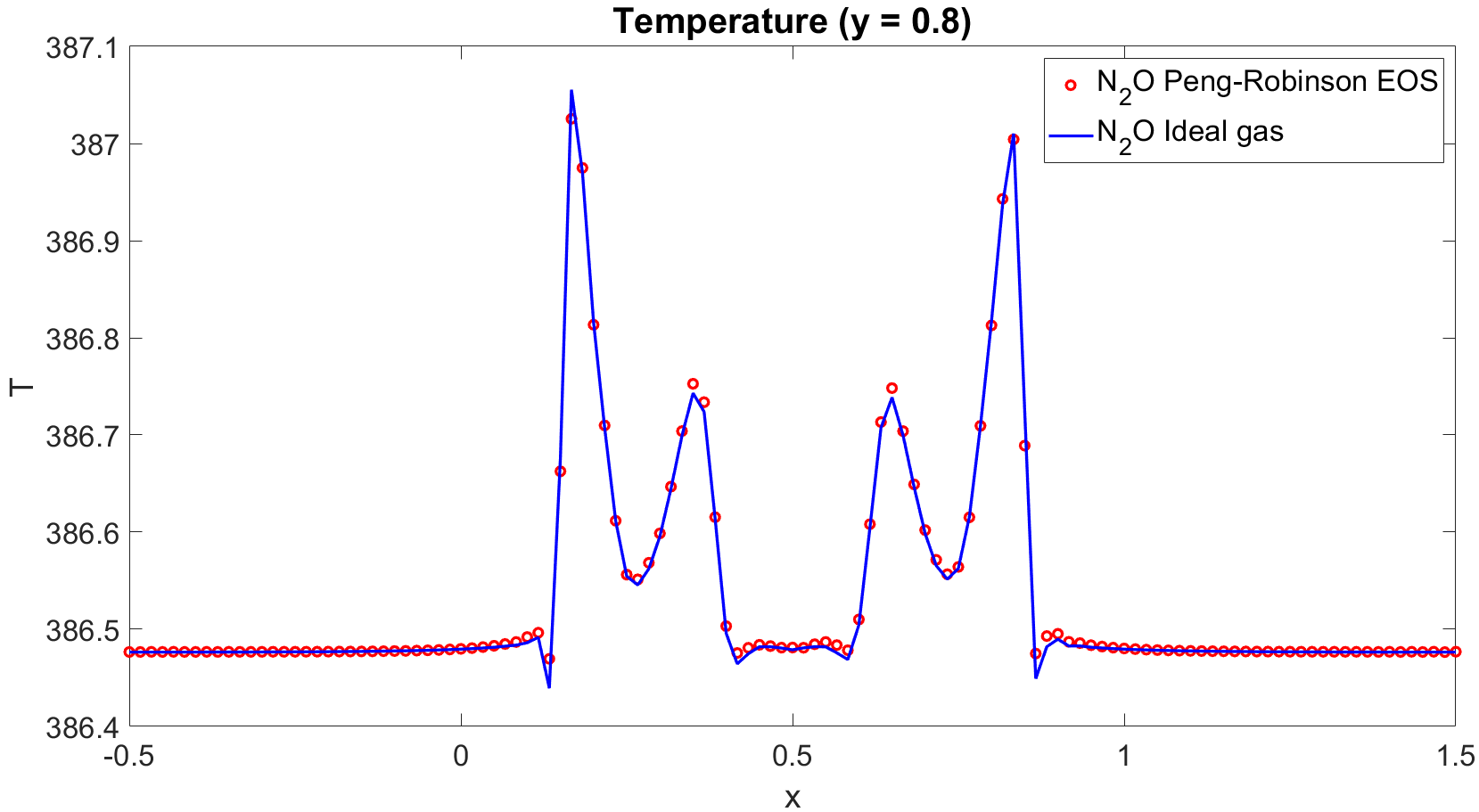}
	\caption{Warm bubble test case for \(N_{2}O\) with Peng-Robinson EOS, temperature profile for \(y = 0.8\) at \(t = \SI{20}{\second}\).}
	\label{fig:Warm_bubble_N2O_temperature_t20_comparison_PG}
\end{figure}

The same test is repeated using data for nitrous oxide \(\left(N_{2}O\right)\) from \cite{lias:2010}, which we report here for the convenience of the reader. At temperature of \(\SI{386.48}{\kelvin}\) and pressure of \(10^5 \hspace{0.05cm}\SI{}{\pascal}\), \(\mu = 1.8884 \cdot 10^{-5} \hspace{0.05cm} \SI{}{\pascal\cdot\second}\) and \(\kappa = 2.4855 \cdot 10^{-2} \hspace{0.05cm} \SI{}{\watt\per\meter\per\kelvin}\), so as to obtain

\[Re \approx 716.1 \qquad Pr \approx 0.73.\]
We consider the Peng-Robinson EOS, for which the expressions of \(\tilde a(T)\) and \(\tilde b\) are the following \cite{fernandez:2009}:

\begin{equation}
\begin{cases}
\tilde a(T) &= 0.45724\frac{\tilde R_{g}^2 \tilde T_{c}^2}{\tilde p_{c}}\alpha(T)^2 \\
\alpha(T) &= 1 + \Gamma\left(1 - \sqrt{\frac{T}{\tilde T_c}}\right) \\
\Gamma &= 0.37464 + 1.54226\omega - 0.26992\omega^2 \\
\tilde b &= 0.0778\frac{\tilde R_{g} \tilde T_{c}}{\tilde p_{c}},
\end{cases}
\end{equation}
where \(\tilde T_{c}\) denotes the non-dimensional critical temperature, \(\tilde p_{c}\) the non-dimensional critical pressure and \(\omega\) the acentric factor. For what concerns \(N_{2}O\), we find from \cite{lias:2010} \(\tilde T_{c} = 309.52, \tilde p_{c} = 7.2450 \cdot 10^{6}\) and \(\omega = 0.1613\). Finally, the function \(\tilde c_{v}(T)\) is computed using the following polynomial from \cite{lias:2010}:

\begin{eqnarray}\label{eq:cv_N2O}
\tilde c_{v}\left(T\right) &=& \frac{1}{T}\left[\left[A \frac{T}{1000} + \frac{1}{2}B \left(\frac{T}{1000}\right)^2 + \frac{1}{3}C\left(\frac{T}{1000}\right)^3 + \right.\right. \nonumber \\
&& \hspace{0.75cm} \left.\left. \frac{1}{4}D\left(\frac{T}{1000}\right)^4 - E\frac{1000}{T}\right]\frac{10^6}{M_{w}} - \tilde R_{g,N_{2}O}T\right],
\end{eqnarray}
with \(\tilde R_{g,N_{2}O} = 188.91\), \(M_{w} = 44.0128\) and \(A,B,C,D,E\) denoting suitable coefficients whose values are reported in Table \ref{tab:cv_N2O}. It is worth to mention once more that \(\tilde c_{v}(T)\) is not a proper specific heat at constant volume, but it denotes the non-dimensional counterpart of \(\frac{e^{\#}(T)}{T}\) from \eqref{eq:caloric_cubic_T_dip}, as shown in \eqref{eq:caloric_cubic_T_dip_adim}. The fluid is initialized with the same temperature and the same pressure as the ideal gas. The same mesh and the time step of the previous case are used, yielding to \(C \approx 92\) and \(C_{u} \approx 0.03\). Figure \ref{fig:Warm_bubble_N2O_t20_PG} shows the temperature, the horizontal and the vertical velocity at \(t = \SI{20}{\second}\). One can easily notice that a good qualitative agreement compared with the results in Figure \ref{fig:Warm_bubble_t20} is obtained. For a more quantitative point of view, since an explicit solution cannot be computed in view of the very large acoustic Courant number and the compressibility factor is \(z \approx 0.997\), a simulation with the ideal gas law \eqref{eq:EOS} is performed, using \(\gamma = 1.2879\), which corresponds to \(\frac{\tilde c_{v}\left(386.48\right)}{\tilde R_{g}} + 1\), so that the internal energy of the ideal gas at \(T = 386.48\) is the same as in the case \(e^{\#}\left(386.48\right)\). The temperature profile at \(y = 0.8\) shown in Figure \ref{fig:Warm_bubble_N2O_temperature_t20_comparison_PG} confirms the good quality of the solution, with only slight differences due to the different equations of state.  

\begin{table}[h!]
	\begin{center}
		\begin{tabular}{|c|c|}
			\hline
			\(A\) & \(27.67988\)   \\
			\hline
			\(B\) & \(51.14898\) \\ 
			\hline
			\(C\) & \(-30.64544\) \\
			\hline
			\(D\) & \(6.847911\) \\
			\hline
			\(E\) & \(-0.157906\) \\
			\hline
		\end{tabular}
	\end{center}
	\caption{Values for polynomial \eqref{eq:cv_N2O}}
	\label{tab:cv_N2O}
\end{table}

\begin{figure}[H]
	\begin{subfigure}{0.5\textwidth}
		\centering
		\includegraphics[width=0.9\textwidth]{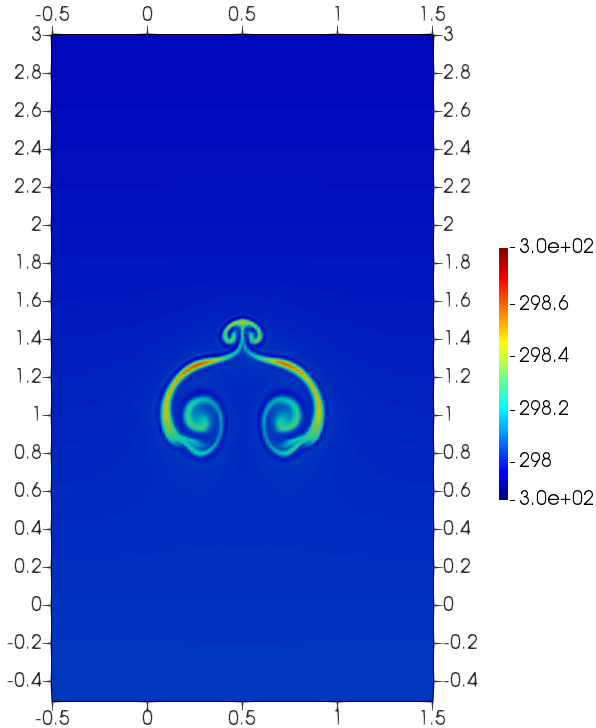} a)
	\end{subfigure}\hspace*{\fill}
	\begin{subfigure}{0.5\textwidth}
		\centering
		\includegraphics[width=0.9\textwidth]{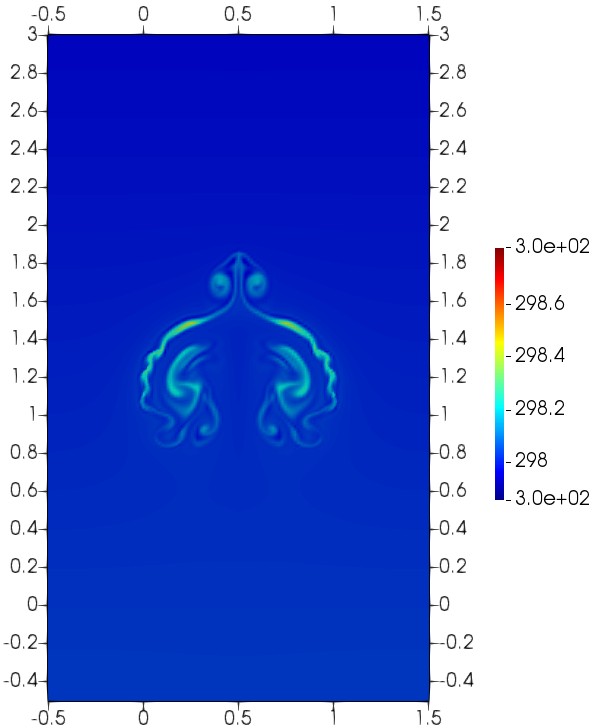} b)
	\end{subfigure}
	\caption{Warm bubble test case for \(N_{2}O\) with Peng-Robinson EOS, a) temperature at \(t = \SI{15}{\second}\), b) temperature at \(t = \SI{20}{\second}\).}
	\label{fig:Warm_bubble_N2O_sat_PG}
\end{figure}

In order to consider also non-ideal effects, we focus on conditions closer to the vapor-liquid phase equilibrium conditions, which are definitely more challenging. More in detail, we set \(p = 4 \cdot 10^{6}\) and we consider the following temperature profile

\begin{equation}
T(\mathbf{x},0) = \begin{cases}
298 \qquad &\text{if } \tilde r > r_0 \\
\frac{\tilde p_{0}}{\tilde R_{g,air}\left(1 - 0.1e^{\frac{\tilde r^2}{\sigma^2}}\right)} - 88.48 \qquad &\text{if } \tilde r \le r_0,
\end{cases}
\end{equation}
which corresponds to a translation with respect to \eqref{eq:warm_bubble_T0_ideal}, yielding \(z \approx 0.72\). We would like to point out that these conditions of pressure and temperature are close to the vapour-liquid phase transition curve of \(N_{2}O\). The maximum acoustic Courant is \(C \approx 74.5\), whereas the maximum advective Courant number is \(C_{u} \approx 0.06\). Figure \ref{fig:Warm_bubble_N2O_sat_PG} shows the contour plots of the temperature at \(t = \SI{15}{\second}\) and \(t = \SI{20}{\second}\). For these conditions of temperature and pressure, we obtain from \cite{lias:2010} \(\mu = 1.6680 \cdot 10^{-5} \hspace{0.05cm} \SI{}{\pascal\cdot\second}\), \(\kappa = 2.1201 \cdot 10^{-2} \hspace{0.05cm} \SI{}{\watt\per\meter\per\kelvin}\) and \(c_{p} = 1.5150 \cdot 10^{3} \hspace{0.05cm} \SI{}{\joule\per\kilogram\per\kelvin}\), so as to obtain

\[Re \approx 810.7 \quad Pr \approx 1.19\]
One can easily notice the full development of the Kelvin-Helmholtz instability with the formation of secondary vortices. The bubble reaches a higher altitude with respect to the previous case.

\begin{figure}[H]
	\begin{subfigure}{0.5\textwidth}
		\centering
		\includegraphics[width=0.9\textwidth]{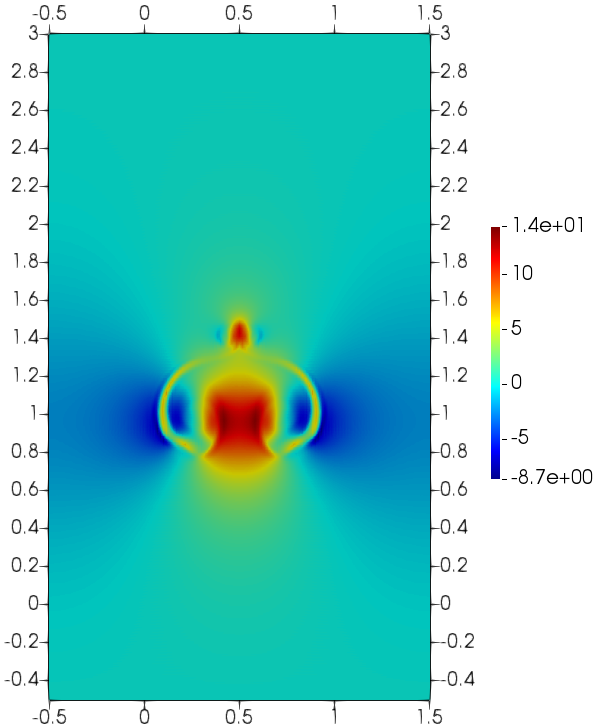} a)
	\end{subfigure}\hspace*{\fill}
	\begin{subfigure}{0.5\textwidth}
		\centering
		\includegraphics[width=0.9\textwidth]{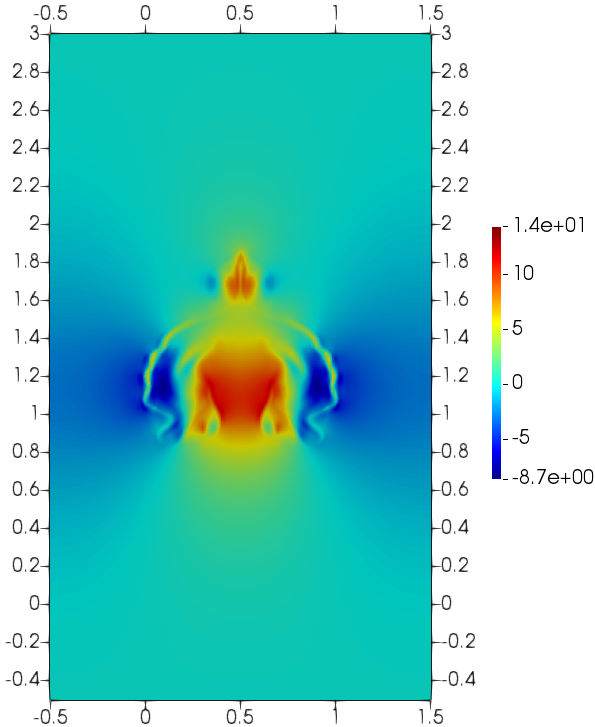} b)
	\end{subfigure}
	\caption{Warm bubble test case for \(N_{2}O\) with Peng-Robinson EOS, a) vertical velocity at \(t = \SI{15}{\second}\), b) vertical velocity at \(t = \SI{20}{\second}\).}
	\label{fig:Warm_bubble_N2O_sat_PG_v}
\end{figure}

\begin{figure}[H]
	\begin{subfigure}{0.5\textwidth}
		\centering
		\includegraphics[width=0.9\textwidth]{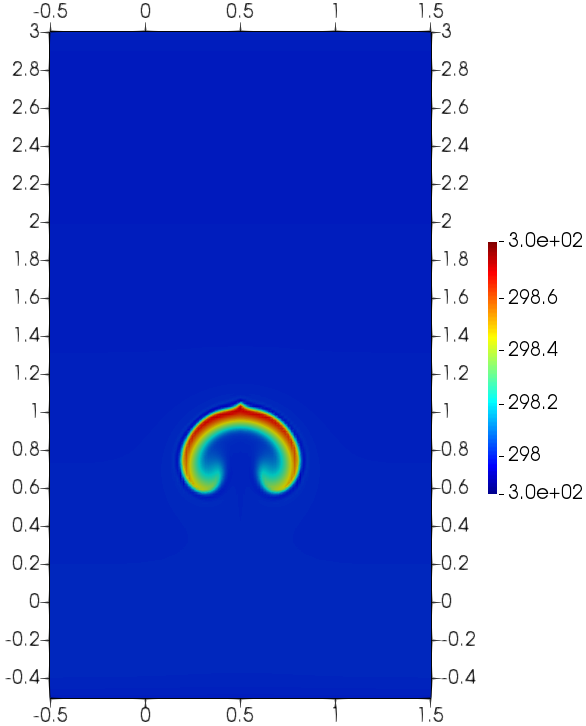} a)
	\end{subfigure}\hspace*{\fill}
	\begin{subfigure}{0.5\textwidth}
		\centering
		\includegraphics[width=0.9\textwidth]{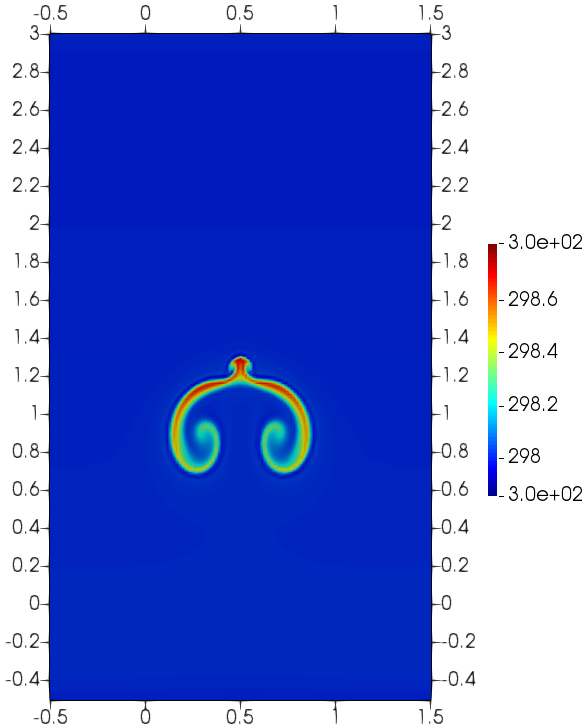} b)
	\end{subfigure}
	\caption{Warm bubble test case for \(N_{2}O\) with SG-EOS, a) temperature at \(t = \SI{15}{\second}\), b) temperature at \(t = \SI{20}{\second}\).}
	\label{fig:Warm_bubble_N2O_sat_SG}
\end{figure}

\begin{figure}[H]
	\begin{subfigure}{0.5\textwidth}
		\centering
		\includegraphics[width=0.9\textwidth]{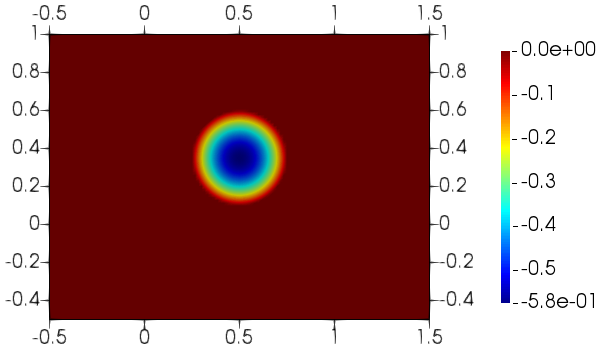} a)
	\end{subfigure}\hspace*{\fill}
	\begin{subfigure}{0.5\textwidth}
		\centering
		\includegraphics[width=0.9\textwidth]{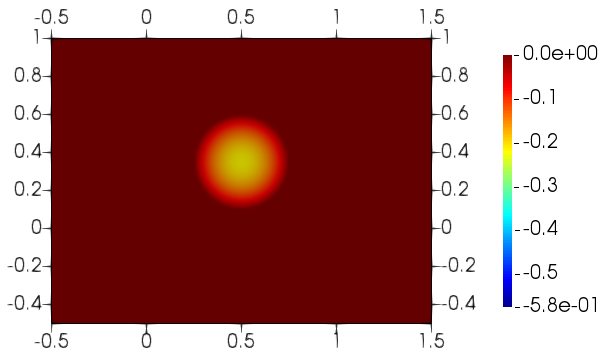} b)
	\end{subfigure}
	\caption{Warm bubble test case for \(N_{2}O\), density deviation from background state at \(t = 0\), a) Peng-Robinson EOS, b) SG-EOS.}
	\label{fig:Warm_bubble_N2O_rho_prime}
\end{figure}

Finally, we consider the SG-EOS with \(\gamma = 1.0936, \tilde c_{v} = 1453.91\) and \(\tilde \pi = \tilde q = 0\). The values for \(\gamma\) and \(\tilde c_{v}\) are computed using the procedure described in \cite{gandolfi:2019}. The maximum acoustic Courant is \(C \approx 67\), whereas the maximum advective Courant number is \(C_{u} \approx 0.04\). Figure \ref{fig:Warm_bubble_N2O_sat_SG} shows the contour plots of the temperature at \(t = \SI{15}{\second}\) and \(t = \SI{20}{\second}\). The different behaviour between the two equations of state can be readily explained since, in the case of Peng-Robinson EOS, the difference between the density of the bubble and the background density is bigger with respect to SG-EOS, as evident from Figure \ref{fig:Warm_bubble_N2O_rho_prime} and, therefore, a bigger upward buoyant force is exerted on the bubble. For this reason, in the simulation with Peng-Robinson EOS reaches a higher level compared to that in  the SG-EOS simulation.

\begin{figure}[H]
	\begin{subfigure}{0.5\textwidth}
		\centering
		\includegraphics[width=0.9\textwidth]{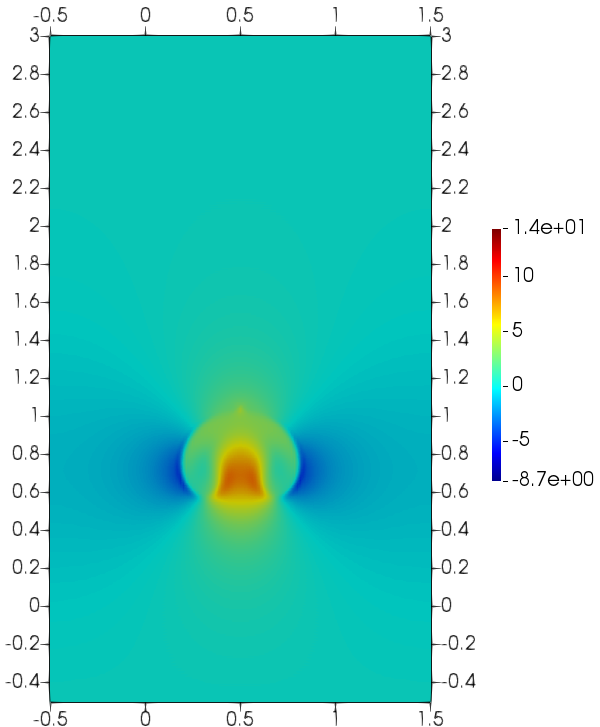} a)
	\end{subfigure}\hspace*{\fill}
	\begin{subfigure}{0.5\textwidth}
		\centering
		\includegraphics[width=0.9\textwidth]{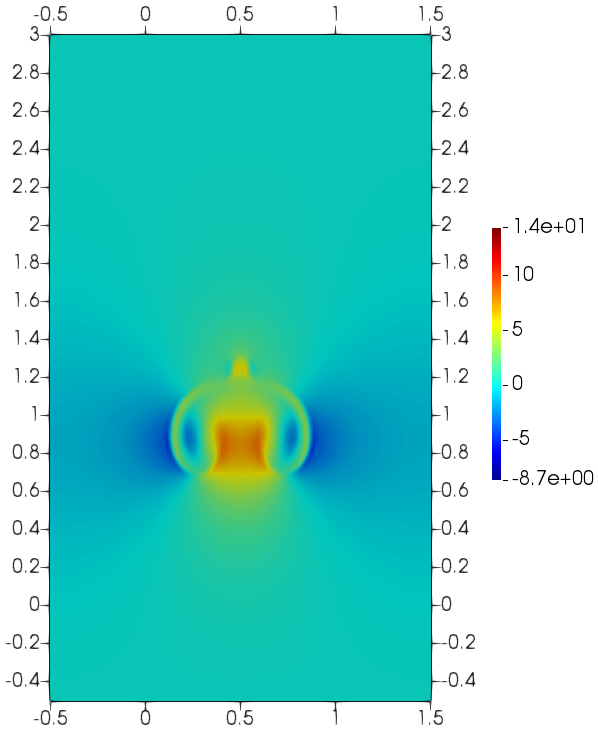} b)
	\end{subfigure}
	\caption{Warm bubble test case for \(N_{2}O\) with SG-EOS, a) vertical velocity at \(t = \SI{15}{\second}\), b) vertical velocity at \(t = \SI{20}{\second}\).}
	\label{fig:Warm_bubble_N2O_sat_SG_v}
\end{figure}

\section{Conclusions and future perspectives}
\label{sec:conclu}

We have proposed an efficient, $h-$adaptive IMEX-DG solver for the compressible Navier-Stokes equations with non-ideal EOS. The solver combines ideas from the discretization approaches in \cite{casulli:1984, dumbser:2016, giraldo:2013, kuehnlein:2019} and proposes an improvement in the choice of the free parameter employed by the explicit part of the IMEX scheme described in \cite{giraldo:2013}. The resulting method achieves full second order accuracy also including viscous terms and implements an $h-$adaptive approach in the framework of the numerical library \textit{deal.II}. A number of physically based adaptation criteria have been proposed, which allow to exploit the full efficiency of the adaptive technique also for non-ideal gas simulations. A number of numerical experiments validate the proposed method and show its potential for low Mach number problems. In future work, we plan to extend the scheme to multiphase flows and to demonstrate its potential for application to atmospheric flows.

\section*{Acknowledgements}
We thank M. Tavelli for providing the original data of the cavity flow simulation discussed in Section \ref{sec:tests}. We also gratefully acknowledge several useful discussions with A. Della Rocca on numerical methods related to those presented here. The simulations have been partly run at CINECA thanks to the computational resources made available through the SIDICoNS - HP10CLPLXI ISCRA-C project. 

\appendix
\section{Stability and monotonicity of the explicit time discretization}
\label{sec:mono_analysis}
In this Appendix, we study the stability and monotonicity of the explicit part of the IMEX scheme applied in the paper. We recall that the Butcher tableaux for the explicit part of the method is given by

\begin{figure}[H]
	\includegraphics[width=0.45\textwidth]{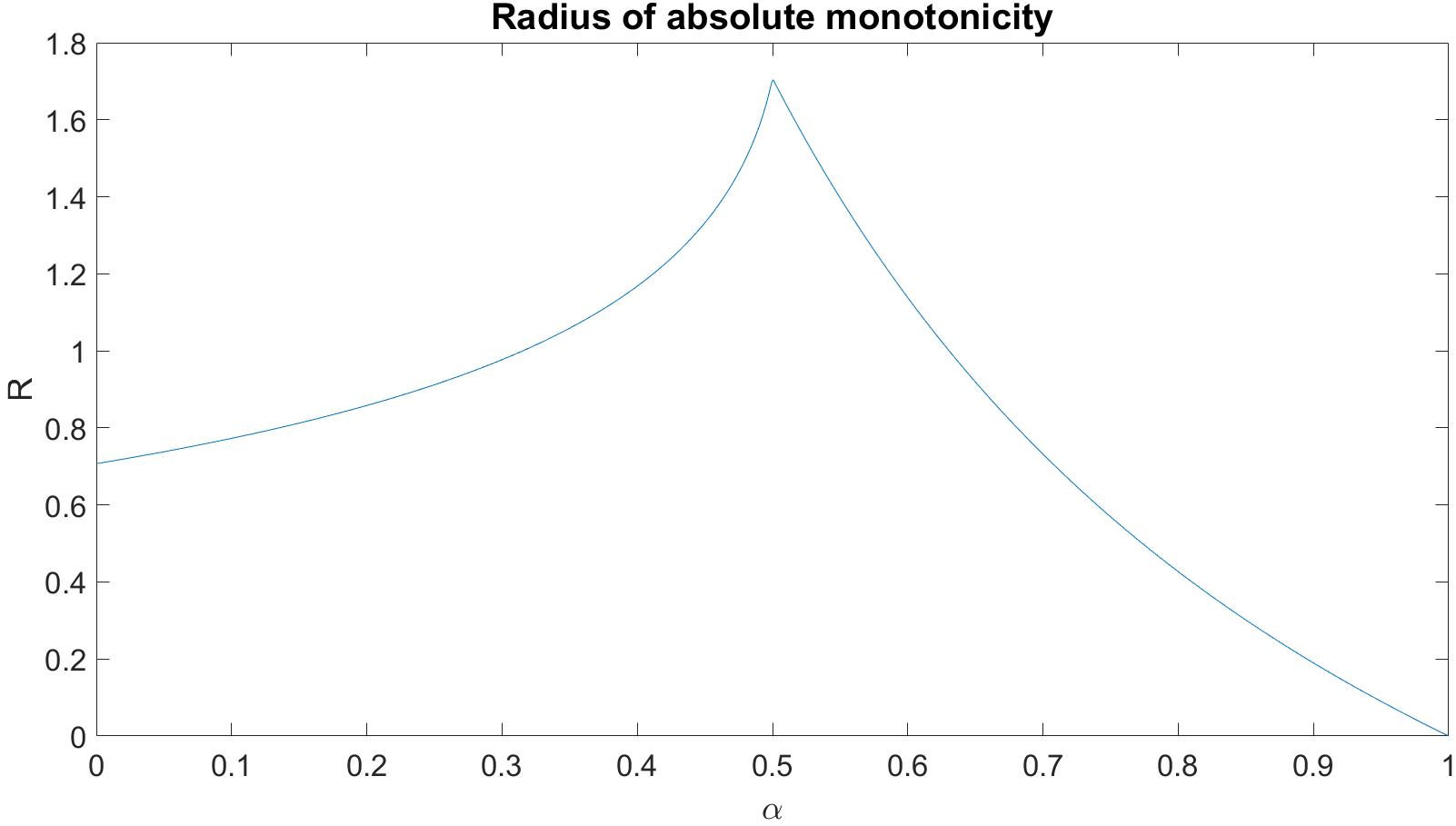} a)
	\includegraphics[width=0.45\textwidth]{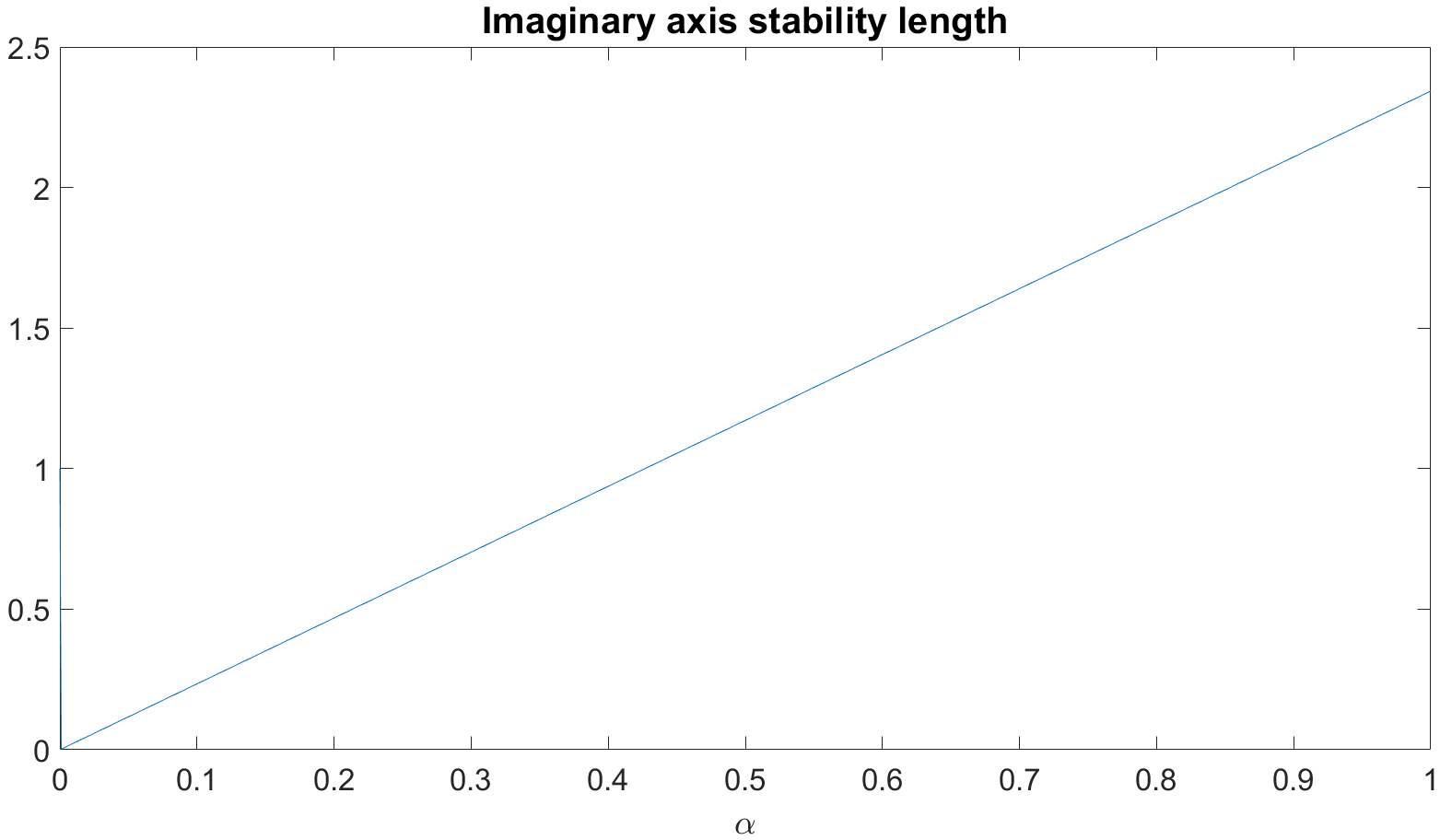} b)
	\caption{Analysis of the explicit part of IMEX scheme: a) Radius of absolute monotonicity as function of \(\alpha\), b) Size of stability region along the imaginary axis as \(\alpha\) varies.}
	\label{fig:Monotonicity}
\end{figure}

\begin{table}[h!]
	\begin{center}
		\begin{tabular}{c|ccc}
			0 & 0 & &  \\
			$\gamma$ & $\gamma$ & 0 & \\
			1 & $ 1 - \alpha$ & $\alpha$ & 0 \\
			\hline
			& $\frac{1}{2}-\frac{\gamma}4$ & $ \frac{1}{2}-\frac{\gamma}4$ & $\frac{\gamma}2$
		\end{tabular}
	\end{center}
\end{table}
In \cite{giraldo:2013}, the choice $ \alpha = \frac{7 - 2\gamma}{6} $ was made to maximize the stability region of the resulting scheme, but this coefficient is indeed a free parameter
and can also be chosen in different ways, as long as stability is not compromised. In order to identify possible alternative choices, we perform an analysis using the concepts introduced in \cite{kraaijevanger:1991}, \cite{higueras:2004}, \cite{ferracina:2004}
(see also the review in \cite{gottlieb:2001}). A similar analysis  for the implicit part of the IMEX scheme was carried out in \cite{bonaventura:2017}, to which we refer for a summary of the related theoretical results. We then define

\begin{equation*}
A = \begin{bmatrix}
0 & 0 & 0 \\
\gamma & 0 & 0 \\
1-\alpha & \alpha & 0
\end{bmatrix} \qquad
b^T = \begin{bmatrix}
\frac{1}{2}-\frac{\gamma}4 & \frac{1}{2}-\frac{\gamma}4 & \frac{\gamma}2 
\end{bmatrix}
\end{equation*}
with \(\gamma = 2 - \sqrt{2}\). We define for \(\xi \in \mathbb{R}\) the quantities

\begin{align}
A(\xi) &= A\left(I - \xi A\right)^{-1} \qquad b^{T}(\xi) = b^{T}\left(I - \xi A\right)^{-1} \nonumber \\
e(\xi) &= \left(I - \xi A\right)^{-1}e \qquad\hspace{0.27cm} \varphi(\xi) = 1 + \xi b^{T}\left(I - \xi A\right)^{-1}e
\end{align}
where \(I\) is the \(3\times 3\) identity matrix and \(e\) is a vector whose all components are equal to \(1\). Therefore, for the specific scheme, we obtain

\begin{equation*}
A(\xi) = \begin{bmatrix}
0 & 0 & 0 \\
\gamma & 0 & 0 \\
1 + \alpha\left(\gamma\xi - 1\right) & \alpha & 0
\end{bmatrix} 
\end{equation*}

\begin{equation*}
b^T(\xi) = \begin{bmatrix}
\frac{1}{4}\left[2+\gamma\left(-1 + \xi\left(4-\gamma+2\alpha\left(\gamma\xi - 1\right)\right)\right)\right] & \\ 
\frac{1}{4}\left[2+\gamma\left(2\alpha\xi - 1\right)\right] & \\
\frac{\gamma}2 
\end{bmatrix}
\end{equation*}

\begin{equation*}
e(\xi) = \begin{bmatrix}
1 & \\ 
1 + \gamma\xi & \\
1 + \xi + \alpha\gamma\xi^2
\end{bmatrix}
\end{equation*}

\begin{equation*}
\varphi(\xi) = 1 + \xi + \frac{\xi^2}{2} + \left(3 - 2\sqrt{2}\right)\alpha\xi^3.
\end{equation*}
A method with tableaux \(\left(A,b^{T}\right)\) is absolutely monotone at \(\xi \in \mathbb{R}\) if \(A(\xi) \ge 0\), \(b^{T}(\xi) \ge 0\), \(e(\xi) \ge 0\) and \(\varphi(\xi)\ge 0\) elementwise; moreover the radius of absolute monotonicity is defined for all \(\xi\) in \(-r \le \xi \le 0\) as 

\[R(a,b) = \sup\left[r|r\ge 0, A(\xi) \ge 0,b^{T}(\xi) \ge 0,e(\xi) \ge 0, \varphi(\xi)\ge 0\right].\] 
Figure \ref{fig:Monotonicity} shows the behaviour of the radius of absolute monotonicity as \(\alpha\) varies, along with the behaviour of the stability region along the imaginary axis. As already mentioned before, \(\alpha = \frac{7 - 2\gamma }{6}\) was chosen originally to maximize the stability region, but in this case \(R = \frac{2\sqrt{2} - 3}{2 + \sqrt{2}} \approx 0.05\), so that the region of absolute monotonicity is quite small. It can be shown that the region of absolute stability is given by

$$ S=\left\{z \in \mathbb{C} : \left|1 + z + \alpha\gamma z^2\right| < 1 \right\}. $$
The alternative value \(\alpha = 0.5\) maximizes the region of absolute monotonicity without compromising too much the stability. The impact of this alternative choice on numerical results is discussed in Section \ref{sec:tests}.

\section{Eigenvalues of 1D Euler equations}
\label{sec:eigenvalues}

In this Appendix we compute the eigenvalues for the Euler equations in non-dimensional form for a general equation of state. For the sake of simplicity, we focus on 1D case and so the equations can be written as follows

\begin{eqnarray}
\frac{\partial \rho}{\partial t}  
+\frac{\partial}{\partial x}\left(\rho u\right)&=& 0   \nonumber \\
\frac{\partial \rho u}{\partial t}  
+\frac{\partial}{\partial x}\left(\rho u^2\right) 
+ \frac{1}{M^2}\frac{\partial p}{\partial x} &=& 0 \\
\frac{\partial \rho E}{\partial t}  
+\frac{\partial}{\partial x}\left[\left(\rho E + p\right)u \right] &=& 0. \nonumber
\end{eqnarray}
This is equivalent to

\begin{eqnarray}
\frac{\partial \rho}{\partial t} + u\frac{\partial \rho}{\partial x} + \rho \frac{\partial u}{\partial x} &&= 0 \nonumber \\
\frac{\partial \rho}{\partial t}u + \frac{\partial u}{\partial t}\rho  
+ u^2\frac{\partial \rho}{\partial x} + 2\rho u\frac{\partial u}{\partial x} 
+ \frac{1}{M^2}\frac{\partial p}{\partial x} &&= 0 \\
\frac{\partial \rho}{\partial t}E + \frac{\partial E}{\partial t}\rho  
+ \left(\rho E + p\right)\frac{\partial u}{\partial x} + u\left(\frac{\partial \rho}{\partial x}E + \frac{\partial E}{\partial x}\rho + \frac{\partial p}{\partial x}\right) &&= 0. \nonumber
\end{eqnarray}
Thanks to the continuity equation and to the relation \(E = e + \frac{1}{2}M^2u^2\), we obtain 

\begin{eqnarray}
\frac{\partial \rho}{\partial t} + u\frac{\partial \rho}{\partial x} + \rho \frac{\partial u}{\partial x} &=& 0 \nonumber \\
\frac{\partial u}{\partial t} + u\frac{\partial u}{\partial x} 
+ \frac{1}{\rho M^2}\frac{\partial p}{\partial x} &=& 0 \\
\frac{\partial e}{\partial t} + \frac{p}{\rho}\frac{\partial u}{\partial x} + u\frac{\partial e}{\partial x} &=& 0. \nonumber
\end{eqnarray}
In general \(e = e(p, \rho)\), so that

\begin{equation}
\frac{\partial e}{\partial \square} = \frac{\partial e}{\partial \rho}\frac{\partial \rho}{\partial \square} + \frac{\partial e}{\partial p}\frac{\partial p}{\partial \square},
\end{equation}
and, hence, the system reduces to

\begin{eqnarray}
\frac{\partial \rho}{\partial t} + u\frac{\partial \rho}{\partial x} + \rho \frac{\partial u}{\partial x} &=& 0 \nonumber \\
\frac{\partial u}{\partial t} + u\frac{\partial u}{\partial x} 
+ \frac{1}{\rho M^2}\frac{\partial p}{\partial x} &=& 0 \\
\frac{\partial p}{\partial t} + \frac{\left(\frac{p}{\rho} - \rho\frac{\partial e}{\partial \rho}\right)}{\frac{\partial e}{\partial p}}\frac{\partial u}{\partial x} + u \frac{\partial p}{\partial x} &=& 0, \nonumber
\end{eqnarray}
which can be thought in the following vector form as:

\begin{equation}
\frac{\partial\mathbf{Q}}{\partial t} + \mathbf{A}\frac{\partial\mathbf{Q}}{\partial x} = \mathbf{0}
\end{equation}
with

\begin{equation}
\mathbf{Q} = \begin{bmatrix}
\rho \\
u \\
p
\end{bmatrix} 
\end{equation}
and

\begin{equation}\label{eq:matrix_Euler}
\mathbf{A} = \begin{bmatrix}
u & \rho & 0 \\
0 & u & \frac{1}{\rho M^2} \\
0 & \frac{\left(\frac{p}{\rho} - \rho\frac{\partial e}{\partial \rho}\right)}{\frac{\partial e}{\partial p}} & u
\end{bmatrix}.
\end{equation}
The eigenvalues of \eqref{eq:matrix_Euler} are \(u - \frac{1}{M}\frac{1}{\rho}\sqrt{\frac{p - \frac{\partial e}{\partial \rho}\rho^2}{\frac{\partial e}{\partial p}}}\), \(u\) and \(u + \frac{1}{M}\frac{1}{\rho}\sqrt{\frac{p - \frac{\partial e}{\partial \rho}\rho^2}{\frac{\partial e}{\partial p}}}\). The first law of thermodynamics, already recalled in Section \ref{sec:EOS}, provides us the following relation

\begin{equation}
Tds = de - \frac{p}{\rho^2}d\rho = \left(\frac{\partial e}{\partial \rho} - \frac{p}{\rho^2}\right)d\rho + \frac{\partial e}{\partial p}dp, 
\end{equation}
or, equivalently,

\begin{equation}
dp = \frac{\frac{p}{\rho^2} - \frac{\partial e}{\partial \rho}}{\frac{\partial e}{\partial p}}d\rho + \frac{T}{\frac{\partial e}{\partial p}}ds.
\end{equation}
Hence, following \cite{vidal:2001}, we have

\begin{equation}
c^2 = \frac{\partial p}{\partial \rho}\bigg\rvert_{s} = \frac{\frac{p}{\rho^2} - \frac{\partial e}{\partial \rho}}{\frac{\partial e}{\partial p}}
\end{equation}
and, therefore, the eigenvalues of \eqref{eq:matrix_Euler} are 
$$ u + \frac{c}{M} \ \ \ u \ \ \ u + \frac{c}{M}$$
also for a generic equation of state, and not only in the case of an ideal gas, as already discussed in \cite{munz:2003}. This justifies the definition \eqref{eq:Courant} also in case of non-ideal gases.

\bibliographystyle{plain}
\bibliography{DG_NS_COMP_IMEX}

\end{document}